\documentclass{article}
\usepackage{a4}

\input amssym.def
\input amssym

\font\zwm=msbm10 at 12pt
\font\zwe=eusm10 at 12pt

\font\teneusm=eusm10
\font\seveneusm=eusm7
\font\fiveeusm=eusm5
\newfam\eusmfam
\textfont\eusmfam=\teneusm
\scriptfont\eusmfam=\seveneusm
\scriptscriptfont\eusmfam=\fiveeusm
\def\script#1{{\fam\eusmfam\relax#1}}

\newtheorem{thm}{Theorem}[section]
\newtheorem{lem}[thm]{Lemma}

\newtheorem{prp}[thm]{Proposition}

\newtheorem{thmsub}{Theorem}[subsection]
\newtheorem{lemsub}[thmsub]{Lemma}
\newtheorem{prpsub}[thmsub]{Proposition}
\newtheorem{corsub}[thmsub]{Corollary}
\newtheorem{dfnsub}[thmsub]{Definition}
\newtheorem{ldcsub}[thmsub]{Computational lemma}

\newtheorem{thmsubf}{Th\'eor\`eme}[subsection]
\newtheorem{lemsubf}[thmsubf]{Lemme}
\newtheorem{prpsubf}[thmsubf]{Proposition}
\newtheorem{corsubf}[thmsubf]{Corollaire}
\newtheorem{dfnsubf}[thmsubf]{D\'efinition}

\def\dem{{\it Proof.\/ }} 
\def\rem{\refstepcounter{thm}
{\bf Remark \thethm\hskip1ex}}
\def\remsub{
\refstepcounter{thmsub}{\bf Remark \thethmsub}\hskip1ex}

\def\exasub{
\addtocounter{thmsub}{1}{\bf Example \thethmsub}\hskip1ex}

\def\questsub{
\refstepcounter{thmsub}{\bf Question \thethmsub}\hskip1ex}
\def\questf{
\refstepcounter{thmf}{\bf Question \thethmf}\hskip1ex}
\def\questsubf{
\refstepcounter{thmsubf}{\bf Question \thethmsubf}\hskip1ex}

\def\vskipa{\vskip\abovedisplayskip}
\def\vskipb{\vskip\belowdisplayskip}
\def\ie{{\it i.\,e.}\ } 
\def\vs{{\it vs.}\ } 
\def\af{{\it a fortiori}\/\ }

\def\toll#1{\vtop{\ialign {##\crcr \rightarrowfill \crcr 
\noalign{\kern -1pt \nointerlineskip \vskip2pt
} 
$\hfil \scriptstyle{\ #1\ } \hfil $\crcr }}}

\def\({{\fam0\rm (}}
\def\){\/{\fam0\rm )}} 
\def\]{\mathopen]}
\def\[{\mathclose[}
\mathchardef\hook="312C

\def\tol{\mathop{-\mkern-4mu\smash-\mkern-4mu\smash\longrightarrow}\limits}

\def\imp{\Rightarrow}
 
\def\ssi{\Leftrightarrow}
\def\Ssi{\Longleftrightarrow}
\def\ecks{\rule{2mm}{2mm}}
\def\bull{\leavevmode\kern .1ex\vrule height 1ex width .9ex depth
-.1ex \kern .8ex} 
\def\bloc{\bull}
\def\eck{\nolinebreak\hspace{\fill}\ecks}
\def\barre{\rule[-2.5pt]{0.8pt}{10pt}}
\def\lmes{\kern0.7pt\barre\kern0.7pt}
\def\rmes{\kern0.7pt\/\barre\kern0.7pt}
\def\mes#1{\mbox{$\lmes #1 \rmes$}}
\def\smes#1{|#1|}
\def\bi#1#2{\bigg(\kern-0.5pt{{#1}\atop{#2}}\kern-0.5pt\bigg)}
\def\bip#1#2{\left(\kern-3pt {{#1}\atop{#2}}\kern-2pt \right)}
\def\se{\subseteq}
\def\HM#1{\setbox0=\hbox{#1}\dimen0=\wd0 #1 \kern-\dimen0
\setbox0\hbox{\raise3pt\hbox{$^\frown$}}\advance\dimen0 by -\wd0
\box0\kern\dimen0}

\def\sgn{\mathop{\fam0 sgn}\nolimits}

\def\osc{\mathop{\fam0 osc}}
\def\im{\mathop{\fam0 im}\nolimits} 
\def\suml{\sum\limits}

\def\maxl{\max\limits} 
\def\P#1{\mathop{\Bbb P}\nolimits
\left[\kern0.5pt{#1}\kern0.5pt\right]}
 
\def\e{\mathop{\fam0 e}\nolimits}
\def\Id{\hbox{\fam0\rm Id}}
\def\ii{{\fam0 i}} 
\def\Alpha{{\fam0 A}}

\def\Zeta{{\fam0 Z}} 
\def\eps{\varepsilon}

\def\Z{\hbox{${\Bbb Z}$}}
\def\U{\hbox{${\Bbb S}$}}
\def\T{\hbox{${\Bbb T}$}}
\def\N{\hbox{${\Bbb N}$}}
\def\R{\hbox{${\Bbb R}$}}
\def\D{\hbox{${\Bbb D}$}}
\def\C{\hbox{${\Bbb C}$}}

\def\sZ{{\Bbb Z}} 
\def\sU{{\Bbb S}}
\def\sT{{\Bbb T}}
\def\sN{{\Bbb N}}

\def\ds{\displaystyle} \def\hoch#1{\mathop{#1}\limits}

\def\lst#1#2{#1_1,\dots,#1_{#2}}
\def\lstf#1#2{#1_{#2},#1_{#2+1},\dots}
\def\lstp#1#2#3{#1_{#2},\dots,#1_{#3}}
\def\sm#1#2{#1_1+\dots+#1_{#2}} \def\bip#1#2{\left(\kern-3pt
{{#1}\atop{#2}}\kern-2pt \right)}
\def\smp#1#2#3{#1_{#2}+\dots+#1_{#3}} 

\def\lap#1{$\ell_{#1}$-$(\!${\it ap}$)$}\def\lpap/{\lap{p}}
\def\lmap#1{$\ell_{#1}$-$(\!${\it map}$)$}\def\lpmap/{\lmap{p}}
\def\uap/{\mbox{$(\!${\it uap}$)$}}\def\ubp/{\mbox{$(\!${\it ubp}$)$}}
\def\as/{\mbox{a.s.}}\def\sbd/{\mbox{$(\!${\it sbd\/}$)$}}
\def\ubs/{\mbox{$(\!${\it ubs}$)$}}\def\umbs/{\mbox{$(\!${\it umbs}$)$}}
\def\umap/{\mbox{$(\!${\it umap}$)$}}
\def\ap/{\mbox{$(\!${\it ap}$)$}}\def\map/{\mbox{$(\!${\it map}$)$}}
\def\fdd/{\mbox{$(${\it fdd\/}$)$}}\def\umfdd/{$(${\it umfdd\/}$)$}

\def\SC#1{\script{C}_{#1}(\T)}
\def\SCE{\SC{E}}
\def\SL#1#2{{\fam0 L}^{#1}_{#2}(\T)}
\def\SLE#1{\SL{#1}{E}}
\def\SLP#1{\SL{p}{#1}}
\def\SLPE{\SLE{p}}

\def\DX/{$X\in\{\SC{},\SLP{}:1\le p<\infty\}$}
\def\PT#1{\mbox{$\script{P}_{#1}(\T)$}}
\def\PTE/{\PT{E}}
\def\J#1{\mbox{$(\script{J}_{#1})$}}
\def\I#1{\mbox{$(\script{I}_{#1})$}}
\def\UP/{\mbox{$(\script{U})$}}
\def\X#1{{\langle\zeta,#1\rangle}}\def\XE{{\X E}}
\def\EL#1{$\Lambda(#1)$}
\def\ER/{\mbox{Rosenthal}} 

\def\DE/{\mbox{$E\subseteq\Z$}}
\def\DEE/{\mbox{$E=\{n_k\}\subseteq\Z$}}
\def\DEEE/{\mbox{$E=\{n_k\}_{k\ge1}\subseteq\Z$}}
\def\DL/{\mbox{$\L/\subseteq\Z$}}
\def\DLL/{\mbox{$\L/=\{\lambda_k\}\subseteq\Z$}}
\def\DLLL/{\mbox{$\L/=\{\lambda_k\}_{k\ge1}\subseteq\Z$}}
\def\Ref#1{\hbox{$(\ref{#1})$}}

\def\DI/{$I\subseteq\N\times\N$}

\def\ud/{equidistributed}  
\def\wud/{weakly equidistributed}  

\def\pstar{^{\scriptscriptstyle(\kern-1pt\lower0.5pt
\hbox{$\scriptstyle *$}\kern-1pt)}}
\def\ppstar{^{\scriptscriptstyle(\kern-1pt*\kern-1pt)}}

\author{Stefan Neuwirth}
\title{Metric unconditionality and Fourier analysis}
\date{}

\begin{document}\parindent=0pt

\maketitle\vfill
{\bf R\'esum\'e } 
Nous \'etudions plusieurs propri\'et\'es fonctionnelles
d'inconditionnalit\'e iso\-m\'e\-trique et presqu'iso\-m\'etrique en les
exprimant \`a l'aide de multiplicateurs. Parmi ceux-ci, la notion la plus
g\'en\'erale est celle de ``propri\'et\'e d'approximation
inconditionnelle m\'etrique''. Nous la caract\'erisons parmi les
espaces de Banach de cotype fini par une propri\'et\'e simple
d'``inconditionnalit\'e par blocs''. En nous ramenant \`a des
multiplicateurs de Fourier, nous \'etudions cette propri\'et\'e dans
les sous-espaces des espaces de Banach de fonctions sur le cercle qui
sont engendr\'es par une suite de ca\-rac\-t\`eres $\e^{\ii nt}$. Nous
\'etudions aussi les suites basiques inconditionnelles iso\-m\'etriques
et presqu'isom\'etriques de caract\`eres, en particulier les ensembles
de Sidon de constante asymptotiquement $1$. Nous obtenons dans chaque
cas des propri\'et\'es combinatoires sur la suite. La propri\'et\'e
suivante des normes ${\fam0 L}^p$ est cruciale pour notre \'etude: si $p$ est
un entier pair,
$\int|f|^p=\int{|f^{p/2}|}^2=\sum|\widehat{f^{p/2}}(n)|{\vphantom{|f^{p/2}|}}^2$
est une
expression polyno\-miale en les coefficients de Fourier de $f$ et
$\bar f$. Nous proposons d'ailleurs une estimation pr\'ecise de la
constante de Sidon des ensembles \`a la Hadamard. 
\vskip\baselineskip

{\bf Zusammenfassung } Verschiedene funktionalanalytische isometrische
und fast iso\-me\-tri\-sche Unbedingtheitseigenschaften werden mittels
Multiplikatoren untersucht. Am allgemeinsten ist die metrische
unbedingte Approximationseigenschaft gefasst. Wir charakterisieren
diese f\"ur Banach\-r\"aume mit endlichem Kotyp durch eine einfache
``blockweise'' Unbedingtheit. Daraufhin betrachten wir genauer den
Fall von Funktionenr\"aumen auf dem Einheitskreis, die durch eine
Folge von Frequenzen $\e^{\ii nt}$ aufgespannt werden. Wir untersuchen
isometrisch und fast isometrisch un\-be\-ding\-te Basisfolgen von
Frequenzen, unter anderem Sidonmengen mit einer Konstante asymptotisch
zu $1$.  F\"ur jeden Fall erhalten wir kombinatorische Eigenschaften
der Folge. Die folgende Eigenschaft der ${\fam0 L}^p$ Normen ist
entscheidend f\"ur diese Arbeit: Ist $p$ eine gerade Zahl, so ist
$\int|f|^p=\int{|f^{p/2}|}^2=\sum|\widehat{f^{p/2}}(n)|
{\vphantom{|f^{p/2}|}}^2$ ein polynomialer Ausdruck der
Fourierkoeffizienten von $f$ und $\bar f$. Des weiteren erhalten wir
eine genaue Absch\"atzung der Sidonkonstante von Hadamardfolgen.
\vskip\baselineskip

{\bf Abstract }
We study several functional properties of isometric
and almost isometric unconditionality and state
them as a property of families of multipliers. The
most general such notion is that of ``metric unconditional
approximation property''. We characterize this ``\umap/'' by a
simple property of ``block unconditionality'' for spaces with
nontrivial cotype. We focus on subspaces of Banach spaces of functions
on the circle spanned by a sequence of characters $\e^{\ii nt}$. There
\umap/ may be stated in terms of Fourier multipliers.  We express
\umap/ as a simple combinatorial property of this sequence. We
obtain a corresponding result for isometric and almost isometric basic
sequences of characters. Our study uses the following crucial property
of the ${\fam0 L}^p$ norm for even $p$:
$\int|f|^p=\int{|f^{p/2}|}^2=\sum|\widehat{f^{p/2}}(n)|{\vphantom{|f^{p/2}|}}^2$
is a
polynomial expression in the Fourier coefficients of $f$ and $\bar f$.
As a byproduct, we get a sharp estimate of the Sidon constant of sets
{\it \`a la}\/ Hadamard. \pagebreak

\tableofcontents\vfill\pagebreak

\section{A general introduction in French}

\subsection{Position du probl\`eme}
Ce travail se situe au croisement de l'analyse fonctionnelle et de
l'analyse harmonique.  Nous allons donner des \'el\'ements de
r\'eponse \`a la question g\'en\'erale suivante.

\questf
Quelle est la validit\'e de la repr\'esentation
\begin{equation}\label{intro-these:q1}
f\sim\sum \varrho_q\e^{\ii\vartheta_q}\e_q
\end{equation}
de la fonction $f$
comme s\'erie de fr\'equences $\e_q$ d'intensit\'e $\varrho_q$ et de
phase $\vartheta_q$~?

Les r\'eponses seront donn\'es en termes de l'espace de fonctions
$X\ni f$ et du spectre $E\supseteq\{q:\varrho_q>0\}$.

Consid\'erons par exemple les deux questions classiques suivantes dans
le cadre des espaces de Banach homog\`enes de fonctions sur le tore
$\T$, des fr\'equences de Fourier $\e_q(t)=\e^{\ii qt}$ et des
coefficients de Fourier
$$\varrho_q\e^{\ii\vartheta_q}=\int\e_{-q}f=\widehat{f}(q).$$

\questsubf\label{intro:q1}
Est-ce que pour les fonctions $f\in X$ \`a spectre
dans $E$
$$\Bigl\| f-\sum_{|q|\le n}\varrho_q\e^{\ii\vartheta_q}\e_q
\Bigr\|_X\tol_{n\to\infty}0\hbox{ ?}$$
Cela revient \`a demander:
est-ce que la suite $\{\e_q\}_{q\in E}$ rang\'ee par valeur absolue
$|q|$ croissante est une base de $X_E$~?  En d'autres termes, la suite
des multiplicateurs idempotents relatifs $T_n:X_E\to X_E$ d\'efinie
par
$$T_n\e_q=\left\{\begin{array}{cl}
    \e_q&\hbox{si }|q|\le n\\
    0&\hbox{sinon}
\end{array}\right.$$
est-elle uniform\'ement born\'ee sur $n$~? Soit $E=\Z$. Un \'el\'ement
de r\'eponse classique est le suivant. 
$$
\|T_n\|_{{\fam0 L}^2(\sT)\to{\fam0 L}^2(\sT)}=1\ ,\ 
\|T_n\|_{{\fam0 L}^1(\sT)\to{\fam0 L}^1(\sT)}= 
\|T_n\|_{\script{C}(\sT)\to\script{C}(\sT)}\asymp\log n.
$$
On sait de plus que les $T_n$ sont aussi uniform\'ement born\'es sur
${\fam0 L}^p(\T)$, $1<p<\infty$.

\questsubf 
Est-ce que la somme de la s\'erie $\sum
\varrho_q\e^{\ii\vartheta_q}\e_q$ d\'epend de l'ordre dans lequel on
somme les fr\'equences~? Cette question est \'equivalente \`a la
suivante: la nature de $\sum \varrho_q\e^{\ii\vartheta_q}\e_q$
d\'epend-elle des phases $\vartheta_q$~? En termes fonctionnels,
$\{\e_q\}_{q\in E}$ forme-t-elle une suite basique inconditionnelle
dans $X$~? Cette question s'\'enonce aussi en termes de
multiplicateurs relatifs: la famille des $T_\epsilon:X_E\to X_E$ avec
$$T_\epsilon\e_q=\epsilon_q\e_q\hbox{ et }\epsilon_q=\pm1$$
est-elle
uniform\'ement born\'ee sur les choix de signes $\epsilon$~?  Un
\'el\'ement de r\'eponse classique est le suivant. Soit $E=\Z$. Alors
$$\|T_\epsilon\|_{{\fam0 L}^2(\sT)\to {\fam0 L}^2(\sT)}=1;$$
si $p\ne2$, il existe un
choix de signes $\epsilon$ tel que $T_\epsilon$ n'est pas born\'e sur
${\fam0 L}^p(\T)$.

\questsubf
Peut-on am\'eliorer ce ph\'enom\`ene en restreignant
le spectre $E$~?  Cette question m\`ene \`a l'\'etude des
sous-ensembles lacunaires de $\Z$, et a \'et\'e trait\'ee en d\'etail
par Walter Rudin.

Nous choisissons la notion de multiplicateur relatif comme
dictionnaire entre l'a\-na\-ly\-se harmonique et l'analyse fonctionnelle.
Nous d\'eveloppons une technique pour le calcul de la norme de
familles $\{T_\epsilon\}$ de multiplicateurs relatifs. Celle-ci nous
permet de traiter les questions suivantes.

\questsubf\label{intro:q4} 
Est-ce que la norme de $f\in X_E$ d\'epend seulement
de l'intensit\'e $\varrho_q$ de ses fr\'equences $\e_q$, et non pas de
leur phase $\vartheta_q$~? Cela revient \`a demander si $\{\e_q\}_{q\in
  E}$ est une suite basique $1$-inconditionnelle complexe dans $X$.

\questsubf\label{intro:q5}
Est-ce que l'on a pour tout choix de signes ``r\'eel''
$\pm$
$$
\Bigl\|\sum_{q\in E}\pm a_q\e_q\Bigr\|_X= \Bigl\|\sum _{q\in E}
a_q\e_q\Bigr\|_X\hbox{ ?}
$$
En d'autres mots, est-ce que $\{\e_q\}_{q\in E}$ est une suite
basique $1$-inconditionnelle r\'eelle dans $X$~?

La r\'eponse est d\'ecevante dans le cas des espaces ${\fam0 L}^p(\T)$, $p$
non entier pair: seules les fonctions dont le spectre a au plus deux
\'el\'ements v\'erifient ces deux propri\'et\'es. Pour mieux cerner le
ph\'enom\`ene, nous proposons d'introduire la question
presqu'iso\-m\'e\-tri\-que suivante.

\questsubf\label{intro:q6}
Est-ce que la norme de $f\in X_E$ d\'epend
arbitrairement peu de la phase $\vartheta_q$ de ses fr\'equences
$\e_q$~? De mani\`ere pr\'ecise, dans quel cas existe-t-il, pour
chaque $\eps>0$, un sous-ensemble $F\se E$ fini tel que
$$\Bigl\| \sum_{q\in E\setminus F}\varrho_q\e^{\ii\vartheta_q}\e_q
\Bigr\|_X\le(1+\eps)\Bigl\| \sum_{q\in E\setminus F}\varrho_q\e_q
\Bigr\|_X\hbox{ ?}
$$
Dans le cas $X=\script{C}(\T)$, cela signifiera que $E$ est un
ensemble de constante de Sidon ``asymptotiquement $1$''.  De m\^eme,
peut-on choisir pour chaque $\eps>0$ un ensemble fini $F$ tel que pour
tout choix de signe ``r\'eel'' $\pm$
$$
\Bigl\| \sum_{q\in E\setminus F}\pm a_q\e_q
\Bigr\|_X\le(1+\eps)\Bigl\| \sum_{q\in E\setminus F}a_q\e_q
\Bigr\|_X\hbox{ ?}
$$

Toutes ces questions s'agr\`egent autour d'un fait bien connu: sommer
la s\'erie de Fourier de $f$ est une tr\`es mauvaise mani\`ere
d'approcher la fonction $f$ d\`es que l'erreur consid\'er\'ee n'est
pas quadratique. On sait qu'il est alors utile de rechercher des
m\'ethodes de sommation plus lisses, c'est-\`a-dire d'autres suites
approximantes plus r\'eguli\`eres. Il s'agit l\`a de suites
d'op\'erateurs de rang fini sur $X_E$ qui approchent ponctuellement
l'identit\'e de $X_E$.  Nous pourrons toujours supposer que ces
op\'erateurs sont des multiplicateurs. Une premi\`ere question est la
suivante.

\questsubf\label{intro:q7} \
Existe-t-il une suite approximante $\{T_n\}$ de
multiplicateurs idempotents~? Cela revient \`a demander: existe-t-il
une d\'ecomposition de $X_E$ en sous-espaces $X_{E_k}$ de dimension
finie avec
\begin{equation}\label{intro:fdd}
X_E=\bigoplus X_{E_k}\quad\hbox{et}\quad A_k:X_E\to X_{E_k}\ ,\ \e_q\mapsto
\left\{\begin{array}{cl}
    \e_q&\hbox{si }q\in E_k\\
    0&\hbox{sinon}
\end{array}\right.
\end{equation}
telle que la suite des $T_n=A_1+\dots+A_n$ est uniform\'ement
born\'ee sur $n$~?  Soit $E=\Z$. Alors la r\'eponse est identique \`a
la r\'eponse de la question \ref{intro:q1}.

Mais nous pouvons produire dans ce cadre plus g\'en\'eral des
d\'ecompositions inconditionnelles de $X_E$ en r\'eponse \`a la
question suivante.

\questsubf
Pour quels espaces $X$ et spectres $E$ existe-t-il une
d\'ecomposition comme ci-dessus telle que la famille des multiplicateurs
\begin{equation}\label{intro:fddi}
\sum_{k=1}^n\epsilon_kA_k\quad\hbox{avec }
n\ge1\hbox{ et }\epsilon_k=\pm 1
\end{equation}
est uniform\'ement born\'ee~?  Littlewood et Paley ont
montr\'e que la partition de $\Z$ en $\Z=\bigcup E_k$ avec $E_0=\{0\}$ et
$E_k=\{j:2^{k-1}\le |j|<2^k\}$ donne une d\'ecomposition
inconditionnelle des espaces ${\fam0 L}^p(\T)$ avec $1<p<\infty$. D'apr\`es la
r\'eponse \`a la question \ref{intro:q7}, ce n'est pas
le cas \af des espaces ${\fam0 L}^1(\T)$ et $\script{C}(\T)$. 
Une \'etude fine de telles partitions a \'et\'e entreprise par Kathryn
Hare et Ivo Klemes.

Notre technique permet de traiter la question suivante.

\questsubf\label{intro:q9} 
Pour quels espaces $X$ et spectres $E$ existe-t-il une
d\'ecomposition du type \Ref{intro:fdd} telle que 
$$
\Bigl\|\sum\epsilon_kA_kf\Bigr\|_X=\|f\|_X \hbox{ pour
  tout choix de signes }\epsilon_k\hbox{ ?}
$$
La r\'eponse d\'ependra de la nature du choix de signes, qui peut \^etre r\'eel ou
complexe. 

Il est instructif de noter que l'espace de Hardy $H^1(\T)$ n'admet pas de
d\'ecomposition du type \Ref{intro:fdd}. $H^1(\T)$ admet n\'eanmoins
des suites approximantes de multiplicateurs et il existe m\^eme des
suites approximantes de multiplicateurs inconditionnelles au sens o\`u
la famille \Ref{intro:fddi} est uniform\'ememt born\'ee. Cela motive
la question suivante, qui est la plus g\'en\'erale dans notre
contexte.

\questsubf\label{intro:q10}
Quels sont les espaces $X$ et spectres $E$ tels que pour chaque
$\eps>0$ il existe une suite approximante $\{T_n\}$ sur $X_E$ telle
que
$$
\sup_{\hbox{\scriptsize signes }\epsilon_n}
\Bigl\|\sum\epsilon_n(T_n-T_{n-1})\Bigr\|_X\le1+\eps
$$
En termes fonctionnels, $X_E$ a-t-il la propri\'et\'e d'approximation
inconditionnelle m\'etrique~? Il faudra distinguer le cas des signes
complexes et r\'eels.

\subsection{Propri\'et\'e d'approximation inconditionnelle
  m\'e\-tri\-que}

Comme nos questions distinguent les choix de signe r\'eel et complexe,
nous proposons pour la fluidit\'e de l'expos\'e de fixer un choix de
signes $\U$ qui sera $\U=\T=\{\epsilon\in\C:|\epsilon|=1\}$ dans le
cas complexe et $\U=\D=\{-1,1\}$ dans le cas r\'eel.

Seule la question \ref{intro:q10} n'impose pas au pr\'ealable de forme
particuli\`ere \`a la suite de multiplicateurs qui est cens\'ee
r\'ealiser la propri\'et\'e consid\'er\'ee. Afin d'\'etablir un lien
entre la \umap/ et la structure du spectre $E$,
nous faisons le d\'etour par une \'etude g\'en\'erale de cette
propri\'et\'e dans le cadre des espaces de Banach s\'eparables.

\subsubsection{Amorce et queue d'un espace de Banach}

Peter G.\ Casazza et Nigel J.\ Kalton ont d\'ecouvert le crit\`ere suivant:
\begin{prpsubf}
Soit $X$ un espace de Banach s\'eparable. $X$ a la \umap/ si et
seulement s'il existe une suite approximante $\{T_k\}$ telle que 
$$
\sup_{\epsilon\in\sU}\|T_k+\epsilon(\Id-T_k)\|_{\script{L}(X)}
\toll{k\to\infty}1.
$$
\end{prpsubf}
Ceci exprime que la constante d'inconditionnalit\'e entre l'amorce
$T_kX$ et la queue $(\Id-T_k)X$ de l'espace $X$
s'am\'eliore asymptotiquement jusqu'\`a l'optimum pour $k\to\infty$.

La \umap/ s'exprime de mani\`ere plus \'el\'ementaire encore si l'on
choisit d'autres notions adapt\'ees d'amorce et de queue. Nous
proposons en particulier la d\'efinition suivante.
\begin{dfnsubf}
Soit $\tau$ une topologie d'espace vectoriel topologique sur $X$. $X$
a la propri\'et\'e $(u(\tau))$ de $\tau$-inconditionnalit\'e si pour
chaque $x\in X$ et toute suite born\'ee $\{y_j\}$ $\tau$-nulle
l'oscillation 
$$
\osc_{\epsilon\in\sU}\|x+\epsilon y_j\|_X=
\sup_{\delta,\epsilon\in\sU}
\bigl(\|x+\epsilon y_j\|-\|x+\delta y_j\|\bigr)
$$
forme elle-m\^eme une suite nulle.
\end{dfnsubf}

Nous avons alors le th\'eor\`eme suivant.
\begin{thmsubf}\label{intro:THM}
Soit $X$ un espace de Banach s\'eparable de cotype fini avec la
propri\'et\'e $(u(\tau))$. Si $X$ admet une suite approximante
$\{T_k\}$ inconditionnelle et commutative 
telle que $T_kx\mathop{\to}\limits^\tau x$ uniform\'ement
sur la boule unit\'e $B_X$, alors des combinaisons convexes
successives $\{U_j\}$ de $\{T_k\}$ r\'ealisent la \umap/.
\end{thmsubf}
{\it Esquisse de preuve.\/ }
On construit ces combinaisons convexes successives par le biais de
d\'ecompositions skipped blocking. En effet, la propri\'et\'e
$(u(\tau))$ a l'effet suivant sur $\{T_k\}$. Pour chaque $\eps>0$, il
existe une sous-suite $\{S_k=T_{n_k}\}$ telle que toute suite de blocs
$S_{b_k}-S_{a_k}$ obtenue en sautant les blocs $S_{a_{k+1}}-S_{b_k}$
se somme de mani\`ere $(1+\eps)$-inconditionnelle.\\
Soit $n\ge1$. Pour chaque $j$, $1\le j\le n$, la suite de blocs
obtenue en sautant $S_{kn+j}-S_{kn+j-1}$ pour $k\ge0$ est
$(1+\eps)$-inconditionnelle. Il s'agit alors d'estimer la moyenne sur
$j$ de ces suites de blocs. On obtient une suite approximante et
l'hypoth\`ese de cotype fini permet de contr\^oler l'apport des blocs
saut\'es.\\
Alors $X$ a la \umap/ parce que $n$ et $\eps$ sont arbitraires.\eck

\subsubsection{Amorce et queue en termes de spectre de Fourier}

Lorsqu'on consid\`ere l'espace invariant par translation $X_E$, une
amorce et une queue naturelle sont les espaces $X_F$ et $X_{E\setminus
  G}$ pour $F$ et $G$ des sous-ensembles finis de $E$. Nous avons
concr\`etement le lemme suivant.
\begin{lemsubf}
$X_E$ a $(u(\tau_f))$, o\`u $\tau_f$ est la topologie
$$ f_n\mathop{\to}\limits^{\tau_f}0\ \Ssi\ \forall k\ \widehat{f_n}(k)\to0 $$
de convergence simple des coefficients de Fourier, si et seulement si
$E$ est bloc-inconditionnel dans $X$ au sens suivant: quels que soient
$\eps>0$ et $F\se E$ fini, il existe $G\se E$ fini tel que pour
$f\in B_{X_F}$ et $g\in B_{X_{E\setminus G}}$
$$\osc_{\epsilon\in\sU}\|f+\epsilon g\|_X=
\sup_{\delta,\epsilon\in\sU}
\bigl(\|f+\epsilon g\|-\|f+\delta g\|\bigr)\le\eps.$$
\end{lemsubf}

Le th\'eor\`eme \ref{intro:THM}  s'\'enonce donc ainsi dans 
ce contexte particulier.
\begin{thmsubf}\label{intro:bloc}
Soit $E\se\Z$ et $X$ un espace de Banach homog\`ene de fonctions sur
le tore $\T$. Si $X_E$ a la \umap/, alors $E$ est bloc-inconditionnel
dans $X$. Inversement, si $E$ est bloc-inconditionnel dans $X$ et de
plus $X_E$ a la propri\'et\'e d'approximation inconditionnelle et un
cotype fini, alors $X_E$ a la \umap/. En particulier, on a
\begin{itemize}
\item[$(i)$]Soit $1<p<\infty$. ${\fam0 L}^p_E(\T)$ a la \umap/ si et seulement si
$E$ est bloc-inconditionnel dans ${\fam0 L}^p(\T)$.
\item[$(ii)$]${\fam0 L}^1_E(\T)$ a la \umap/ si et seulement si ${\fam0 L}^1_E(\T)$ a la
propri\'et\'e d'approximation inconditionnelle et $E$ est
bloc-inconditionnel dans ${\fam0 L}^1(\T)$.
\item[$(iii)$]Si $E$ est bloc-inconditionnel dans $\script{C}(\T)$ et $E$
est un ensemble de Sidon, alors $\script{C}_E(\T)$ a la \umap/.
\end{itemize}
\end{thmsubf}

Donnons une application de ce th\'eor\`eme.
\begin{prpsubf}
Soit $E=\{n_k\}\se\Z$. Si $n_{k+1}/n_k$ est un entier impair pour
tout $k$, alors $\script{C}_E(\T)$ a la \umap/ r\'eelle.
\end{prpsubf}
{\it Preuve.\/ }
Comme $E$ est n\'ecessairement un ensemble de Sidon, il suffit de
v\'erifier que $E$ est bloc-inconditionnel. Soient
$\eps>0$ et $F\se E\cap[-n,n]$. Soit $l$ tel que $|n_l|\ge \pi n/\eps$
et $G=\{\lst n{l-1}\}$. Soit $f\in B_{\script{C}_F}$ et $g\in
B_{\script{C}_{E\setminus G}}$. Alors $g(t+\pi/n_l)=-g(t)$ par
hypoth\`ese et 
$$ |f(t+\pi/n_l)-f(t)|\le \pi/|n_l|\cdot\|f'\|_\infty\le\pi
n/|n_l|\le\eps $$
par l'in\'egalit\'e de Bernstein. Alors, pour un certain $u\in\T$
\begin{eqnarray*}
\|f-g\|_\infty&=&|f(u)+g(u+\pi/n_l)|\\
&\le&|f(u+\pi/n_l)+g(u+\pi/n_l)|+\eps\\
&\le&\|f+g\|_\infty+\eps.
\end{eqnarray*}
Donc $E$ est bloc-inconditionnel au sens r\'eel.\eck

En particulier, soit la suite g\'eom\'etrique $G=\{3^k\}$. Alors
$\script{C}_G(\T)$ et $\script{C}_{G\cup-G}(\T)$ ont la \umap/
r\'eelle.

\questsubf\label{intro:q2.1}
Qu'en est-il de la \umap/ complexe et qu'en est-il de la suite
g\'eom\'etrique $G=\{2^k\}$~?
 
\subsection{Norme de multiplicateurs et conditions
  combinatoires}\label{intro:nmcc}

Nous proposons ici une m\'ethode uniforme pour r\'epondre aux
questions \ref{intro:q4}, \ref{intro:q5}, \ref{intro:q6},
\ref{intro:q9} et \ref{intro:q10}. En effet, les questions
\ref{intro:q4}, \ref{intro:q5} et \ref{intro:q6} reviennent \`a
\'evaluer l'oscillation de la norme 
$$
\Theta(\epsilon,a)=
\|\epsilon_0a_0\e_{r_0}+\dots+\epsilon_ma_m\e_{r_m}\|_X.
$$
La question \ref{intro:q9} revient \`a \'evaluer l'oscillation de la
norme 
\begin{eqnarray*}
\Psi(\epsilon,a)&=&\Theta((\overbrace{1,\dots,1}^j,
\overbrace{\epsilon,\dots,\epsilon}^{m-j}),a)\\
&=&\|a_0\e_{r_0}+\dots+a_j\e_{r_j}+\epsilon a_{j+1}\e_{r_{j+1}}+\dots+
\epsilon a_m\e_{r_m}\|_X
\end{eqnarray*}
Par le th\'eor\`eme \ref{intro:bloc}, la question \ref{intro:q10}
revient \`a \'etudier cette m\^eme expression dans le cas particulier
o\`u on fait un saut de grandeur arbitraire entre $r_j$ et
$r_{j+1}$.\\
Dans le cas des espaces $X={\fam0 L}^p(\T)$, $p$ entier pair, ces normes sont
des polyn\^omes en $\epsilon$, $\epsilon^{-1}$, $a$ et $\bar a$. Dans
le cas des espaces $X={\fam0 L}^p(\T)$, $p$ non entier pair, elles s'expriment
comme des s\'eries. Il n'y a pas moyen d'exprimer ces normes comme
fonction $\script{C}^\infty$ pour $X=\script{C}(\T)$.

Soit $X={\fam0 L}^p(\T)$. D\'eveloppons $\Theta(\epsilon,a)$. Posons
$q_i=r_i-r_0$. On peut supposer $\epsilon_0=1$ et $a_0=1$. Nous
utilisons la notation suivante:
$$
{x\choose\alpha}=
\frac{x(x-1)\cdots(x-n+1)}{\alpha_1!\alpha_2!\dots}\quad\hbox{pour
  }\alpha\in\N^m\hbox{ tel que }\sum\alpha_i=n
$$
Alors, si $|a_1|,\dots,|a_m|<1/m$ lorsque $p$ n'est pas un entier pair
et sans restriction sinon, 
\begin{eqnarray}
\nonumber
\Theta(\epsilon,a)
&=&\int\biggl|\sum_{n\ge0}\bip{p/2}n
\biggl(\sum_{i=1}^m\epsilon_ia_i\e_{q_i}
\biggr)^n\biggr|^2\\
\nonumber&=&\int\biggl|
\sum_{n\ge0}\bip{p/2}n
\sum_{\scriptstyle\alpha:\lst\alpha m\ge0
 \atop\scriptstyle\sm\alpha m=n}
\bi n\alpha
\epsilon^\alpha a^\alpha
\e_{\lower1pt\hbox{$\Sigma$}\alpha_iq_i}
\biggr|^2\\
\nonumber&=&\int\biggl|\sum_{\alpha\in\sN^m}
\bip{p/2}\alpha
\epsilon^\alpha a^\alpha
\e_{\lower1pt\hbox{$\Sigma$}\alpha_iq_i}\biggr|^2\\
\nonumber&=&\sum_{R\in\script{R}}
\biggl|\sum_{\alpha\in R}
\bip{p/2}\alpha\epsilon^\alpha a^\alpha
\biggr|^2\\
\nonumber&=&
\sum_{\alpha\in\sN^m}
{\bip{p/2}\alpha}^2|a|^{2\alpha}+
\sum_{\scriptstyle\alpha\ne\beta\in\sN^m\atop\scriptstyle\alpha\sim\beta}
\bip{p/2}\alpha\bip{p/2}\beta
\epsilon^{\alpha-\beta}a^\alpha \bar a^\beta\label{intro:gc:res}
\end{eqnarray}
o\`u $\script{R}$ est la partition de $\N^m$ induite par la relation
d'\'equivalence 
$$\alpha\sim\beta\ssi\sum\alpha_iq_i=\sum\beta_iq_i.$$ 

Nous pouvons r\'epondre imm\'ediatement aux questions \ref{intro:q4}
et \ref{intro:q5} pour $X={\fam0 L}^p(\T)$.

\subsubsection[Suites basiques $1$-inconditionnelles complexes]
{Question \ref{intro:q4}: suites basiques $1$-inconditionnelles complexes}
\label{sec:intro:q4}

Soient $r_0,\dots r_m$ sont choisis dans $E$, alors \Ref{intro:gc:res}
doit \^etre constante pour $a\in\{|z|<1/m\}^m$ et $\epsilon\in\T^m$. Cela
veut dire que pour tous $\alpha\ne\beta\in\N^m$,
$$
\sum\alpha_iq_i\ne\sum\beta_iq_i\quad\hbox{ou}\quad
\bip{p/2}\alpha\bip{p/2}\beta=0.
$$

\bloc
Si $p$ n'est pas un entier pair, alors
$\bip{p/2}\alpha\bip{p/2}\beta\ne0$ pour tous $\alpha,\beta\in\N^m$ et
on a les relations arithm\'etiques suivantes sur $q_1,q_2,0$:
\begin{eqnarray*}
\overbrace{q_1+\dots+q_1}^{|q_2|}=
\overbrace{q_2+\dots+q_2}^{|q_1|}&&\mbox{si }q_1q_2>0;\\
\overbrace{q_1+\dots+q_1}^{|q_2|}+\overbrace{q_2+\dots+q_2}^{|q_1|}
=0&&\mbox{sinon.}
\end{eqnarray*}
Il suffit donc de prendre $\alpha=(|q_2|,0,\dots)$,
$\beta=(|q_1|,0,\dots)$ et $\alpha=(|q_2|,|q_1|,0,\dots)$,
$\beta=(0,\dots)$ respectivement pour conclure que $\{r_0,r_1,r_2\}$
n'est pas une suite basique $1$-inconditionnelle complexe dans ${\fam0
  L}^p(\T)$ si $p$ n'est pas un entier pair.

\bloc
Si $p$ est un entier pair, $\bip{p/2}\alpha\bip{p/2}\beta=0$ si et
seulement si 
$$\sum\alpha_i>p/2\quad\hbox{ou}\quad\sum\beta_i>p/2.$$ 
On obtient que
$E$ est une suite basique $1$-inconditionnelle dans ${\fam0 L}^p(\T)$ si et
seulement si $E$ est ``$p$-ind\'ependant'', c'est-\`a-dire que
$\sum\alpha_i(r_i-r_0)\ne\sum\beta_i(r_i-r_0)$ pour tous
$r_0,\dots,r_m\in E$ et $\alpha\ne\beta\in\N^m$ tels que
$\sum\alpha_i,\sum\beta_i\le p/2$. Cette condition est \'equivalente
\`a: tout entier $n\in\Z$ s'\'ecrit de mani\`ere au plus unique comme
somme de $p/2$ \'el\'ements de $E$.

\subsubsection[Suites basiques $1$-inconditionnelles r\'eelles]
{Question \ref{intro:q5}: suites basiques $1$-inconditionnelles r\'eelles}

Les suites basiques $1$-inconditionnelles r\'eelles et complexes
co\"\i ncident et la r\'eponse \`a la question \ref{intro:q5} est
identique \`a la r\'eponse \`a la question \ref{intro:q4}. En effet,
d\`es qu'une relation arithm\'etique $\sum(\alpha_i-\beta_i)q_i$
p\`ese sur $E$, on peut supposer que $\alpha_i-\beta_i$ est impair
pour au moins un $i$ en simplifiant la relation par le plus grand
diviseur commun des $\alpha_i-\beta_i$. Mais alors \Ref{intro:gc:res}
n'est pas une fonction constante pour $\epsilon_i$ r\'eel.

Cette propri\'et\'e est propre au tore $\T$. En effet, par exemple la
suite des fonctions de Rademacher est $1$-inconditionnelle r\'eelle
dans $\script{C}(\D^\infty)$, alors que sa constante
d'inconditionnalit\'e complexe est $\pi/2$.

\subsubsection[Suites basiques inconditionnelles m\'etriques]
{Question \ref{intro:q6}: suites basiques inconditionnelles m\'etriques}

On peut m\^eme tirer des cons\'equences utiles du calcul de
\Ref{intro:gc:res} dans le cas pres\-qu'iso\-m\'e\-tri\-que. Il faut pour cela
prendre la pr\'ecaution suivante qui permet un passage \`a la
limite. Soit $0<\varrho<1/m$. Alors 
$$\big\{\Theta\colon\U^m\times\{|z|\le\varrho\}^m\to\R^+:q_1,\dots,
q_m\in\Z^m\big\}$$
est un sous-ensemble relativement compact de
$\script{C}^\infty(\U^m\times\{|z|\le\varrho\}^m)$. Il en d\'ecoule
que si $E$ est une suite basique inconditionnelle m\'etrique, alors
certains coefficients de \Ref{intro:gc:res} deviennent arbitrairement
petits lorsque $q_1,\dots,q_m$ sont choisis grands.

Donnons deux cons\'equences de ce raisonnement.
\begin{prpsubf}
Soit $E\se\Z$.
\begin{itemize}
\item[$(i)$]Soit $p$ un entier pair. Si $E$ est une suite basique inconditionnelle
m\'etrique r\'eelle, alors $E$ est en fait une suite basique
$1$-inconditionnelle complexe \`a un ensemble fini pr\`es.
\item[$(ii)$]Si $E$ est un ensemble de Sidon de constante asymptotiquement $1$,
alors 
$$\XE=\sup_{G\se E\hbox{\scriptsize\ fini}}\,
\inf\bigl\{|\zeta_1p_1+\dots+\zeta_mp_m|:
\lst p m\in E\setminus G\hbox{ distincts}\bigr\}>0$$
pour tout $m\ge1$ et $\zeta\in{\Z^*}^m$.
\end{itemize}
\end{prpsubf}
On peut exprimer cette derni\`ere propri\'et\'e en disant que la
relation arithm\'etique $\zeta$ ne persiste pas sur $E$.

\subsubsection[Propri\'et\'e d'approximation inconditionnelle m\'etrique]
{Question \ref{intro:q10}: propri\'et\'e d'approximation inconditionnelle m\'e\-tri\-que}

On peut appliquer la technique du paragraphe pr\'ec\'edent en
observant que si $X_E$ a la \umap/, alors
$$
\osc_{\epsilon\in\sU}\Psi(\epsilon,a)
\toll{\lstp r{j+1}m\in E\to\infty}0.
$$

\begin{dfnsubf}
$E$ a la propri\'et\'e $\J{n}$ de bloc-ind\'ependance si pour tout
$F\se E$ fini il existe $G\se E$ fini tel que si un $k\in\Z$ admet
deux repr\'esentations comme somme de $n$ \'el\'ements de
$F\cup(E\setminus G)$
$$\sm pn=k=\sm{p'}n,$$ 
alors
$$\mes{\{j:p_j\in F\}}\quad\hbox{et}\quad\mes{\{j:p'_j\in F\}}$$
sont \'egaux \(choix de signes complexe $\U=\T$\) ou de m\^eme parit\'e
\(choix de signes r\'eel $\U=\D$\).
\end{dfnsubf}
\begin{thmsubf}
Soit $E\se\Z$.
\begin{itemize}
\item [$(i)$]Si $X={\fam0 L}^p(\T)$, $p$ entier pair, alors ${\fam0 L}^p_E(\T)$ a la \umap/ si et
seulement si $E$ satisfait $\J{p/2}$.
\item [$(ii)$]Si $X={\fam0 L}^p(\T)$, $p$ non entier pair, ou $X=\script{C}(\T)$, alors
$X_E$ a la \umap/ seulement si $E$ satisfait
$$\XE=\sup_{G\se E\hbox{\scriptsize\ fini}}\,
\inf\bigl\{|\zeta_1p_1+\dots+\zeta_mp_m|:
\lst p m\in E\setminus G\hbox{ distincts}\bigr\}>0$$
pour tout $m\ge1$ et $\zeta\in{\Z^*}^m$ tel que $\sum\zeta_i$ est non
nul \(cas complexe\) ou impair \(cas r\'eel\).
\end{itemize}
\end{thmsubf}

On obtient la hi\'erarchie suivante.
$$
\parbox{13mm}{$\SCE$ a\\\umap/}\imp
\parbox{25mm}{$\SLPE$ a \umap/,\\$p$ non entier pair}
\imp\dots\imp
\parbox{15mm}{$\SLE{2n+2}$\\ a \umap/}\imp
\parbox{14mm}{$\SLE{2n}$ a\\\umap/}
\imp\dots\imp\parbox{13mm}{$\SLE{2}$ a\\\umap/.}
$$

Nous pouvons r\'epondre \`a la question \ref{intro:q2.1}. Soit
$G=\{j^k\}$ avec $j\in\Z\setminus\{-1,0,1\}$ et consid\'erons
$\zeta=(j,-1)$. Alors $\langle\zeta,G\rangle=0$. Donc
$\script{C}_G(\T)$ n'a pas la \umap/ complexe. $\script{C}_G(\T)$ n'a
pas la \umap/ r\'eelle si $j$ est pair.

\subsubsection{Deux exemples}

\`A l'aide de nos conditions arithm\'etiques, nous sommes \`a m\^eme
de prouver la proposition suivante.
\begin{prpsubf}
Soit $\sigma>1$ et $E$ la suite des parties enti\`eres de
$\sigma^k$. Alors les assertions suivantes sont \'equivalentes.
\begin{itemize}
\item [$(i)$]$\sigma$ est un nombre transcendant.
\item [$(ii)$]${\fam0 L}^p_E(\T)$ a la \umap/ complexe pour tout $p$ entier pair.
\item [$(iii)$]$E$ est une suite basique inconditionnelle m\'etrique dans chaque
${\fam0 L}^p(\T)$, $p$ entier pair.
\item [$(iv)$]Pour chaque $m$ donn\'e, la constante de Sidon des sous-ensembles 
\`a $m$ \'el\'ements de queues de $E$ est asymptotiquement $1$.
\end{itemize}
\end{prpsubf}

Nous obtenons aussi la proposition suivante.
\begin{prpsubf}
Soit $E$ la suite des bicarr\'es. ${\fam0 L}^p_E(\T)$ a la \umap/ r\'eelle
seulement si $p=2$ ou $p=4$.
\end{prpsubf}
{\it Preuve.\/ }
$E$ ne satisfait pas la propri\'et\'e de bloc-ind\'ependance $\J{3}$
r\'eelle. En effet, Ramanujan a d\'ecouvert l'\'egalit\'e suivante
pour tout $n$:
$$
(4n^5-5n)^4+(6n^4-3)^4+(4n^4+1)^4=
(4n^5+n)^4+(2n^4-1)^4+3^4.\eqno{\ecks}
$$

\subsection{Impact de la croissance du spectre}
\label{sec:intro:croiss}

Nous d\'emontrons de mani\`ere directe le r\'esultat positif suivant.
\begin{thmsubf}
Soit $E=\{n_k\}\se\Z$ tel que $n_{k+1}/n_k\to\infty$. Alors la suite
des projections associ\'ee \`a $E$ r\'ealise la \umap/ complexe dans
$\SCE$ et $E$ est un ensemble de Sidon de constante asymptotiquement
$1$. Dans l'hypoth\`ese o\`u les rapports $n_{k+1}/n_k$ sont tous
entiers, la r\'eciproque vaut.
\end{thmsubf}
\begin{corsubf}
Alors $E$ est une suite basique inconditionnelle m\'etrique dans tout
espace de Banach homog\`ene $X$ de fonctions sur $\T$. De plus, $X_E$ a la
\umap/ complexe.
\end{corsubf}
{\it Esquisse de preuve.\/ }
Nous prouvons concr\`etement que si $n_{k+1}/n_k\to\infty$, alors quel
que soit $\eps>0$ il existe $l\ge1$ tel que pour toute fonction
$f=\sum a_k\e_{n_k}$
\begin{equation}\label{p6}
\|f\|_\infty\ge(1-\eps)\Bigl(\Bigl\|\sum_{k\le
  l}a_k\e_{n_k}\Bigr\|_\infty+\sum_{k>l}|a_k|\Bigr).
\end{equation}
Cela revient \`a dire que la suite $\{\pi_k\}$ de projections
associ\'ee \`a la base $E$ r\'ealise la $1/(1-\eps)$-\uap/. Pour
obtenir l'in\'egalit\'e $(\ref{p6})$, on utilise une r\'ecurrence
bas\'ee sur l'id\'ee suivante. 

Soit $u\in\T$ tel que $\|\pi_kf\|_\infty=|\pi_kf(u)|$. 
Il existe alors $v\in\T$ tel que
$$ 
|u-v|\le\pi/|n_{k+1}|\quad\mbox{et}\quad
|\pi_kf(u)+a_{k+1}\e_{n_{k+1}}(v)|=\|\pi_kf\|_\infty+|a_{k+1}|.$$
De plus, dans ce cas, 
$$
|\pi_kf(u)-\pi_kf(v)|\le
|u-v|\,\|\pi_kf'\|_\infty\le
\pi|n_k/n_{k+1}|\,\|\pi_kf\|_\infty.
$$ 

En r\'esum\'e, $a_{k+1}\e_{n_{k+1}}$ a le m\^eme argument que $\pi_kf$
tr\`es pr\`es du maximum de $|\pi_kf|$, et $\pi_kf$ varie peu. 

Mais alors 
\begin{eqnarray*}
\|\pi_kf(t)+a_{k+1}\e_{n_{k+1}}\|_\infty&\ge&|\pi_kf(v)+a_{k+1}\e_{n_{k+1}}(v)|\\
&\ge&\|\pi_kf\|_\infty+|a_{k+1}|-\pi|n_k/n_{k+1}|\|\pi_kf\|_\infty\\
&=&(1-\pi|n_k/n_{k+1}|)\|\pi_kf\|_\infty+|a_{k+1}|.
\end{eqnarray*}
On obtient $(\ref{p6})$ en r\'eit\'erant cet argument.\eck

Notre technique donne d'ailleurs l'estimation suivante de la constante
de Sidon des ensembles de Hadamard.
\begin{corsubf}
  Soit $E=\{n_k\}\se\Z$ et $q>\sqrt{\pi^2/2+1}$. Si $|n_{k+1}|\ge
  q|n_k|$, alors la constante de Sidon de $E$ est inf\'erieure ou
  \'egale \`a $1+\pi^2/(2q^2-2-\pi^2)$.
\end{corsubf}
Nous prouvons que cette estimation est optimale au sens o\`u
l'ensemble $E=\{0,1,q\}$, $q\ge2$, a pour constante
d'inconditionnalit\'e r\'eelle dans $\SC{}$ 
$$\bigl(\cos(\pi/(2q)\bigr)^{-1}\ge1+\pi^2/8\,q^{-2}.$$

\section{Introduction}
We study isometric and almost isometric counterparts to the following
two properties of a separable Banach space $Y$:

\vskipa {\bf (ubs) } $Y$ is the closed span of an unconditional basic
sequence;

\vskipa{\bf (uap) } $Y$ admits an unconditional finite dimensional
expansion of the identity.\vskipb

We focus on the case of translation invariant spaces of functions on
the torus group $\T$, which will provide us with a bunch of natural
examples. Namely, let $E$ be a subset of $\Z$ and $X$ be one of the
spaces $\SLP{}$ $(1\le p<\infty)$ or $\SC{}$. If $\{\e^{{\rm
 i}nt}\}_{n\in E}$ is an unconditional basic sequence (\ubs/ for
short) in $X$, then $E$ is known to satisfy strong conditions of
lacunarity: $E$ must be in Rudin's class \EL{q}, $q=p\vee2$, and a
Sidon set respectively. We raise the following question: what kind of
lacunarity is needed to get the following stronger property:

\vskipa{\bf (umbs) } $E$ is a metric unconditional basic sequence in
$X$: for any $\eps>0$, one may lower its unconditionality constant to
$1+\eps$ by removing a finite set from it.\vskipb

In the case of $\SC{}$, $E$ is a \umbs/ exactly when $E$ is a Sidon set
with constant asymptotically $1$.

In the same way, call $\{T_k\}$ an approximating sequence (\as/ for
short) for $Y$ if the $T_k$'s are finite rank operators that tend
strongly to the identity on $Y$; if such a sequence exists, then $Y$
has the bounded approximation property. Denote by $\Delta
T_k=T_k-T_{k-1}$ the difference sequence of $T_k$. Following
Rosenthal\index{Rosenthal, Haskell Paul}
(see \cite[\S1]{fe80}), we then say that $Y$ has the unconditional
approximation property (\uap/ for short) if it admits an \as/
$\{T_k\}$ such that for some $C$
\begin{equation}\label{intro:uap}
\biggl\|\sum_{k=1}^n
\epsilon_k\Delta T_k
\biggr\|_{\script{L}(Y)} \le C
\qquad\mbox{for all $n$ and scalar $\epsilon_k$ with 
$|\epsilon_k|=1$.} 
\end{equation}
By the uniform boundedness principle, $(\ref{intro:uap})$ means
exactly that $\sum\Delta T_ky$ converges unconditionally for all $y\in
Y$. We now ask the following question: which conditions on $E$ do
yield the corresponding almost isometric (metric for short) property,
first introduced by Casazza and Kalton \cite[\S3]{ck91}~?

\vskipa{\bf (umap) } The span $Y=X_E$ of $E$ in $X$ has the metric
unconditional approximation property: for any $\eps>0$, one may lower
the constant $C$ in $(\ref{intro:uap})$ to $1+\eps$ by choosing an
adequate \as/ $\{T_k\}$.\vskipb

Several kinds of metric, \ie almost isometric properties have been
investigated in the last decade (see \cite{hww93}). There is a common
feature to these notions since 
Kalton's\index{Kalton, Nigel J.}
\cite{ka93}: they can be reconstructed from a 
corresponding interaction between some break and some tail of the
space. We prove that \umap/ is characterized by almost
$1$-unconditionality between a specific break and tail, that we coin
``block unconditionality''.\vskipb

Property \umap/ has been studied by Li \cite{li96} for $X=\SC{}$. He
obtains remarkably large examples of such sets $E$, in particular
Hilbert sets. Thus, the second property seems to be much weaker than
the first (although we do not know whether $\SCE$ has \umap/ for all \umbs/
$E$ in $\SC{}$: for sets of the latter kind, the natural sequence of
projections realizes \uap/ in $\SCE$, but we do not know whether it achieves
\umap/).\vskipb
 
In fact, both problems lead to strong arithmetical conditions on $E$
that are somewhat complementary to the property of quasi-independence
(see \cite[\S3]{pi81}). In order to obtain them, we apply
Forelli's\index{Forelli, Frank} \cite[Prop.\ 2]{fo64} and
Plotkin's\index{Plotkin, A. I.} \cite[Th.\ 1.4]{pl74} techniques 
in the study of isometric operators on ${\fam0 L}^p$: see Theorem
\ref{mub:thm} and Lemma \ref{umap:lem}. This may be done at once for
the projections associated to basic sequences of characters. In the
case of general metric unconditional approximating sequences, however,
we need a more thorough knowledge of their connection with the structure
of $E$: this is the duty of Theorem \ref{sbd:thm}. As in
Forelli's
and Plotkin's results, we obtain that the spaces
$X=\SLP{}$ with $p$ 
an even integer play a special r\^ole. For instance, they are the
only spaces which admit $1$-unconditional basic sequences \DE/ with
more than two elements: see Proposition \ref{mub:isom}.\vskipb

There is another fruitful point of view: we may consider elements of
$E$ as random variables on the probability space $(\T,dm)$. They have
uniform distribution and if they were independent, then our questions
would have trivial answers. In fact, they are strongly dependent: for
any $k,l\in\Z$, Rosenblatt's\index{Rosenblatt, Murray} \cite{ro56}
strong mixing coefficient\index{strong mixing}
$$
\sup\bigl\{|m[A\cap B]-m[A]m[B]|: A\in\sigma(\e^{\ii kt})\mbox{ and
 }B\in\sigma(\e^{\ii lt}) \bigr\}
$$
has its maximum value, $1/4$. But lacunarity of $E$ enhances their
independence in several weaker senses (see \cite{be90}). Properties
\umap/ and \umbs/ can be seen as an expression of almost independence
of elements of $E$ in the ``additive sense'', \ie when appearing in
sums. We show their relationship to the notions of pseudo-independence
(see \cite[\S4.2]{mu82}) and almost i.i.d.\ sequences (see
\cite{be87}).\vskipb

The gist of our results is the following: almost isometric properties
for spaces $X_E$ in ``little'' Fourier analysis may be read as a
smallness property of $E$. They rely in an essential way on the
arithmetical structure of $E$ and distinguish between real and complex
properties. In the case of $\SL{2n}{}$, $n$ integer, these arithmetical
conditions are in finite number and turn out to be sufficient, because
the norm of trigonometric polynomials is a polynomial expression in
these spaces. Furthermore, the number of conditions increases with
$n$ in that case. In the remaining cases of $\SLP{}$, $p$ not an even
integer, and $\SC{}$, these arithmetical conditions are infinitely many
and become much more coercive. In particular, if our properties are
satisfied in $\SC{}$, then they are satisfied in all spaces $\SLP{}$,
$1\le p<\infty$. \vskipb
 
We now turn to a detailed discussion of our results: in Section
$\ref{sect:mub}$, we first characterize the sets $E$ and values $p$
such that $E$ is a $1$-unconditional basic sequence in $\SLP{}$ (Prop.\ 
\ref{mub:isom}). Then we show how to treat similarly the almost
isometric case and obtain a range of arithmetical conditions \I{n} on
$E$ (Th.\ \ref{mub:thm}). These conditions turn out to be identical
whether one considers real or complex unconditionality: this is
surprising and in sharp contrast to what happens when $\T$ is replaced
by the Cantor group. They also do not distinguish amongst $\SLP{}$
spaces with $p$ not an even integer and $\SC{}$, but single out $\SLP{}$
with $p$ an even integer: this property does not ``interpolate''. This
is similar to the phenomena of
equimeasurability\index{equimeasurability} (see 
\cite[introduction]{ko91}) and $\script{C}^\infty$-smoothness of
norms\index{smoothness} 
(see \cite[Chapter V]{dgz93}). These facts may also be appreciated
from the point of view of natural renormings\index{renormings} of the
Hilbert space 
$\SLE{2}$.\vskipb

In Section $\ref{sect:ex}$, of purely arithmetical nature, we give
many examples of $1$-uncon\-ditional and metric unconditional basic
sequences through an investigation of property \I{n}. As expected
with lacunary series, number theoretic conditions show up (see
especially Prop.\ \ref{mub:trans}).\vskipb

In Section $\ref{sect:block}$, we first return to the general case of
a separable Banach space $Y$ and show how to connect the metric
unconditional approximation property with a simple property of ``block
unconditionality''. Then a skipped blocking technique invented by
Bourgain\index{Bourgain, Jean} and 
Rosenthal\index{Rosenthal, Haskell Paul} \cite{br80} 
gives a canonical way to construct an \as/ that realizes
\umap/ (Th.\ \ref{block:thm}).\vskipb

In Section \ref{sect:lpmap}, we introduce the $p$-additive approximation
property \lpap/ and its metric counterpart, \lpmap/. It may be
described as simply as \umap/. Then we connect \lpmap/ with the work
of Godefroy, Kalton, Li and Werner \cite{kw95,gkl96} on
subspaces of ${\fam0 L}^p$ which are almost isometric to $\ell_p$.\vskipb

Section \ref{sect:tis} focusses on \uap/ and \umap/ in the case of
translation invariant subspaces $X_E$. The property of block
unconditionality may then be expressed in terms of ``break'' and
``tail'' of $E$: see Theorem \ref{sbd:thm}.\vskipb

In Section $\ref{sect:umap}$, we proceed as in Section
$\ref{sect:mub}$ to obtain a range of arithmetical conditions \J{n}
for \umap/ and metric unconditional \fdd/ (Th.\ \ref{umap:thm} and
Prop.\ \ref{umap:prp:fdd}). These conditions are similar to \I{n},
but are decidedly weaker: see Proposition \ref{arith:csq}$(i)$. This
time, real and complex unconditionality differ; again spaces $\SLP{}$
with even $p$ are singled out.\vskipb

In Section $\ref{sect:arith}$, we continue the arithmetical
investigation begun in Section \ref{sect:ex} with property \J{n} and
obtain many examples for the $1$-un\-con\-di\-tio\-nal and the metric
unconditional approximation property.\vskipb

However, the main result of Section $\ref{sect:positif}$, Theorem
\ref{positif:thm}, shows how a rapid (and optimal) growth condition on
$E$ allows avoiding number theory in any case considered. We therefore
get a new class of examples for \umbs/, in particular Sidon sets of
constant asymptotically $1$, and \umap/. We also prove that
$\SC{\{3^k\}}$ has real \umap/ and that this is due to the oddness of
$3$ (Prop.\ \ref{res:geo}). A sharp estimate of the Sidon constant of
Hadamard sets is obtained as a byproduct (Cor.\ \ref{positif:cor}). 
\vskipb

Section $\ref{sect:comb}$ uses combinatorial tools to give some rough
information about the size of sets $E$ that satisfy our arithmetical
conditions. In particular, we answer a question of Li \cite{li96}: for
$X=\SC{}$ and for $X=\SLP{}$, $p\ne2,4$, the maximal density $d^*$ of
$E$ is zero if $X_E$ has \umap/ (Prop.\ \ref{comb:thm}). For
$X=\SL{4}{}$, our technique falls short of the expected result: we
just know that if $\SL{4}{E\cup\{a\}}$ has \umap/ for every $a\in\Z$,
then $d^*(E)=0$.\vskipb

Section $\ref{sect:proba}$ is an attempt to describe the relationship
between these notions and probabilistic independence. Specifically the
Rademacher and Steinhaus sequences show the way to a connection
between metric unconditionality and the almost i.i.d.\ sequences of
\cite{be87}. We note further that the arithmetical property
\I{\infty} of Section \ref{sect:mub} is equivalent to Murai's
\cite[\S4.2]{mu82} property of pseudo-independence.\vskipb

In Section \ref{sect:resume}, we collect our results on metric
unconditional basic sequences of characters and \umap/ in translation
invariant spaces. We conclude with open questions. \vskipb

{\bf Notation and definitions } Sections $\ref{sect:mub}$,
\ref{sect:tis}, \ref{sect:umap} and \ref{sect:positif} will take place
in the following framework. $(\T,dm)$ denotes the compact abelian
group $\{z\in\C:|z|=1\}$ endowed with its Haar measure $dm$; $m[A]$ is
the measure of a subset $A\se\T$. Let $\D=\{-1,1\}$. $\U$ will denote
either the complex ($\U=\T$) or real ($\U=\D$) choice of signs. For a
real function $f$ on $\U$, the oscillation\index{oscillation} of $f$ is
$$
\osc_{\epsilon\in\sU}f(\epsilon)= \sup_{\epsilon\in\sU}f(\epsilon)-
\inf_{\epsilon\in\sU}f(\epsilon).
$$
We shall study homogeneous\index{homogeneous Banach space} Banach
spaces $X$ of functions on $\T$ 
\cite[Chapter I.2]{ka68}, and especially the peculiar behaviour of the
following ones: $\SLP{}$ ($1\le p<\infty$), the space of $p$-integrable
functions with the norm $\|f\|_p=(\int|f|^pdm)^{1/p}$, and $\SC{}$, the
space of continuous functions with the norm
$\|f\|_\infty=\max\{|f(t)|:t\in\T\}$. $\script{M}(\T)$ is the dual of
$\SC{}$ realized as Radon measures on \T.

The dual group $\{\e_n\colon z\mapsto z^n:n\in\Z\}$ of $\T$ is
identified with $\Z$. We write $\mes{B}$ for the cardinal of a set
$B$. For a not necessarily increasing sequence \DEEE/, let \PTE/ be
the space of trigonometric polynomials spanned by [the characters in]
$E$. Let $X_E$ be the translation invariant subspace of those
elements in $X$ whose Fourier transform vanishes off $E$: for all
$f\in X_E$ and $n\notin E$, $\widehat{f}(n)=\int
f(t)\e_{-n}(t)dm(t)=0$. $X_E$ is also the closure of \PTE/ in
homogeneous $X$ \cite[Th.\ 2.12]{ka68}. Denote by $\pi_k:X_E\to X_E$
the orthogonal projection onto $X_{\{\lst n k\}}$. It is given by
$$
\pi_k(f)= \widehat{f}(n_1)\e_{n_1}+\dots+\widehat{f}(n_k)\e_{n_k}.
$$
Then the $\pi_k$ commute. They form an \as/ for $X_E$ if and only
if $E$ is a basic sequence. For a finite or cofinite $F\se E$,
$\pi_F$ is similarly the orthogonal projection of $X_E$ onto $X_F$.

Sections \ref{sect:block} and \ref{sect:lpmap} consider the general
case of a separable Banach space $X$. $B_X$ is the unit ball of $X$
and $\Id$ denotes the identity operator on $X$. For a given sequence
$\{U_k\}$, its difference sequence is $\Delta U_k=U_k-U_{k-1}$ (where
$U_0=0$).

The functional notions of \ubs/, \umbs/
are defined in \ref{mub:def}. The functional notions of \as/, \uap/
and \umap/ are defined in \ref{block:def}. Properties \lpap/ and
\lpmap/ are defined in \ref{str:def}. The functional property \UP/ of
block unconditionality is defined in \ref{block:block:def}. The sets
of arithmetical relations $\Zeta^m$ and $\Zeta_n^m$ are defined before
\ref{mub:isom}. The arithmetical properties \I{n} of almost
independence and \J{n} of block independence are defined in
\ref{mub:def:ar} and \ref{arith:def} respectively. The pairing $\XE$
is defined before \ref{mub:lim}.

\section{Metric unconditional basic sequences of characters \umbs/}
\label{sect:mub}

\subsection{Definitions. Isomorphic case}

We start with the definition of metric 
unconditional basic sequences (\umbs/ for short). 
$\U=\T=\{\epsilon\in\C:|\epsilon|=1\}$ (\vs $\U=\D=\{-1,1\}$) 
is the complex (\vs real) choice of signs.
\begin{dfnsub}\label{mub:def}
  Let \DE/ and $X$ be a homogeneous Banach space on \T.
  \begin{itemize}
  \item [$(i)$]\cite{ka48} $E$ is an unconditional basic
    sequence\index{unconditional basic sequence of characters} \ubs/
    in $X$ if there is a constant $C$ such that
\begin{equation}\label{umbp:dfn}
\biggl\|\sum_{q\in G}\epsilon_qa_q\e_q\biggr\|_X\le
C\biggl\|\sum_{q\in G}a_q\e_q\biggr\|_X
\end{equation}
for all finite subsets $G\se E$, coefficients $a_q\in\C$ and signs
$\epsilon_q\in\T$ \(\vs $\epsilon_q\in\D$\). The infimum of such $C$
is the complex \(\vs real\) unconditionality\index{unconditionality constant}
constant of $E$ in $X$. If $C=1$
works, then $E$ is a complex \(\vs real\)
$1$-\ubs/\index{$1$-unconditional basic sequence of characters} in
$X$.
  \item [$(ii)$]$E$ is a complex \(\vs real\) metric unconditional basic
sequence\index{metric unconditional basic sequence} \umbs/ in $X$ if
for each $\eps>0$ there is a finite set $F$ such that
the complex \(\vs real\) unconditionality constant of $E\setminus F$
is less than $1+\eps$.
\end{itemize}
\end{dfnsub}
Note that \Z\ itself is an \ubs/ in $\SLP{}$ if and only if $p=2$ by
Khinchin's inequality. The same holds in the
framework of the Cantor group\index{Cantor group} $\D^\infty$ and its dual group of Walsh
functions: their common feature with the $\e_n$ is that their modulus
is everywhere equal to $1$ (see \cite{ke81}). 

The following facts are folklore. 
\begin{prpsub}
  Let $Y$ be a Banach space.
  \begin{itemize}
  \item [$(i)$]If $\bigl\|\sum\epsilon_ky_k\bigr\|_Y\le C\bigl\|\sum
    y_k\bigr\|_Y$ for all $\epsilon_k\in\T$ \(\vs $\epsilon_k\in\D$\),
    then this holds automatically for all complex \(\vs real\)
    $\epsilon_k$ with $|\epsilon_k|\le 1$.
  \item [$(ii)$]Real\index{real vs.\ complex}\index{complex vs.\ real} and
  complex unconditionality are isomorphically $\pi/2$-equivalent.
  \end{itemize}
\end{prpsub}
\dem
$(i)$ follows by convexity. $(ii)$ Let us use the fact that the 
complex unconditionality constant of the Rademacher sequence is $\pi/2$ 
\cite{se97}: 
\begin{eqnarray*}
\sup_{\delta_k\in\sT}\Bigl\|\sum\delta_ky_k\Bigr\|_Y
&=&\sup_{y^*\in Y^*}\sup_{\delta_k\in\sT}\sup_{\epsilon_k=\pm1}
\Bigl|\sum\delta_k\langle y^*,y_k\rangle\epsilon_k\Bigr|\\&\le&
\pi/2\sup_{y^*\in Y^*}\sup_{\epsilon_k=\pm1}
\Bigl|\sum\langle y^*,y_k\rangle\epsilon_k\Bigr|=
\pi/2\sup_{\epsilon_k=\pm1}\Bigl\|\sum \epsilon_ky_k\Bigr\|_Y.
\end{eqnarray*}
Taking the Rademacher sequence in $\script{C}(\D^\infty)$, 
we see that $\pi/2$ is optimal.\eck\vskipb

In fact, if $(\ref{umbp:dfn})$ holds, then $E$ is a basis 
of its span in $X$, which is $X_E$
\cite[Th.\ 2.12]{ka68}. 
We have the following relationship between 
the unconditionality 
constants of $E$ in $\SC{}$ and in a homogeneous Banach space $X$ on
$\T$. 
\begin{prpsub}\label{mub:pinf}
  Let \DE/ and $X$ be a homogeneous Banach space on $\T$.
  \begin{itemize}
  \item [$(i)$]The complex \(\vs real\) unconditionality constant of
    $E$ in $X$ is at most the complex \(\vs real\)
unconditionality constant of $E$ in $\SC{}$.
\item [$(ii)$]If $E$ is a \ubs/ \(\vs $1$-\ubs/, \umbs/\) in $\SC{}$,
  then $E$ is a
\ubs/ \(\vs $1$-\ubs/, \umbs/\) in $X$.
\end{itemize}
\end{prpsub}
This follows from the well-known (see {\it e.g.\/} 
\cite{ha87})
\begin{lemsub}\label{sbd:multiplier}\index{relative multipliers}
  Let \DE/ and $X$ be a homogeneous Banach space on $\T$. Let $T$ be a
  multiplier on $\SCE$. Then $T$ is also a multiplier on $X_E$ and
  $$\|T\|_{\script{L}(X_E)}\le \|T\|_{\script{L}(\script{C}_E)}.$$
%
\end{lemsub}
\dem
The linear functional $f\mapsto Tf(0)$ on $\SCE$ extends to a 
measure $\mu\in\script{M}(\T)$ such that 
$\|\mu\|_\script{M}=\|T\|_{\script{L}(\script{C}_E)}$. Let 
$\check{\mu}(t)=\mu(-t)$. Then $Tf=\check{\mu}*f$ for $f\in\PTE/$ and
$$
\|T\|_{\script{L}(X_E)}
\le\|\check{\mu}\|_\script{M}
=\|T\|_{\script{L}(\script{C}_E)}.\eqno{\ecks}$$

\questsub
There is no interpolation 
theorem for such relative multipliers. The forthcoming 
Theorem \ref{mub:thm} shows that there can be no 
metric interpolation\index{relative multipliers!interpolation}. 
Is it possible that one
cannot interpolate 
multipliers at all between $\SLPE$ and $\SLE{q}$~? \vskipb

Note that conversely, \cite{fo82} furnishes the example 
of an \DE/ such that 
the $\pi_k$ are uniformly 
bounded on $\SLE{1}$ but not on $\SCE$.

It is known that $E$ is an \ubs/ 
in $\SC{}$ (\vs in $\SLP{}$) if and only if 
it is a Sidon (\vs \EL{2\vee p}) set. 
To see this, let us recall the relevant definitions.
\begin{dfnsub}\label{mub:sido}
  Let \DE/.
  \begin{itemize}
  \item [$(i)$]\cite{ka57} $E$ is a Sidon\index{Sidon set} set if there is a
  constant $C$ such that
  $$
  \sum_{q\in G}|a_q|\le C\biggl\|\sum_{q\in
    G}a_q\e_q\biggr\|_\infty \mbox{ for all finite $G\se E$ and
    $a_q\in\C$.}
  $$
  The infimum of such $C$ is $E$'s Sidon 
  constant\index{Sidon set!constant}.
  \item [$(ii)$]\cite[Def.\ 1.5]{ru60} Let $p>1$. $E$ is a
  \EL{p}\index{Lambda(p) set@\EL{p} set} set if there is a constant
  $C$ such that $\|f\|_p\le C\|f\|_1$ for $f\in\PTE/$.
\end{itemize}
\end{dfnsub}
In fact, the Sidon constant of $E$ is the complex unconditionality
constant of $E$ in $\SC{}$. Thus $E$ is a complex \umbs/ 
in $\SC{}$ if and only if tails of $E$ 
have their Sidon 
constant arbitrarily close to $1$. We may also 
say: $E$'s Sidon constant is 
asymptotically $1$.

Furthermore, $E$ is a \EL{2\vee p} set if and only if 
$\SLPE=\SLE{2}$. Therefore \EL{2\vee p} sets are 
\ubs/ in $\SLP{}$. Conversely, if $E$ is an 
\ubs/ in $\SLP{}$, then by 
Khinchin's
inequality
$$
\biggl\|\sum_{q\in G} a_q\e_q\biggr\|_p^p\approx
\mbox{average}
\biggl\|\sum_{q\in G}\pm a_q\e_q\biggr\|_p^p\approx
\Bigl(\sum_{q\in G}|a_q|^2\Bigr)^{p/2}=
\biggl\|\sum_{q\in G} a_q\e_q\biggr\|_2^p
$$
for all finite $G\se E$ (see 
\cite[proof of Th.\ 3.1]{ru60}). This shows also that the \EL{p} 
set constant and the unconditionality constant in $\SLP{}$ are connected
{\it via}\/ the constants in Khinchin's inequality; whereas Sidon sets
have their unconditionality constant in $\SLP{}$ uniformly bounded, the
\EL{p} set constant\index{Lambda(p) set@\EL{p} set!constant} 
of infinite sets grows at least like $\sqrt{p}$ 
\cite[Th.\ 3.4]{ru60}.

\subsection{Isometric case: $1$-unconditional basic sequences of
 characters}\label{ss:isom}

The corresponding isometric question: when is $E$ a complex
$1$-\ubs/~? admits a rather easy answer. To 
this end, introduce the following notation for arithmetical relations\index{arithmetical relation}:
let 
$\Alpha_n=\bigl\{\alpha=\{\alpha_p\}_{p\ge1}:
\alpha_p\in\N\ \&\ 
\alpha_1+\alpha_2+\dots=n\bigr\}$. If 
$\alpha\in \Alpha_n$, all but a finite number of the 
$\alpha_p$ vanish and the multinomial number 
$$\ds\bi n\alpha=\frac{n!}{\alpha_1!\alpha_2!\dots}$$
is well defined. Let 
$\Alpha_n^m=\{\alpha\in \Alpha_n: \alpha_p=0\mbox{ for }p>m\}$.
Note that $\Alpha_n^m$ is finite.
We call $E$ $n$-independent\index{independent set of integers} if 
every integer admits at most one representation as the sum of $n$
elements of $E$, up to a permutation. In terms of arithmetical
relations, this yields
$$
\sum\alpha_ip_i=\sum\beta_ip_i\imp\alpha=\beta\mbox{ for
 $\alpha,\beta\in \Alpha_n^m$ and distinct $\lst pm\in E$}.
$$
This notion is studied in \cite{dr75} where it is called
birelation\index{birelation}. In Rudin's\index{Rudin, Walter}
\cite[\S$1.6(b)$]{ru60} 
notation, the number $r_n(E;k)$ of
representations of $k\in\Z$ as a sum of $n$ elements of $E$ is at most
$n!$ for all $k$ if $E$ is $n$-independent (the converse if false).
This may also be expressed in the framework of arithmetical 
relations 
$$\Zeta^m=\{\zeta\in{\Z^*}^m:\sm\zeta m=0\}\quad\&\quad\Zeta_n^m=\{\zeta\in\Zeta^m:|\zeta_1|+\dots+|\zeta_m|\le 2n\}.$$
Note that $\Zeta_n^m$ is finite, and void if $m>2n$. 
Then $E$ is $n$-independent if and only if 
$$
\sum\zeta_ip_i\ne0\quad\hbox{for all }\zeta\in\Zeta_n^m
\hbox{ and distinct }\lst pm\in E.
$$ 
We shall prefer 
to treat arithmetical relations in terms of $\Zeta_n^m$ rather 
than $\Alpha_n^m$. 
\begin{prpsub}\label{mub:isom}
\index{$1$-unconditional basic sequence of characters!in spaces $\SLP{}$, $p$ even} 
\index{$1$-unconditional basic sequence of characters!in $\SC{}$ and $\SLP{}$, $p\notin2\N$} 
  Let \DE/.
  \begin{itemize}
  \item [$(i)$]$E$ is a complex $1$-\ubs/ in $\SLP{}$, $p$ not an even
  integer, or in $\SC{}$, if and only if $E$ has at most two elements.
  \item [$(ii)$]If $p$ is an even integer, then $E$ is a complex $1$-\ubs/
  in $\SLP{}$ if and only if $E$ is $p/2$-independent. There is a
  constant $C_p>1$ depending only on $p$, such that either $E$ is a
  complex $1$-\ubs/ in $\SLP{}$ or the complex unconditionality
  constant of $E$ in $\SLP{}$ is at least $C_p$.
  \end{itemize}
\end{prpsub}
\dem
$(i)$ By Proposition \ref{mub:pinf}$(ii)$, if $E$ is not a complex
$1$-\ubs/ in some $\SLP{}$, then neither in $\SC{}$. 
Let $p$ be not an even integer. We may suppose 
$0\in E$; let $\{0,k,l\}\se E$. If we had
$\|1+\mu a\e_k+\nu b\e_l\|_p=\|1+a\e_k+b\e_l\|_p$ 
for all $\mu,\nu\in\T$, then 
\begin{eqnarray*}
\int|1+a\e_{k}+b\e_{l}|^pdm&=&
\int|1+\mu a\e_{k}+\nu b\e_{l}|^pdm(\mu)dm(\nu)dm\\
&=&
\int|1+\mu a+\nu b|^pdm(\mu)dm(\nu).
\end{eqnarray*}
Denoting by $\theta_i\colon(\epsilon_1,\epsilon_2)\mapsto\epsilon_i$ 
the projections of $\T^2$ onto $\T$, this would mean that
$\|1+a\e_{k}+b\e_{l}\|_p
=\|1+a\theta_1+b\theta_2\|
_{{\fam0 L}^p(\sT^2)}$ 
for all $a,b\in\C$. By 
\cite[Th.\ I]{ru76}, 
$(\e_{k},\e_{l})$ and $(\theta_1,\theta_2)$ 
would have the
same distribution. This is false, since 
$\theta_1$ and $\theta_2$ are independent random variables 
while $\e_{k}$ and $\e_{l}$ are not. 

$(ii)$
Let $\lst qm\in E$ be distinct and 
$\lst\epsilon m\in\T$. 
By the multinomial formula for the power $p/2$ and 
Bessel--Parseval's formula, we get
\begin{eqnarray}
\nonumber
\lefteqn{\Biggl\|\sum_{i=1}^m\epsilon_ia_i\e_{q_i}\Biggr\|_p^p
=\int\Biggl|\sum_{\alpha\in \Alpha_{p/2}^m}
\bip{p/2}\alpha\prod_{i=1}^m(\epsilon_ia_i)^{\alpha_i}
\e_{\lower1pt\hbox{$\Sigma$}\alpha_iq_i}\Biggr|^2dm}\\
\nonumber&=&\sum_{A\in\script{R}_q}\Biggl|
\sum_{\alpha\in A}
\bip{p/2}\alpha
\prod_{i=1}^m(\epsilon_ia_i)^{\alpha_i}\Biggr|^2\\
\nonumber&=&\sum_{\alpha\in \Alpha_{p/2}^m}
{\bip{p/2}\alpha}^2\prod_{i=1}^m|a_i|^{2\alpha_i}
\label{cp>1}+
\sum_{\scriptstyle\alpha\ne\beta\in \Alpha_{p/2}^m\atop
\scriptstyle\alpha\sim\beta}
\bip{p/2}\alpha\bip{p/2}\beta
\prod_{i=1}^m\epsilon_i^{\alpha_i-\beta_i}
a_i^{\alpha_i}\overline{a_i}^{\beta_i},
\end{eqnarray}
where $\script{R}_q$ is the partition of $\Alpha_{p/2}^m$ 
induced by the equivalence relation 
$\alpha\sim\beta\Leftrightarrow
\sum\alpha_iq_i=\sum\beta_iq_i$. If $E$ is 
$p/2$-independent, the second sum in $(\ref{cp>1})$ 
is void and $E$ is a 
$1$-\ubs/.

Furthermore, suppose $E$ is not $p/2$-independent and 
let $\lst qm\in E$ be a minimal set of distinct elements of $E$ 
such that there are $\alpha,\beta\in \Alpha_{p/2}^m$ with 
$\alpha\sim\beta$. Then $m\le p$. Take $a_i=1$ in 
the former computation: then the clearly nonzero 
oscillation of 
$(\ref{cp>1})$ for $\lst\epsilon m\in\T$ does only 
depend on $\script{R}_q$ and thus is finitely valued. 
This yields $C_p$.
\eck\vskipb

\exasub\index{unconditionality constant!in $\SL{4}{}$}
Let us treat explicitly the case $p=4$. If $E$ is not $2$-independent,
then one of the two following arithmetic relations occurs on $E$:
$$
2q_1=q_2+q_3\quad\hbox{or}\quad q_1+q_2=q_3+q_4.
$$
In the first case, we may assume $q_2<q_1<q_3$ and thus
$$
2q_2<q_1+q_2<2q_1=q_2+q_3<q_1+q_3<2q_3.
$$
Let $\varrho>0$. Then
$$
\int{|\e_{q_1}+\varrho\e_{q_2}+\epsilon\varrho\e_{q_3}|}^4dm=
1+6\varrho^4+4\varrho^2(2+\Re\epsilon).
$$
Taking $\epsilon=-1$ and $\epsilon=1$, $\varrho=6^{-1/4}$, we see that
$E$'s real 
unconditionality constant is at least the fourth root of $2\sqrt{6}-3$. In fact, $E$'s real and
complex unconditionality constants coincide with this value.

In the second case, we may assume $q_1<q_3<q_4<q_2$ and thus
$$
2q_1<q_1+q_3<q_1+q_4,2q_3<q_1+q_2=
q_3+q_4<q_2+q_3,2q_4<q_2+q_4<2q_2.
$$
We may further assume $q_1+q_4\ne2q_3$ and $q_2+q_3\ne2q_4$: 
otherwise the first case occurs. Then
$$
\int{|\e_{q_1}+\e_{q_2}+\e_{q_3}+\epsilon\e_{q_4}|}^4dm=28+8\Re\epsilon.
$$
Thus $E$'s real
unconditionality constant must be at least $(9/5)^{1/4}$. In fact,
$E$'s real and 
complex unconditionality constants coincide with this value.

From these two cases we conclude that $C_2=(9/5)^{1/4}\approx1.16$ is the optimal
choice 
for the constant in Proposition \ref{mub:isom}$(ii)$.\vskipb

\remsub
We shall compute explicitly the Sidon constant of sets with three
elements and show that it is equal to the real unconditionality
constant in that case. This provides an alternative proof and a
generalization of Prop.\ \ref{mub:isom} $(i)$ for $\SC{}$.\vskipb

\remsub\edef\remre{\thethmsub}
In fact the conclusion in $(ii)$ holds also 
if we assume that $E$ is just a
real\index{real vs.\ complex}\index{complex vs.\ real} 
$1$-\ubs/. If we have some arithmetical 
relation $\alpha\sim\beta$, we may 
assume that $\alpha_i-\beta_i$ is 
odd for one $i$ at least. 
Indeed, we may simplify all 
$\alpha_i-\beta_i$ by their 
greatest common divisor and this 
yields another arithmetical 
relation $\sum(\alpha'_i-\beta'_i)q_i=0$. 
But then the 
oscillation of 
$(\ref{cp>1})$ is again clearly nonzero for 
$\lst\epsilon m\in\D$. \vskipb

\remsub
We shall see in Remark \ref{umbs:remRE} that $(i)$ also holds 
in the real\index{real vs.\ complex}\index{complex vs.\ real} setting. This is a property of $\T$ 
and fails for the Cantor\index{Cantor group} group $\D^\infty$: 
the Rademacher
sequence forms a real $1$-\ubs/ 
in $\script{C}(\D^\infty)$ but is clearly not complex 
$1$-unconditional in any space ${\fam0 L}^p(\D^\infty)$, $p\ne2$: 
see Section \ref{sect:proba} and \cite{se97}. \vskipb

\questsub
There are nevertheless subspaces of $\SLP{}$, $p$ not an even integer,
and $\SC{}$ with $1$-unconditional
bases\index{$1$-unconditional basic sequence of characters!in $\SC{}$ and $\SLP{}$, $p\notin2\N$}, 
like sequences of functions with disjoint support. What
about spaces $\SLPE$ and $\SCE$, in particular when $E$ is finite ? Are there
$1$-unconditional bases that do not consist of characters ?\vskipb

\remsub
For each even integer $p\ge4$, there are $p/2$-independent sets that
are not 
\EL{p+\eps} for any $\eps>0$: such maximal \EL{p} sets are
constructed in \cite{ru60}.\vskipb

\subsection{Almost isometric case. A computation}
\label{ss:almost}

As $1$-\ubs/ are thus a quite exceptional 
phenomenon and distinguish so harshly 
between even integers and all other reals, one may 
wonder what kind of behaviour its almost isometric counterpart will 
bring about. In the proof 
of Proposition \ref{mub:isom}$(i)$, we used the fact that 
the $\e_n$, seen as random variables, are dependent: 
the ${\fam0 L}^p$ norm for even integer $p$ is just somewhat 
blind to this because it keeps the interaction of 
the random variables down to a finite number of 
arithmetical relations\index{arithmetical relation}. 
The contrast with the other 
${\fam0 L}^p$ norms becomes clear when we try to compute 
explicitly an expression of type 
$\bigl\|\sum\epsilon_qa_q\e_q\bigr\|_p$ for any $p\in[1,\infty\mathclose[$. 
This sort of seemingly brutal computation has been 
applied successfully in \cite[Prop.\ 2]{fo64} and 
\cite[Th.\ 1.4]{pl74} to 
study isometric operators\index{isometries on ${\fam0 L}^p$} on
${\fam0 L}^p$, $p$ not 
an even integer.

We now undertake this tedious computation as preparatory 
work for Theorem \ref{mub:thm}, 
Lemma \ref{umap:lem} and Proposition \ref{umap:prp:fdd}. 
Let us fix some more notation:
for $x\in\R$ and $\alpha\in \Alpha_n$, 
put 
$$\bi x\alpha=\bi xn\bi n\alpha.$$
This generalized 
multinomial number is nonzero if and only if $x\ge n$ or 
$x\notin\N$.
\begin{ldcsub}\label{mub:calcul}
  Let $\U=\T$ or $\U=\D$ in the complex and real case respectively.
  Let $1\le p<\infty$ and $m\ge1$. Put
  $$
  \varphi_q(\epsilon,z,t)=\biggl|1+\sum_{i=1}^m
  \epsilon_iz_i\e_{q_i}(t)\biggr|^p\quad,\quad \Phi_q(\epsilon,z)=
  \int\varphi_q(\epsilon,z,t)\,dm(t)
  $$
  for $q=(\lst qm)\in\Z^m$, $\epsilon=(\lst\epsilon m)\in\U^m$ and
  $z=(\lst zm)\in D^m$, where $D$ is the disc $\{|w|\le\varrho\}\se\C$
  for some $0<\varrho<1/m$.  
Define the equivalence relation $\alpha\sim\beta\Leftrightarrow
\sum\alpha_iq_i=\sum\beta_iq_i$. Then
\begin{eqnarray}
\nonumber\Phi_q(\epsilon,z)&=&
\sum_{\alpha\in\sN^m}
{\bip{p/2}\alpha}^2\prod|z_i|^{2\alpha_i}
\label{mub:phi}+
\sum_{\scriptstyle\alpha\ne\beta\in\sN^m\atop\scriptstyle\alpha\sim\beta}
\bip{p/2}\alpha\bip{p/2}\beta
\prod z_i^{\alpha_i}
\overline{z_i}^{\beta_i}
\epsilon_i^{\alpha_i-\beta_i}.
\end{eqnarray}
Furthermore, $\{\Phi_q:q\in\Z^m\}$ is a relatively compact subset of
$\script{C}^\infty(\U^m\times D^m)$.
\end{ldcsub}
\dem 
The function $\Phi_q$ is infinitely 
differentiable on the compact set 
$\U^m\times D^m$. Furthermore the family 
$\{\Phi_q:\lst qm\in\Z\}$ 
is bounded in 
$\script{C}^\infty(\U^m\times D^m)$ and henceforth 
relatively compact by Montel's theorem.
Let us compute $\varphi_q$. 
By the expansion of the function $(1+w)^{p/2}$, 
analytic on the unit disc, and the multinomial formula, we have
\begin{eqnarray*}
\varphi_q(\epsilon,z)
&=&\biggl|\sum_{a\ge0}\bip{p/2}a
\biggl(\sum_{i=1}^m\epsilon_iz_i\e_{q_i}
\biggr)^a\biggr|^2\\
&=&\biggl|
\sum_{a\ge0}\bip{p/2}a
\sum_{\alpha\in \Alpha_a^m}\bi a\alpha
\prod(\epsilon_iz_i)^{\alpha_i}
\e_{\lower1pt\hbox{$\Sigma$}\alpha_iq_i}
\biggr|^2\\
&=&\biggl|\sum_{\alpha\in\sN^m}
\bip{p/2}\alpha
\prod(\epsilon_iz_i)^{\alpha_i}
\e_{\lower1pt\hbox{$\Sigma$}\alpha_iq_i}\biggr|^2.
\end{eqnarray*}
Let $\script{R}_q$ be the partition of $\N^m$ induced by $\sim$. Then,
by Bessel--Parseval's formula 
$$
\Phi_q(\epsilon,z)=\sum_{A\in\script{R}_q}
\biggl|\sum_{\alpha\in A}
\bip{p/2}\alpha\prod(\epsilon_iz_i)^{\alpha_i}
\biggr|^2
$$
and this gives $(\ref{mub:phi})$ by expanding the modulus.
\eck\vskipb

\remsub\edef\remar{\thethmsub}
If $m\ge2$, this expansion has a finite number of terms if and only if $p$ 
is an even integer: then and only then 
$\bip{p/2}\alpha=0$ for $\sum\alpha_i>p/2$, whereas 
$\script{R}_q$ contains clearly some class with 
two elements and thus an infinity thereof. For example,
we have the following arithmetical 
relation\index{arithmetical relation} on $q_1,q_2$ or 
$q_1,q_2,0$ respectively: 
$$
\overbrace{q_1+\dots+q_1}^{|q_2|}=
\overbrace{q_2+\dots+q_2}^{|q_1|}\qquad\mbox{if }\sgn q_1=\sgn q_2;
$$
$$
\overbrace{q_1+\dots+q_1}^{|q_2|}+\overbrace{q_2+\dots+q_2}^{|q_1|}
=0\qquad\mbox{if not.}
$$

\remsub\label{umbs:remRE}
This shows that Proposition \ref{mub:isom}$(i)$ 
holds also in the real setting: we may suppose that 
$0\in E$; take $m=2$ and choose
$q_1,q_2\in E$. One of the two relations 
in Remark \remar\ yields an arithmetical relation on $E$ 
with at least one odd coefficient, 
as done in Remark \remre. But then
$(\ref{mub:phi})$ contains terms nonconstant in 
$\epsilon_1\in\D$ or in $\epsilon_2\in\D$ and thus $E$ cannot be a real
$1$-unconditional basic sequence in $\SLP{}$.\vskipb

We return to our computation. 
\begin{ldcsub}\label{mub:culcul}
  Let $r=(\lstp r0m)\in E^{m+1}$ and put $q_i=r_i-r_0$ $(1\le i\le
  m)$.  Define
\begin{equation}\label{arith:Theta}
\Theta_r(\epsilon,z)=\int
\biggl|\e_{r_0}+\sum_{i=1}^m\epsilon_iz_i\e_{r_i}\biggr|^p=
\Phi_q(\epsilon,z)
\end{equation}
Let $\lstp\zeta0m\in\Z^*$ and
\begin{equation}\label{mub:gammadelta}
(\gamma_i,\delta_i)=
(-\zeta_i\vee0,\zeta_i\vee0)\qquad(1\le i\le m).
\end{equation} 
If the arithmetical relation
\begin{equation}\label{mub:ArithRel}
\zeta_0r_0+\dots+\zeta_mr_m=0
\qquad\mbox{while}\qquad\smp\zeta0m=0
\end{equation}
holds, then the coefficient of $\prod z_i^{\gamma_i}
\overline{z_i}^{\delta_i} \epsilon_i^{\gamma_i-\delta_i}$ in
$(\ref{mub:phi})$ is $\bip{p/2}\gamma\bip{p/2}\delta$ and thus
independent of $r$. If $\sum|\zeta_i|\le p$ or $p$ is not an even
integer, this coefficient is nonzero.
\end{ldcsub}
\dem We have 
$\delta_i-\gamma_i=\zeta_i$, 
$\sum\gamma_i-\sum\delta_i=\zeta_0$ and 
$\sum\gamma_i+\sum\delta_i=|\zeta_1|+\dots+|\zeta_m|$, 
so that $\sum\gamma_i\vee\sum\delta_i=\frac12\sum|\zeta_i|$. 
Moreover $\sum(\delta_i-\gamma_i)q_i=\sum\zeta_ir_i=0$, 
so that $\gamma\sim\delta$.
\eck\vskipb

\subsection{Almost independent sets of integers. Main theorem}

The Computational lemmas suggest the following definition. 
\begin{dfnsub}\label{mub:def:ar}
  Let \DE/.
  \begin{itemize}
  \item [$(i)$]$E$ enjoys the property \I{n} of almost
  $n$-independence\index{almost independence} provided there is a
  finite subset $G\se E$ such that $E\setminus G$ is $n$-independent,
  \ie $\zeta_1r_1+\dots+\zeta_mr_m\ne0$ for all $\zeta\in\Zeta_n^m$
  and $\lst rm\in E\setminus G$.
  \item [$(ii)$]$E$ enjoys exactly \I{n} if furthermore it fails \I{n+1}.
  \item [$(iii)$]$E$ enjoys \I{\infty} if it enjoys \I{n} for all $n$, \ie for any
  $\zeta\in\Zeta^m$ there is a finite set $G$ such that
  $\zeta_1r_1+\dots+\zeta_mr_m\ne0$ for $\lst rm\in E\setminus G$.
  \end{itemize}
\end{dfnsub}
Note that property \I{1} is void and that $\I{n+1}\imp\I{n}$. 
This property is also stable under unions with a finite set. 
The preceding computations yield 
\begin{thmsub}\label{mub:thm}
  \index{metric unconditional basic sequence} Let \DEE/ and $1\le
  p<\infty$.
  \begin{itemize}
  \item [$(i)$]Suppose $p$ is an even integer. Then $E$ is a real, and at
  the same times complex, \umbs/ in $\SLP{}$ if and only if $E$ enjoys
  \I{p/2}. If \I{p/2} holds, there is in fact a finite $G\se E$ such
  that $E\setminus G$ is a $1$-\ubs/ in $\SLP{}$.
  \item [$(ii)$]If $p$ is not an even integer and $E$ is a real or complex
  \umbs/ in $\SLP{}$, then $E$ enjoys \I{\infty}.
  \end{itemize}
\end{thmsub}
\dem
Sufficiency in $(i)$ follows directly from 
Proposition \ref{mub:isom}: 
if $E\setminus G$ is $p/2$-independent, then 
$E\setminus G$ is a real and complex $1$-\ubs/.

Let us prove the necessity of the arithmetical property.
We keep the notation of Computational lemmas 
\ref{mub:calcul} and \ref{mub:culcul}. 
Assume $E$ fails \I{n} and let 
$\lstp\zeta0m\in\Z^*$ with $\sum\zeta_i=0$ 
and $\sum|\zeta_i|\le2n$ such that for each $l\ge1$ there 
are distinct $\lstp{r^l}0m\in E\setminus\{\lst nl\}$ with 
$\zeta_0r_0^l+\dots+\zeta_mr_m^l=0$. One may 
furthermore assume that at least one of the $\zeta_i$ 
is not even. 

Assume $E$ is a \umbs/ in $\SLP{}$. Then the oscillation of 
$\Theta_r$ in $(\ref{arith:Theta})$ satisfies 
\begin{equation}\label{mub:ApplDef}
\osc_{\epsilon\in\sU^m}
\Theta_{r^l}(\epsilon,z)\tol_{l\to\infty}0
\end{equation}
for each $z\in D^m$. We may assume that 
the sequence of functions $\Theta_{r^l}$ converges 
in $\script{C}^\infty(\U^m\times D^m)$ to a function 
$\Theta$. Then by $(\ref{mub:ApplDef})$,
$\Theta(\epsilon,z)$ is constant in $\epsilon$ 
for each $z\in D^m$: in particular, its coefficient of 
$\prod z_i^{\gamma_i}
\overline{z_i}^{\delta_i}
\epsilon_i^{\gamma_i-\delta_i}$ is zero. 
(Note that at least one of the $\gamma_i-\delta_i$ is not 
even). 
This is impossible by Computational lemma \ref{mub:culcul} 
if $p$ is either not an even integer or if $p\ge2n$.
\eck
\begin{corsub}\label{arith:cor}
  \index{Sidon set!with constant asymptotically $1$} Let \DE/. If $E$
  is a \umbs/ in $\SC{}$, that is $E$'s Sidon constant is
  asymptotically $1$, then $E$ enjoys \I{\infty}. The converse does
  not hold.
\end{corsub}
\dem
Necessity follows from Theorem \ref{mub:thm} and 
Proposition \ref{mub:pinf}$(ii)$. 
There is a counterexample to the converse in 
\cite[Th.\ 4.11]{ru60}: Rudin constructs a set $E$ that 
enjoys \I{\infty} while $E$ is not even a Sidon\index{Sidon set} set.
\eck\vskipb


For $p$ an even integer, Sections \ref{sect:ex} 
and \ref{sect:comb} will provide various examples of 
\umbs/ in $\SLP{}$. 
Proposition \ref{comb:grow} gives a 
general growth condition on $E$ under which it is an \umbs/. 

As we do not know any partial converse to 
Theorem \ref{mub:thm}$(ii)$ and Corollary \ref{arith:cor}, the sole 
known examples 
of \umbs/ in $\SLP{}$, $p$ not an 
even integer, and $\SC{}$ are those given by 
Theorem \ref{positif:thm}. This theorem will 
therefore provide us with Sidon sets of constant 
asymptotically $1$. Note, however, that Li 
\cite[Th.\ 4]{li96} already constructed 
implicitly such a Sidon set by using 
Kronecker's theorem.

\section{Examples of metric unconditional basic 
s\-e\-q\-u\-e\-n\-c\-e\-s}\label{sect:ex}

After a general study of the arithmetical property 
of almost independence \I{n}, we shall investigate 
three classes of subsets of $\Z$: integer 
geometric sequences, more generally integer parts of 
real geometric sequences, and polynomial sequences.

\subsection{General considerations}
\label{ss:mubs:gen}

The quantity
\begin{eqnarray*}
\XE&=&\sup_{G\se E\hbox{\scriptsize\ finite}}\,
\inf\bigl\{|\zeta_1p_1+\dots+\zeta_mp_m|:
\lst p m\in E\setminus G\hbox{ distinct}\bigr\}\\
&=&\lim_{l\to\infty}\inf\bigl\{|\zeta_1p_1+\dots+\zeta_mp_m|:
\lst p m\in\{\lstf nl\}\hbox{ distinct}\bigr\},
\end{eqnarray*}
where $\{n_k\}=E$, plays a key r\^ole. We have
\begin{prpsub}\label{mub:lim}
  Let \DEE/.
  \begin{itemize}
  \item [$(i)$]$E$ enjoys \I{n} if and only if $\XE\ne0$ for all
  $\zeta\in\Zeta_n^m$.  If $\XE<\infty$ for some $\lst\zeta m\in\Z^*$,
  then $E$ fails \I{|\zeta_1|+\dots+|\zeta_m|}. Thus $E$ enjoys
  \I{\infty} if and only if $\XE=\infty$ for all $\lst\zeta m\in\Z^*$.
  \item [$(ii)$]Suppose $E$ is an increasing sequence. If $E$ enjoys \I{2},
  then the pace $n_{k+1}-n_k$ of $E$ tends to infinity.
  \item [$(iii)$]Suppose $jF+s,kF+t\se E$ for an infinite $F$, $j\ne
  k\in\Z^*$ and $s,t\in\Z$. Then $E$ fails \I{|j|+|k|}.
  \item [$(iv)$]Let $E'=\{n_k+m_k\}$ with $\{m_k\}$ bounded. Then
  $\XE=\infty$ if and only if $\X{E'}=\infty$. Thus \I{\infty} is
  stable under bounded perturbations of $E$.
  \end{itemize}
\end{prpsub}
\dem
$(i)$
Suppose $\XE<\infty$. Then there is an $h\in\Z$ such that there are 
sequences $\lst{p^l}m\in\{n_k\}_{k\ge l}$ with 
$\sum\zeta_ip_i^l=h$ and $\{\lst{p^{l+1}}m\}$ is disjoint from
$\{\lst{p^l}m\}$ for all $l\ge1$. As
$\sum\zeta_ip_i^l-\sum\zeta_ip_i^{l+1}=0$ for $l\ge1$,
$E$ fails \I{|\zeta_1|+\dots+|\zeta_m|}.

$(ii)$
Indeed, $\langle (1,-1),E\rangle=\infty$.

$(iii)$
Put $\zeta=(j,-k)$. Then $\XE<\infty$.
\eck\vskipb

\subsection{Geometric sequences}\index{geometric sequences}

Let $G=\{j^k\}_{k\ge0}$ with
$j\in\Z\setminus\{-1,0,1\}$.
Then $G,jG\se G$: so $G$
fails \I{|j|+1}. 
In order to check \I{|j|} for $G$, 
let us study more carefully the following Diophantine equation:
\begin{equation}\label{arith:dioph}
%
%
\sum_{i=1}^m\zeta_ij^{k_i}=0\quad\hbox{with}\quad
\zeta\in\N^*\times{\Z^*}^{m-1}\
\&\ \sum_{i=1}^m|\zeta_i|\le2|j|\ \&\ k_1<\dots<k_m.
\end{equation}
Suppose $(\ref{arith:dioph})$ holds.
Then necessarily $m\ge2$ and
$\zeta_1+\sum_{i=2}^m\zeta_ij^{k_i-k_1}=0$. Hence
$j\mid\zeta_1$ and $\zeta_1\ge|j|$. As $\zeta_1<2|j|$, $\zeta_1=|j|$. Then 
$\sgn j+\sum_{i=2}^m\zeta_ij^{k_i-k_1-1}=0$. Hence
$k_2=k_1+1$ and $j\mid\sgn j+\zeta_2$. As $|\zeta_2|\le|j|$,
$\zeta_2\in\{-\sgn j,j-\sgn j\}$. If $\zeta_2=j-\sgn j$, then
$m=3$, $k_3=k_1+2$ and $\zeta_3=-1$. If $\zeta_2=-\sgn j$,
then $m=2$: otherwise, $j\mid\zeta_3$ as before and
$|\zeta_1|+|\zeta_2|+|\zeta_3|>2|j|$. Thus
$(\ref{arith:dioph})$ has exactly two solutions:
\begin{equation}\label{arith:sol}
|j|\cdot j^k+(-\sgn j)\cdot j^{k+1}=0\ \&\ 
|j|\cdot j^k+(j-\sgn j)\cdot j^{k+1}+
(-1)\cdot j^{k+2}=0.
\end{equation}
If $j$ is positive, this shows that $G$ enjoys \I{j}: 
both solutions yield $\sum\zeta_i\ne0$. If $j$
is negative, $G$ enjoys \I{|j|-1}, but the second
solution of $(\ref{arith:dioph})$ shows that $G$
fails \I{|j|}.\vskipb

\subsection{Algebraic and transcendental numbers}
\index{transcendental numbers}

An interesting feature of property \I{\infty}
is that it distinguishes between algebraic and
transcendental numbers. A similar fact has
already been noticed by Murai\index{Murai, Takafumi}
\cite[Prop.\ 26, Cor.\ 28]{mu82}.
\begin{prpsub}\label{mub:trans}
  Let \DEE/.
  \begin{itemize}
  \item [$(i)$]If $n_{k+1}/n_k\to\sigma$ where $\sigma>1$ is transcendental,
  then $\XE=\infty$ for any $\lst\zeta m\in\Z^*$. Thus $E$ enjoys
  \I{\infty}.
  \item [$(ii)$]Write $[x]$ for the integer part of a real $x$.  Let
  $n_k=[\sigma^k]$ with $\sigma>1$ algebraic. Let
  $P(x)=\zeta_0+\dots+\zeta_dx^d$ be the corresponding polynomial of
  minimal degree. Then $\XE<\infty$ and $E$ fails
  \I{|\zeta_0|+\dots+|\zeta_d|}.
  \end{itemize}
\end{prpsub}
Note that part $(ii)$ is very restrictive on the speed of
convergence of $n_{k+1}/n_k$ to $\sigma$: even if we take into account
Proposition \ref{mub:lim}$(iv)$, it requires that
$$|n_{k+1}/n_k-\sigma|\preccurlyeq\sigma^{-k}.$$

\dem
$(i)$
Suppose on the contrary that we have $\zeta$ and
sequences $p_1^l<\dots<p_m^l$ in $E$ that tend to
infinity such that $\zeta_1p_1^l+\dots+\zeta_mp_m^l=0$.
As the sequences $\{p_i^l/p_m^l\}_l$ $(1\le i\le m)$
are bounded, we
may assume they are converging ---~and by hypothesis,
they converge either to $0$, say for $i<j$, or to
$\sigma^{-d_i}$ for $d_i\in\N$ and $i\ge j$. But then
$\zeta_j\sigma^{-d_j}+\dots+\zeta_m\sigma^{-d_m}=0$
and $\sigma$ is algebraic.

$(ii)$ Apply Proposition \ref{mub:lim}$(i)$ with
$\zeta$: 
$$
|\zeta_0[\sigma^k]+\dots+\zeta_d[\sigma^{k+d}]|=
|\zeta_0([\sigma^k]-\sigma^k)+\dots
+\zeta_d([\sigma^{k+d}]-\sigma^{k+d})|
\le\sum|\zeta_i|.\eqno{\eck}
$$

\subsection{Polynomial sequences}\index{polynomial sequences}

Let us first give some numerical evidence for the classical case of
sets of $d$th powers. The table below reads as 
follows: ``the set $E=\{k^d\}$ for $d$ the value in the first column
fails the property in the second column by the counterexample given in
the third column.'' Indeed, each such counterexample to $n$-independence
yields arbitrarily large counterexamples.\vskipb

\hfill\vbox{\offinterlineskip\halign{
\vrule\hbox to 2pt{}#\hfil\hbox to 1pt{}\vrule&\strut\lower4pt\vbox to 15pt{}
\hfil#\hfil\vrule&\ #\hfil\hbox to 1pt{}\vrule\cr
\noalign{\hrule}
$\{k^d\}$&\ fails\ &by counterexample\cr
\noalign{\hrule}
d=2\index{polynomial sequences!squares}&\I{2}& $7^2+1^2=2\cdot5^2$ 
(or $18^2+1^2=15^2+10^2$ \cite[book II, problem 9]{di59})\cr
d=3\index{polynomial sequences!cubes}&\I{2}& $12^3+1^3=10^3+9^3$ 
\cite[due to Fr\'enicle]{br96}\index{Fr\'enicle de Bessy, Bernard}\cr
d=4\index{polynomial sequences!biquadrates}&\I{2}& 
$158^4+59^4=134^4+133^4$ 
(or $12231^4+2903^4=10381^4+10203^4$ \cite{eu72})\cr
d=5&\I{3}& $67^5+28^5+24^5=62^5+54^5+3^5$ 
(another first in \cite{mo39})\cr
d=6&\I{3}& $23^6+15^6+10^6=22^6+19^6+3^6$ \cite{ra34}\cr
d=7&\I{4}& $149^7+123^7+14^7+10^7=146^7+129^7+90^7+15^7$ \cite{ek96}\cr
d=8&\I{5}& $43^8+20^8+11^8+10^8+1^8=41^8+35^8+32^8+28^8+5^8$ (see \cite{ek98})\cr
d=9&\I{6}& $23^9+18^9+14^9+2\cdot13^9+1^9=22^9+21^9+15^9+10^9+9^9+5^9$
\cite{lps67}
\cr
d=10&\I{7}& $38^{10}+33^{10}+2\cdot26^{10}+15^{10}+8^{10}+1^{10}=$\cr
&&\hfill$36^{10}+35^{10}+32^{10}+29^{10}+24^{10}+23^{10}+22^{10}$
(another first in \cite{mo39})\cr 
\noalign{\hrule}
}}\hfill\hbox{}\nopagebreak\smallskip\\\nopagebreak
\hfill\refstepcounter{thmsub}Table \thethmsub\hfill\hbox{}

\vskipa 
Note that a positive answer to Euler's\index{Euler's conjecture}
conjecture
---~for $k\ge5$ $a^k+b^k=c^k+d^k$ has only
trivial solutions
in integers~--- would imply that
the set of $k$th powers has \I{2}. This
conjecture has been neither proved nor disproved for any
value of $k\ge 5$ (see \cite{sw95} and \cite{ek98}).

Now let \DEE/ be a set of polynomial growth: $|n_k|\preccurlyeq k^d$ for
some $d\ge1$. Then $\mes{E\cap[-n,n]}\succcurlyeq n^{1/d}$ and by
\cite[Th.\ 3.6]{ru60}, $E$ fails the \EL{p} property for $p>2d$
and $E$ fails {\it a fortiori}\/ \I{d+1}. In the special case
$E=\{P(k)\}$ for a polynomial $P$ of degree $d$, we can exhibit a huge
explicit arithmetical relation. 
%
%
%
%
Recall that
\begin{equation}\label{mub:poly}
\Delta^jP(k)=\suml_{i=0}^j\bi ji(-1)^iP(k-i)
\ ,\ \suml_{i=0}^j\bi ji(-1)^i=0
\ ,\ \suml_{i=0}^j\bi ji=2^j.
\end{equation}
As $\Delta^{d+1}P(k)=0$, this makes $E$ fail \I{2^d}, which is coarse. 
\vskipb

{\bf Conclusion }
By Theorem \ref{mub:thm},
property \I{n} yields directly \umbs/ in the spaces $\SL{2p}{}$, $p\le n$ 
integer.
But we do not know whether
\I{\infty} ensures \umbs/ in spaces
$\SLP{}$, $p$ not an even integer.

\section[Metric unconditional approximation property \umap/]
        {Metric unconditional approximation property}
\label{sect:block}

As we investigate 
simultaneously real and complex \umap/, 
it is convenient to introduce a 
subgroup $\U$ of $\T$ corresponding to each case. Thus, 
if $\U=\D=\{-1,1\}$, then 
the following applies to real \umap/. 
If $\U=\T=\{\epsilon\in\C:|\epsilon|=1\}$, it applies to complex \umap/.

He who is first and foremost interested in the 
application to harmonic analysis may concentrate on the equivalence 
$(ii)\Leftrightarrow(iv)$ in Theorem \ref{block:thm} and then 
pass on to Section \ref{sect:tis}.

\subsection{Definition}

We start with defining the metric unconditional 
approximation property (\umap/ for short). Recall that 
$\Delta T_k=T_k-T_{k-1}$ (where $T_0=0$).
\begin{dfnsub}\label{block:def}\index{approximation property}
  Let $X$ be a separable Banach space.
  \begin{itemize}
  \item [$(i)$]A sequence $\{T_k\}$ of operators on $X$ is an approximating
  sequence\index{approximating sequence} \(\as/ for short\) if each
  $T_k$ has finite rank and $\|T_kx-x\|\to0$ for every $x\in X$.  If
  $X$ admits an \as/, it has the bounded approximation property.  An
  \as/ of commuting projections is a finite-dimensional
  decomposition\index{finite-dimensional decomposition} \(\fdd/ for
  short\).
  \item [$(ii)$]\cite{fe80} $X$ has the unconditional approximation property
  \index{unconditional approximation property} \uap/ if there are an
  \as/ $\{T_k\}$ and a constant $C$ such that
\begin{equation}\label{block:uapdef}
\biggl\|\sum_{k=1}^n\epsilon_k\Delta T_k\biggr\|
\le C\quad\hbox{for all $n$ and $\epsilon_k\in\U$}.
\end{equation}
The \uap/ constant is the least such $C$.
  \item [$(iii)$]\cite[\S3]{ck91} $X$ has the metric unconditional
approximation 
property\index{metric unconditional approximation property} \umap/ if it has \uap/ with constant $1+\eps$ for any
$\eps>0$.
\end{itemize}
\end{dfnsub}
Property $(ii)$ is the approximation property which most 
appropriately generalizes the 
unconditional basis property. It has first been 
introduced by Pe\l czy\'nski \index{Pe\l czy\'nski, Aleksander} and 
Wojtaszczyk\index{Wojtaszczyk, Przemys\l aw} \cite{pw71}. 
They showed that it holds if and only if $X$ is a complemented 
subspace of a space with an unconditional 
\fdd/\index{unconditional fdd@unconditional \fdd/}. 
By 
\cite[Th.\ 1.g.5]{lt77}, 
this implies that $X$ is subspace of a space with an 
unconditional basis. Thus, neither ${\fam0 L}^1([0,1])$ nor 
$\script{C}([0,1])$ share \uap/.
 
Property $(iii)$ has been introduced by 
Casazza\index{Casazza, Peter G.} and 
Kalton\index{Kalton, Nigel J.} as an 
extreme form of metric approximation. It has been studied in 
\cite[\S3]{ck91}, 
\cite[\S8,9]{gks93}, \cite{gkl96} and \cite[\S IV]{gk95}.

There is a simple and very useful criterion for \umap/: 
\begin{prpsub}
  [\protect{\cite[Th.\ 3.8]{ck91} and \cite[Lemma 8.1]{gks93}}]
\label{block:umap:eq}
Let $X$ be a separable Banach space.  $X$ has \umap/ if and only if
there is an \as/ $\{T_k\}$ such that
\begin{equation}\label{block:umap:def}
\sup_{\epsilon\in\sU}\|(\Id-T_k)+\epsilon T_k\|\tol_{k\to\infty}1.
\end{equation}
\end{prpsub}
If \Ref{block:umap:def} holds, we say that $\{T_k\}$ realizes 
\umap/. A careful reading of the above mentioned proof 
also gives 
the following results for \as/ 
that satisfy $T_{n+1}T_n=T_n$. 
\begin{prpsub}
  Let $X$ be a separable Banach space.
  \begin{itemize}
  \item [$(i)$]Let $\{T_k\}$ be an \as/ for $X$ such that $T_{n+1}T_n=T_n$.
  A subsequence $\{T'_k\}$ of $\{T_k\}$ realizes $1$-\uap/ in $X$ if
  and only if for all $k\ge1$ and $\epsilon\in\U$
  $$
    \|\Id-(1+\epsilon)T'_k\|=1.
  $$
  \item [$(ii)$]$X$ has metric unconditional 
\fdd/\index{metric unconditional fdd@metric unconditional \fdd/} if and only if there
  is an \fdd/ $\{T_k\}$ such that $(\ref{block:umap:def})$ holds.
  \end{itemize}
\end{prpsub}

\subsection{Characterization of \umap/. Block unconditionality}

We want to characterize \umap/ in an even simpler 
way than Proposition \ref{block:umap:eq}. 
Relation $(\ref{block:umap:def})$ and 
the method of \cite[Th.\ 4.2]{kw95}, suggest 
considering some unconditionality condition between 
a certain ``break'' and a certain ``tail'' of $X$. We 
propose two such notions. 
\begin{dfnsub}\label{block:def:u}
  Let $X$ be a separable Banach space.
  \begin{itemize}
  \item [$(i)$]Let $\tau$ be a vector space topology on $X$.  Then $X$ has
  the property $(u(\tau))$ of
  $\tau$-unconditionality\index{tau-unconditionality@$\tau$-unconditionality}
  if for all $u\in X$ and norm bounded sequences $\{v_j\}\se X$ such
  that $v_j\hoch{\to}^{\tau}0$
\begin{equation}\label{block:def:u:1}
\osc_{\epsilon\in\sU}\|\epsilon u+v_j\|\to0.
\end{equation}
  \item [$(ii)$]Let $\{T_k\}$ be a commuting a.s.  $X$ has the property
$(u(T_k))$ of commuting block 
unconditionality\index{commuting block unconditionality} if for all $\eps>0$ and $n\ge1$ we may choose
$m\ge n$ such that for all $x\in T_nB_X$ and $y\in (\Id-T_m)B_X$
\begin{equation}\label{block:def:u:2}
\osc_{\epsilon\in\sU}\|\epsilon x+y\|\le\eps.
\end{equation}
  \end{itemize}
\end{dfnsub}
Thus, given a commuting \as/ 
$\{T_k\}$, $T_nX$ 
is the ``break''\index{break} and $(\Id-T_m)X$ the
``tail''\index{tail} of $X$. 
We have 
\begin{lemsub}\label{block:lm}
  Let $X$ be a separable Banach space and $\{T_k\}$ a commuting \as/
  for $X$. The following are equivalent.
  \begin{itemize}
  \item [$(i)$]$X$ enjoys $(u(\tau))$ for some vector space topology $\tau$
  such that $T_nx\hoch{\to}^{\tau}x$ uniformly for $x\in B_X$;
  \item [$(ii)$]$X$ enjoys $(u(T_k))$.
  \end{itemize}
\end{lemsub}
\dem
$(i)\imp(ii)$. 
Suppose that $(ii)$ fails: there are 
$n\ge1$ and $\eps>0$ such that for each $m\ge n$, there 
are $x_m\in T_nB_X$
and $y_m\in(\Id-T_m)B_X$ 
such that 
$$
\osc_{\epsilon\in\sU}\|\epsilon x_m+y_m\|>\eps.
$$
As $T_nB_X$ is compact, we may suppose by extracting 
a convergent subsequence that $x_m=x$. 
Let $\tau$ be as in $(i)$: then $y_m\hoch{\to}^{\tau}0$ 
and $(u(\tau))$ must fail.

$(ii)\imp(i)$. 
Let us define a vector space topology $\tau$ 
by 
$$
x_n\hoch{\to}^\tau0\quad\Longleftrightarrow\quad
\forall k\ \|T_kx_n\|\to0.
$$
Then $T_nx\hoch{\to}^\tau x$ uniformly on $B_X$. 
Indeed, $T_k(T_nx-x)=(T_n-\Id)T_kx$ and $T_n-\Id$ converges 
uniformly to 0 on $T_kB_X$ which is norm compact.

Let us 
check $(u(\tau))$. Let $u\in B_X$ and 
$\{v_j\}\se B_X$ be such 
that $v_j\hoch{\to}^{\tau}0$. 
Let $\eps>0$. There is $n\ge1$ such that 
$\|T_nu-u\|\le\eps$. 
Choose $m$ such that $(\ref{block:def:u:2})$ holds 
for $x\in T_nB_X$ and $y\in (\Id-T_m)B_X$.
Then choose $k\ge1$ such that 
$\|T_mv_j\|\le\eps$ for 
$j\ge k$. We have, for any $\epsilon\in\U$,
\begin{eqnarray*}
\|\epsilon u+v_j\|&\le&
\|\epsilon T_nu+(\Id-T_m)v_j\|+\|T_nu-u\|+\|T_mv_j\|\\
&\le&\|T_nu+(\Id-T_m)v_j\|+3\eps\le\|u+v_j\|+5\eps.
\end{eqnarray*}
Thus we have $(\ref{block:def:u:1})$.
\eck

In order to obtain \umap/ from block independence, we shall
have to construct unconditional 
skipped blocking decompositions. 
\begin{dfnsub}\label{block:sbddef}
  Let $X$ be a separable Banach space. $X$ admits unconditional
  skipped blocking 
decompositions\index{unconditional skipped blocking decompositions} if for each $\eps>0$, there is an unconditional
  \as/ $\{S_k\}$ such that for all $0\le a_1<b_1<a_2<b_2<\dots$ and
  $x_k\in(S_{b_k}-S_{a_k})X$
  $$
  \sup_{\epsilon_k\in\sU}\Bigl\|\sum
  \epsilon_kx_k\Bigr\|\le(1+\eps)\Bigl\|\sum x_k\Bigr\|.
  $$
\end{dfnsub}

\subsection{Main theorem: convex combinations of multipliers}

We have 
\begin{thmsub}\label{block:thm}
  \index{metric unconditional approximation property} Consider the
  following properties for a separable Banach space $X$.
  \begin{itemize}
  \item [$(i)$]There are an unconditional commuting \as/ $\{T_k\}$ and a
  vector space topology $\tau$ such that $X$ enjoys $(u(\tau))$ and
  $T_kx\hoch{\to}^{\tau}x$ uniformly for $x\in B_X$;
  \item [$(ii)$]$X$ enjoys $(u(T_k))$ for an unconditional commuting \as/
  $\{T_k\}$;
  \item [$(iii)$]$X$ admits unconditional skipped blocking decompositions;
  \item [$(iv)$]$X$ has \umap/.
  \end{itemize}
  Then $(iv)\imp(i)\Leftrightarrow(ii)\Rightarrow(iii)$.  If $X$ has
  finite cotype, then $(iii)\imp(iv)$.
\end{thmsub}
\dem
$(i)\Leftrightarrow(ii)$ holds by Lemma \ref{block:lm}.

$(iv)\imp(ii)$. By 
Godefroy\index{Godefroy, Gilles}--Kalton's\index{Kalton, Nigel J.} 
\cite[Th.\ IV.1]{gk95}, there 
is in fact an \as/ $\{T_k\}$ that satisfies 
$(\ref{block:umap:def})$ such that 
$T_kT_l=T_{\min(k,l)}$ if $k\ne l$. 

Let $C$ be a uniform bound for $\|T_k\|$. 
Let $\eps>0$ and $n\ge1$. There is $m\ge n+2$ such that 
$$
\sup_{\epsilon\in\sU}
\|\epsilon T_{m-1}+(\Id-T_{m-1})\|\le1+\eps/2C.
$$ 
Let $x\in T_nB_X$ and $y\in(\Id-T_m)B_X$. As $x-T_{m-1}x=0$ and 
$T_{m-1}y=0$,
$$
\epsilon x+y=\epsilon T_{m-1}(x+y)+
(\Id-T_{m-1})(x+y),
$$
and, for all $\epsilon\in\U$, 
$$
\|\epsilon x+y\|\le(1+\eps/2C)\|x+y\|\le\|x+y\|+\eps.
$$

$(ii)\imp(iii)$. 
By a perturbation \cite[proof of Lemma III.9.2]{si81}, we may
suppose that $T_kT_l=T_{\min(k,l)}$ if $k\ne l$. 
Let $\eps>0$ and choose a sequence of $\eta_j>0$ 
such that $1+\eps_j=\prod_{i\le j}(1+\eta_i)<1+\eps$ 
for all $j$. 
By $(ii)$, there is a subsequence $\{S_j=T_{k_j}\}$
such that $k_0=0$ and thus $S_0=0$, and 
\begin{equation}\label{block:h}
\sup_{\epsilon\in\sU}\|x+\epsilon y\|\le(1+\eta_j)\|x+y\|
\end{equation}
for $x\in(\Id -S_j)X$ and $y\in S_{j-1}X$. Let us show that 
it is an unconditional skipped blocking decomposition: we 
shall prove by induction that
$$
\left\{\!\!\begin{array}{l}
\ds\sup_{\epsilon_i\in\sU}
\biggl\|x+\sum_{i=1}^n
\epsilon_ix_i\biggr\|\le(1+\eps_j)
\biggl\|x+\sum_{i=1}^nx_i\biggr\|\mbox{ for $x\in(\Id -S_j)X$}\\
\mbox{and $x_i\in(S_{b_i}-S_{a_i})X$ 
$(0\le a_1<b_1<\dots<a_n<b_n\le j-1)$.}
\end{array}\right.\leqno{(H_j)}
$$

\bloc $(H_1)$ trivially holds.

\bloc 
Assume $(H_i)$ holds for $i<j$. Let $x$ and $x_i$ as in $(H_j)$. 
Let $\epsilon_i\in\U$. Then
$$
\Bigl\|x+\sum_{i=1}^n\epsilon_ix_i\Bigr\|\le
(1+\eta_j)\Bigl\|x+\overline{\epsilon_n}
\sum_{i=1}^n\epsilon_ix_i\Bigr\|=
(1+\eta_j)\Bigl\|x+x_n+
\sum_{i=1}^{n-1}\overline{\epsilon_n}\epsilon_ix_i\Bigr\|
$$
by $(\ref{block:h})$. Note that $x+x_n\in(\Id -S_{a_n})X$: 
an application of $(H_{a_n})$ yields $(H_j)$.

$(iii)\imp(iv)$. 
Let $\eps>0$, $n>1$. 
There is an unconditional skipped blocking 
decomposition $\{S_k\}$. Let $C_u$ be the \uap/ constant 
of $\{S_k\}$. Let 
$$
V_{i,j}=S_{in+j-1}-S_{(i-1)n+j}\qquad
\mbox{for $1\le j\le n$ and $i\ge 0$.}
$$
The $j$th skipped blocks are 
$$
U_j=\Id -\sum_i V_{i,j}=\sum_i\Delta S_{in+j};
$$
then 
$\sum_{j=1}^n U_j=\Id $. Let 
$$
R_i=\frac1{n-1}\sum_{j=1}^n V_{i,j};
$$
then $R_i$ has finite rank and 
$$
R_0+R_1+\dots=(n\Id -\Id )/(n-1)=\Id .
$$
Thus $W_j=\sum_{i\le j}R_i$ defines an a.s.
We may bound its \uap/ constant. First, since 
$\{S_k\}$ is a skipped blocking decomposition,
\begin{eqnarray*}
\forall x\in B_X\quad\sup_{\epsilon_i\in\sU}
\Bigl\|\sum\epsilon_iR_ix\Bigr\|
&\le&\frac1{n-1}\sum_{j=1}^n\sup_{\epsilon_i\in\sU}
\Bigl\|\sum_i\epsilon_iV_{i,j}x\Bigr\|\\
&\le&\frac{1+\eps}{n-1}\sum_{j=1}^n\|x-U_jx\|\\
&\le&\frac{1+\eps}{n-1}\Bigl(n+\sum_{j=1}^n\|U_jx\|\Bigr).
\end{eqnarray*}
Let us bound $\sum_1^n\|U_jx\|$. 
Let $q<\infty$ be the cotype of $X$ and 
$C_c$ its cotype constant. 
Then by H\"older's inequality we have for all $x\in B_X$
\begin{eqnarray}\label{block:cotype}
\nonumber\sum\|U_jx\|
&\le&n^{1-1/q}
\Bigl(\sum
\|U_jx\|^q\Bigr)^{1/q}\\
&\le&n^{1-1/q}C_c\cdot\mathop{\hbox{average}}\limits_{\pm}
\Bigl\|\sum\pm U_jx\Bigr\|
\le n^{1-1/q}C_cC_u.
\end{eqnarray}
Thus the \uap/ constant of 
$\{W_j\}$ is at most
$(1+\eps)(n+C_cC_un^{1-1/q})/(n-1)$.
As $\eps$ is arbitrarily little and $n$ arbitrarily 
large, $X$ has \umap/.\eck\vskipb

\remsub
How does Theorem \ref{block:thm} look in the special cases where $\tau$ is
the weak or the weak$^*$ topology ? They correspond to the classical
cases where the \as/ is 
shrinking\index{shrinking approximating sequence} \vs boundedly 
complete\index{boundedly complete approximating sequence}.\vskipb
 
We may remove the cotype\index{cotype} assumption in 
Theorem \ref{block:thm} 
$(iii)\imp(iv)$ if the space has the properties of 
commuting \lap{1} or $\ell_q$-\fdd/ for $q<\infty$, 
which will be introduced in Section \ref{sect:lpmap}:
\begin{thmsub}\label{block:thm'}
  Consider the following properties for a separable Banach space $X$.
  \begin{itemize}
  \item [$(i)$]There are a commuting $\ell_1$-\as/ or an $\ell_q$-\fdd/
  $\{T_k\}$, $q<\infty$, and a vector space topology $\tau$ such that
  $X$ enjoys $(u(\tau))$ and $T_kx\hoch{\to}^{\tau}x$ uniformly for
  $x\in B_X$;
  \item [$(ii)$]$X$ enjoys $(u(T_k))$ for a commmuting $\ell_1$-\as/ or an
  $\ell_q$-\fdd/ $\{T_k\}$, $q<\infty$;
  \item [$(iii)$]$X$ admits unconditional skipped blocking decompositions
  and one may in fact take an $\ell_1$-\as/ or an $\ell_q$-\fdd/
  $\{T_k\}$, $q<\infty$, in its definition \ref{block:sbddef};
  \item [$(iv)$]$X$ has \umap/.
  \end{itemize}  
  Then $(i)\Leftrightarrow(ii)\imp(iii)\imp(iv)$.
\end{thmsub}
\dem
Part $(i)\Leftrightarrow(ii)\imp(iii)$ goes as before. 
To prove $(iii)\imp(iv)$, note that in the proof of 
Theorem \ref{block:thm} 
$(iii)\imp(iv)$, one may replace 
the estimate in $(\ref{block:cotype})$ by 
$$
\forall x\in B_X\quad\sum\|U_jx\|\le 
n^{1-1/q}\Bigl(\sum
\|U_jx\|^q\Bigr)^{1/q}\le n^{1-1/q}C_{\ell},
$$
where $C_{\ell}$ is the $\ell_1$-\ap/ or
the $\ell_q$-\fdd/ constant.
\eck

\section{The $p$-additive approximation property \lpap/}
\label{sect:lpmap}

\subsection{Definition}

\begin{dfnsub}\label{str:def}
  Let $X$ be a separable Banach space.
  \begin{itemize}
  \item [$(i)$]$X$ has the $p$-additive approximation property
  \lpap/\index{p-additive approximation property@$p$-additive approximation property} if there are an \as/ $\{T_k\}$ and a
  constant $C$ such that
\begin{equation}\label{block:def:str}
C^{-1}\|x\|\le
\Bigl(\sum\|\Delta T_kx\|^p\Bigr)^{1/p}\le C\|x\|
\end{equation}
for all $x\in X$. The \lpap/ constant is the least such $C$.
  \item [$(ii)$]$X$ has the metric $p$-additive
approximation 
property\index{metric $p$-additive approximation property} 
\lpmap/ if it has \lpap/ with constant $1+\eps$ for any
$\eps>0$.

  \end{itemize}
\end{dfnsub}
Note that \lpap/ implies \uap/ and
\lpmap/ implies \umap/. 
Note also that in $(\ref{block:def:str})$,
the left inequality is trivial with $C=1$ if $p=1$; 
the right inequality is always achieved for some $C$ if $p=\infty$. 

Property $(ii)$ is implicit in 
Kalton\index{Kalton, Nigel J.}--Werner's\index{Werner, Dirk}
\cite{kw95} investigation of subspaces of ${\fam0 L}^p$ that
are almost isometric to subspaces of $\ell_p$: 
see Section \ref{sect:appen}.

The proof of Proposition \ref{block:umap:eq} 
can be adapted to yield
\begin{prpsub}\label{block:str:eq}
  Let $X$ be a separable Banach space.
  \begin{itemize}
  \item [$(i)$]If there is an \as/ $\{T_k\}$ such that
\begin{equation}\label{block:umapf:def}
\Bigl(\|x-T_kx\|^p+\|T_kx\|^p\Bigr)^{1/p}
\tol_{k\to\infty}1
\end{equation}
uniformly on the unit sphere, then $X$ has \lpmap/.  The converse
holds if $p=1$.
  \item [$(ii)$]$X$ has a metric $\ell_p$-\fdd/ if and only if there is an
\fdd/ $\{T_k\}$ such that $(\ref{block:umapf:def})$ holds.
  \end{itemize}
\end{prpsub}
We shall say that
$\{T_k\}$ realizes \lpmap/
if it satisfies
$(\ref{block:umapf:def})$.\vskipb

\dem
Let $\{T_k\}$ be an \as/ that
satisfies $(\ref{block:umapf:def})$ and $\eps>0$.
By a perturbation \cite[Lemma 2.4]{jrz71},
we may suppose that $T_{k+1}T_k=T_k$.
Choose a sequence of $\eta_j>0$ such that
$1+\eps_k=\prod_{j\le k}(1+\eta_j)\le1+\eps$
for each $k$.
We may assume by taking a subsequence of the $T_k$'s that
for all $k$ and $x\in X$,
\begin{equation}\label{block:equ:dem}
(1+\eta_k)^{-1}\|x\|\le
\Bigl(\|x-T_kx\|^p+\|T_kx\|^p\Bigr)^{1/p}
\le(1+\eta_k)\|x\|.
\end{equation}
We then prove by induction the hypothesis $(H_k)$
$$
\forall x\in X\quad
(1+\eps_k)^{-1}\|x\|\le
\Bigl(\|x-T_kx\|^p+
\sum_{j=1}^k\|\Delta T_jx\|^p\Bigr)^{1/p}
\le(1+\eps_k)\|x\|.
$$

\bloc
$(H_1)$ is true.

\bloc
Suppose $(H_{k-1})$ is true. Let $x\in X$. Note
that
$$
x-T_kx=(\Id -T_k)(x-T_{k-1}x)\quad,\quad
\Delta T_kx=T_k(x-T_{k-1}x).
$$
By $(\ref{block:equ:dem})$, we get
\begin{eqnarray*}
\bigl(\|x-T_kx\|^p+\|\Delta T_kx\|^p\bigr)^{1/p}
\le(1+\eta_k)\|x-T_{k-1}x\|.
\end{eqnarray*}
Hence, by $(H_{k-1})$,
\begin{eqnarray*}
\lefteqn{\Bigl(
\|x-T_kx\|^p+\sum_{j=1}^k\|\Delta T_jx\|^p
\Bigr)^{1/p}\le}\quad\quad\quad&&\\
\quad\quad&\le&(1+\eta_k)\Bigl(\|x-T_{k-1}x\|^p+
\sum_{j=1}^{k-1}\|\Delta T_jx\|^p\Bigr)^{1/p}
\le(1+\eps_k)\|x\|.
\end{eqnarray*}

\bloc
We obtain the lower bound in the same way.
Thus the induction is complete.

Hence $\{T_k\}$ realizes
\lpap/ with constant $1+\eps$. As $\eps$ is arbitrary, 
$X$ has \lpmap/.

If $X$ has \lmap{1}, then for each $\eps>0$,
there is a sequence $\{S_k\}$ such that
$$
\|x\|\le\|x-S_kx\|+\|S_kx\|\le
\sum\|\Delta S_kx\|\le(1+\eps)\|x\|
$$
for all $x\in X$. By a diagonal
argument, this gives an \as/
$\{T_k\}$ satisfying $(\ref{block:umapf:def})$.

$(iii)$
If $X$ has a metric $\ell_p$-\fdd/,
then for each $\eps>0$ there is a 
\fdd/ $\{T_k\}$ such that
$(\ref{block:def:str})$ holds with $C=1+\eps$.
Then, for all $k\ge1$,
$$
(1-\eps)\|T_kx\|\le\biggl(\sum_{j=1}^k\|\Delta T_jx\|^p
\biggr)^{1/p}\le(1+\eps)\|T_kx\|
$$
$$
(1-\eps)\|x-T_kx\|\le
\biggl(\sum_{j=k+1}^\infty\|\Delta T_jx\|^p\biggr)^{1/p}
\le(1+\eps)\|x-T_kx\|.
$$
Thus
$$
(1-\eps)/(1+\eps)\|x\|\le
\Bigl(\|x-T_kx\|^p+\|T_kx\|^p\Bigr)^{1/p}
\le(1+\eps)/(1-\eps)\|x\|.
$$
By a diagonal
argument, this gives an \fdd/
$\{T_k\}$ satisfying $(\ref{block:umapf:def})$.
\eck\vskipb

\questsub
What about the converse in Proposition \ref{block:str:eq}$(i)$ for
$p>1$~?

\subsection{Some consequences of \lpap/}

We start with the simple
\begin{prpsub}\label{appen:l1}
  Let $X$ be a separable Banach space.
  \begin{itemize}
  \item [$(i)$]If $X$ has \lpap/ with constant $C$, then $X$ is
  $C$-isomorphic to a subspace of an $\ell_p$-sum of finite
  dimensional subspaces of $X$.
  \item [$(ii)$]If furthermore $X$ is a subspace of ${\fam0 L}^q$, then $X$ is
  $(C+\eps)$-isomorphic to a subspace of $(\bigoplus \ell_q^n)_p$ for
  any given $\eps>0$.
  \item [$(iii)$]In particular, if a subspace of ${\fam0 L}^p$ has \lpap/ with
  constant $C$, then it is $(C+\eps)$-isomorphic to a subspace of
  $\ell_p$ for any given $\eps>0$. If a subspace of ${\fam0 L}^p$ has \lpmap/,
  then it is almost isometric to subspaces of $\ell_p$.
  \end{itemize}
 \end{prpsub}
\dem
$(i)$
Indeed, $\Phi\colon X\hookrightarrow(\bigoplus\im\Delta T_i)_{p}$, 
$x\mapsto\{\Delta T_ix\}_{i\ge1}$
is an embedding: for all $x\in X$
$$
C^{-1}\|x\|_X\le\|\Phi x\|=
\Bigl(\sum\|\Delta T_ix\|^p_X\Bigr)^{1/p}\le C\|x\|_X.
$$

$(ii\,\&\,iii)$
Recall that, given $\eps>0$, a
finite dimensional subspace of ${\fam0 L}^q$
is $(1+\eps)$-isomorphic to a subspace
of $\ell_q^n$ for some $n\ge1$.
\eck\vskipb

We have in particular
(see \cite[\S VIII, Def.\ 7]{hmp86}
for the definition of Hilbert sets)
\begin{corsub}\label{appen:uap}
  
  Let \DE/ be infinite.
  \begin{itemize}
  \item [$(i)$]No $\SLE{q}$ 
\index{p-additive approximation property@$p$-additive approximation property!for spaces $\SLPE$}
  $(1\le q<\infty)$ has \lpap/ for $p\ne2$.
  \item [$(ii)$]No $\SCE$ 
\index{p-additive approximation property@$p$-additive approximation property!for spaces $\SCE$}
  has \lap{q} for $q\ne1$.  If $E$ is a Hilbert\index{Hilbert set}
  set, then $\SCE$ fails \lap{1}.
  \end{itemize}
\end{corsub}
\dem
This is a consequence of Proposition \ref{appen:l1}$(i)$: every
infinite $E$ contains a Sidon set and
thus a \EL{2\vee p} set. So $\SLPE$ contains $\ell_2$.
Also, if $E$ is a Hilbert set, then $\SCE$ contains $c_0$
by \cite[Th.\ 2]{li95}.
\eck\vskipb

However, there is a Hilbert set $E$ such that
$\SCE$ has complex \umap/: see \cite[Th.\ 10]{li96}.
The class of sets $E$ such that $\SCE$ has \lap{1}
contains the Sidon\index{Sidon set} sets and 
Blei's\index{Blei, Ron C.} 
sup-norm-partitioned sets\index{sup-norm-partitioned sets} 
\cite{bl74}.

\subsection{Characterization of \lpmap/}

Recall \cite[Def.\ 4.1]{kw95}: 
\begin{dfnsub}\label{block:strong:dfn}
  Let $X$ be a separable Banach space.
  \begin{itemize}
  \item [$(i)$]Let $\tau$ be a vector space topology on $X$. $X$ enjoys
  property $(m_p(\tau))$ if for all $x\in X$ and norm bounded
  sequences $\{y_j\}$ such that $y_j\hoch{\to}^{\tau}0$
  $$
  \bigl| \|x+y_j\|-\bigl(\|x\|^p+\|y_j\|^p\bigr)^{1/p} \bigr|\to0.
  $$
  \item [$(ii)$]$X$ enjoys the property $(m_p(T_k))$ for a commuting \as/
  $\{T_k\}$ if for all $\eps>0$ and $n\ge1$ we may choose $m\ge n$
  such that for all $x\in B_X$
  $$
  \bigl|\|T_nx+(\Id -T_m)x\|- \bigl(\|T_nx\|^p+\|(\Id
  -T_m)x\|^p\bigr)^{1/p} \bigr|\le\eps.
  $$
  \end{itemize}
\end{dfnsub}
Then \cite[Th.\ 4.2]{kw95} may be read as follows 
\begin{thmsub}\label{block:strong:thm}
  \index{metric $p$-additive approximation property}
  Let $1\le p<\infty$ and consider the following properties for a
  separable Banach space $X$.
  \begin{itemize}
\item [$(i)$]There are an unconditional commuting \as/ $\{T_k\}$ and a
  vector space topology $\tau$ such that $X$ enjoys $(m_p(\tau))$ and
  $T_kx\hoch{\to}^{\tau}x$ uniformly for $x\in B_X$;
\item [$(ii)$]$X$ enjoys the property $(m_p(T_k))$ for an unconditional
  commuting \as/ $\{T_k\}$.
\item [$(iii)$]$X$ has \lpmap/.
\end{itemize} 
  Then $(i)\Leftrightarrow(ii)$.  If $X$ has finite cotype, then
  $(ii)\imp(iii)$.
\end{thmsub}
As for Theorem \ref{block:thm}, we may remove
the cotype assumption if $X$ has commuting $\ell_1$-\ap/ or
$\ell_p$-\fdd/, $p<\infty$:
\begin{thmsub}\label{block:strong:thm'}
  Let $1\le p<\infty$.  Consider the following properties for a
  separable Banach space $X$.
  \begin{itemize}
  \item [$(i)$]There are an $\ell_p$-\fdd/ \(or just a commuting
  $\ell_1$-\as/ in the case $p=1$\) $\{T_k\}$ and a vector space
  topology $\tau$ such that $X$ enjoys $(m_p(\tau))$ and
  $T_kx\hoch{\to}^{\tau}x$ uniformly for $x\in B_X$;
  \item [$(ii)$]$X$ enjoys $(m_p(T_k))$ for an $\ell_p$-\fdd/ \(or just a
  commuting $\ell_1$-\as/ in the case $p=1$\) $\{T_k\}$.
  \item [$(iii)$]$X$ has \lpmap/.
  \end{itemize}
  Then $(i)\Leftrightarrow(ii)\imp(iii)$.
\end{thmsub}

\subsection{Subspaces of ${\fam0 L}^p$ with \lpmap/}
\label{sect:appen}

Although no translation invariant subspace of
$\SL{p}{}$ has \lap{p} for $p\ne2$, Proposition \ref{appen:l1}
$(iii)$ is not void. By the work of 
Godefroy\index{Godefroy, Gilles},
Kalton\index{Kalton, Nigel J.}, Li\index{Li, Daniel} and
Werner\index{Werner, Dirk}
\cite{kw95,gkl96}, we get examples of subspaces of
${\fam0 L}^p$ with \lmap{p} and even a
characterization of such spaces.

Let us treat the case $p=1$.
Recall first that a space $X$ has the
$1$-strong Schur\index{Schur property, $1$-strong} property
when, given $\delta\in{\mathopen]}0,2]$
and $\eps>0$, any
normalized $\delta$-separated sequence in $X$ contains
a subsequence that is $(2/\delta+\eps)$-equivalent to the
unit vector basis of $\ell_1$ (see \cite{ro79}).
In particular, a gliding hump argument
shows that any subspace of $\ell_1$
shares this property. By
Proposition \ref{appen:l1}$(iii)$, a space $X$
with \lmap{1} also does. Now recall the main theorem
of \cite{gkl96}:\vskipa

{\bf Theorem }
{\it Let $X$ be a subspace of ${\fam0 L}^1$ with the approximation
property. Then the following properties are equivalent:}
\begin{itemize}
\item [$(i)$]{\it The unit ball of $X$ is compact and
locally convex in measure;}
\item [$(ii)$]{\it $X$ has \umap/ and the $1$-strong Schur property;}
\item [$(iii)$]{\it $X$ is $(1+\eps)$-isomorphic to a $w^*$-closed 
subspace $X_\eps$ of $\ell_1$ for any $\eps>0$.}
\end{itemize}
We may then add to these three the fourth equivalent property
\begin{itemize}
\item [$(iv)$]{\it
\index{metric $1$-additive approximation property!for subspaces of ${\fam0 L}^1$} 
$X$ has \lmap{1}.}
\end{itemize}
\dem
We just showed that $(ii)$ holds
when $X$ has \lmap{1}. Now suppose we have $(iii)$
and let $\eps>0$.
Thus there is a quotient $Z$ of $c_0$ such that $Z^*$ has
the approximation property and $Z^*$ is
$(1+\eps)$-isomorphic
to $X$.

Let us show that any such $Z^*$ has \lmap{1}.
$Z$ has beforehand the metric approximation
property, with say
$\{R_n\}$, because $Z^*$ has it as a dual separable space. 
By \cite[\protect{Th.\ 2.2}]{gs88}, 
$\{R_n^*\}$ is a metric \as/ in $Z^*$. 
Let $Q$ be the canonical quotient map from $c_0$ onto $Z$.
Let $\{P_n\}$ be the sequence of projections associated
to the natural basis of $c_0$. Then
$\{P_n^*\}$ is also an \as/ in $\ell_1$. Thus
$$
\|P_n^*Q^*x^*-Q^*R_n^*x^*\|_{\ell_1}\to0\qquad\mbox{for any }
x^*\in Z^*.
$$
By Lebesgue's dominated convergence theorem
(see \cite[Th.\ 1]{ka74}), $QP_n-R_nQ\to0$ weakly in the space
$\script{K}(c_0,Z)$ of
compact operators from $c_0$ to $Z$. By Mazur's theorem,
there are 
convex combinations $\{C_n\}$ of
$\{P_n\}$ and $\{D_n\}$ of $\{R_n\}$ such that
$\|QC_n-D_nQ\|_{\script{L}(c_0,Z)}\to0$. Thus
\begin{equation}\label{append:prox}
\|C_n^*Q^*-Q^*D_n^*\|_{\script{L}(Z^*,\ell_1)}\to0.
\end{equation}
Furthermore $C_n^*:\ell_1\to\ell_1$ has the form
$C_n^*(x_1,x_2,\dots)=(t_1x_1,t_2x_2,\dots)$
with $0\le t_i\le1$. Therefore, defining $Q^*a=(a_1,a_2,\dots)$,
\begin{eqnarray}
\nonumber\lefteqn{\|C_n^*Q^*a\|_1+\|Q^*a-C_n^*Q^*a\|_1=}
\qquad&&\\
\nonumber&=&\|(t_1a_1,t_2a_2,\dots)\|_1+
\|((1-t_1)a_1,(1-t_2)a_2,\dots)\|_1\\
\label{appen:fin}
&=&\sum(|t_i|+|1-t_i|)|a_i|=\sum|a_i|=\|Q^*a\|_1.
\end{eqnarray}
As $\{D_n^*\}$ is
still an \as/ for $Z^*$,
$\{D_n^*\}$ realizes \lmap{1} in $Z^*$ by
$(\ref{appen:fin})$,
$(\ref{append:prox})$ and Proposition \ref{block:str:eq}$(i)$.

Thus $X$ has \lap{1} with constant $1+2\eps$.
As $\eps$ is arbitrary, $X$ has \lmap{1}.
\eck\vskipb

For $1<p<\infty$, we have similarly by \cite[Th.\ 4.2]{kw95}
\begin{prpsub}
\index{metric $p$-additive approximation property!for subspaces of ${\fam0 L}^p$}
  Let $1<p<\infty$ and $X$ be a subspace of ${\fam0 L}^p$ with the
  approximation property. The following are equivalent:
  \begin{itemize}
  \item [$(i)$]$X$ is $(1+\eps)$-isomorphic to a subspace $X_\eps$ of
  $\ell_p$ for any $\eps>0$.
  \item [$(ii)$]$X$ has \lpmap/.
  \end{itemize}
\end{prpsub}
\dem
$(ii)\imp(i)$ is in Proposition \ref{appen:l1}. 
For $(i)\imp(ii)$, it suffices 
to prove that any subspace $Z$ of $\ell_p$ 
with the approximation property has \lpmap/.

As $Z$ is reflexive, $Z$ admits a commuting shrinking \as/
$\{R_n\}$. Let $i$ be the injection of $Z$ into $\ell_p$. Let
$\{P_n\}$ be the sequence of projections associated to the
natural basis of $\ell_p$. It is also an \as/ for $\ell_{p'}$. Thus
$$
\|i^*P_n^*x^*-R_n^*i^*x^*\|_{Z^*}\to0\qquad
\mbox{for any }x^*\in\ell_{p'}.
$$
As before,
there are convex combinations $\{C_n\}$ of
$\{P_n\}$ and $\{D_n\}$ of $\{R_n\}$ such that
$\|C_ni-iD_n\|\to0$. 
The convex combinations are finite and may be chosen not to
overlap, so that for each $n\ge1$ there is
$m>n$ such that 
$$
\|C_nx+(\Id -C_m)x\|=\bigl(\|C_nx\|^p+\|(\Id -C_m)x\|^p\bigr)^{1/p}
$$
for $x\in\ell_p$. 
Thus $Z$ satisfies the property $(m_p(D_n))$. 
Following the lines of \cite[Lemma 1]{fe80}, we observe that
$\{D_n\}$ is a commuting unconditional \as/ since $\{P_n\}$ is.
By Theorem \ref{block:strong:thm}, $Z$
has \lpmap/.\eck

\section{\uap/ and \umap/ 
in translation invariant spaces}\label{sect:tis}

Recall that $\U$ is a subgroup of $\T$. 
If $\U=\D=\{-1,1\}$, the following applies to real \umap/. 
If $\U=\T=\{\epsilon\in\C:|\epsilon|=1\}$, it applies to complex
\umap/.

\subsection{General properties. Isomorphic case}

$\SLP{}$ spaces $(1<p<\infty)$ are known to have an 
unconditional basis; furthermore, they have an 
unconditional \fdd/ in translation invariant subspaces 
$\SLP{I_k}$: this 
is a corollary of Littlewood--Paley
\index{Littlewood--Paley partition} 
theory \cite{lp31}. 
One may choose $I_0=\{0\}$ and 
$I_k=\mathopen]-2^k,-2^{k-1}]\cup
[2^{k-1},2^k\mathclose[$. Thus any 
$\SLPE$ $(1<p<\infty)$ has an unconditional \fdd/ 
in translation invariant subspaces 
$\SLP{E\cap I_k}$. The
spaces $\SL{1}{}$ and $\SC{}$, however, do not even have 
\uap/.
 
%
\begin{prpsub}[see \protect{\cite[Lemma 5, Cor.\ 6, Th.\ 7]
    {li96}}]\label{sbd:ssesp}
  Let $E\se\Z$ and $X$ be a homogeneous Banach space on \T.
  \begin{itemize}
  \item [$(i)$]If $X_E$ has \umap/ \(\vs \uap/, \lap{1} or \lmap{1}\), then
  some \as/ of multipliers realizes it.
  \item [$(ii)$]Let $F\se E$. If $X_E$ has \umap/ \(\vs \uap/, \lap{1} or
  \lmap{1}\), then so does $X_F$.
  \item [$(iii)$]If $\SCE$ has \umap/ \(\vs \uap/\),
then so does $X_E$.
  \end{itemize}
\end{prpsub}
Note the important property that \as/ of multipliers commute and 
commute with one another.

Whereas \uap/ is always satisfied for 
$\SLPE$ 
$(1<p<\infty)$, we have the following generalization of 
\cite[remark after Th.\ 7, Prop.\ 9]{li96} for the spaces 
$\SL{1}{E}$ and $\SCE$. By the method of \cite{gk95},
\begin{lemsub}\label{uap:sep}
  If $X$ has \uap/ with a commuting \as/ and $X\not\supseteq c_0$,
  then $X$ is a dual space.
\end{lemsub}
\dem
Suppose $\{T_n\}$ is a commuting \as/ such that 
$(\ref{block:uapdef})$ holds. As $X\not\supseteq c_0$, 
$Px^{**}=\lim T_n^{**}x^{**}$ is well defined 
for each $x^{**}\in X^{**}$. As 
$\{T_n\}$ is an \as/, $P$ is a projection onto $X$.
Let us show that $\ker P$ is $w^*$-closed. 
Indeed, if $x^{**}\in\ker P$, then 
$$
\|T_n^{**}x^{**}\|=\lim_m\|T^{\vphantom{**}}_mT_n^{**}x^{**}\|=
\lim_m\|T^{\vphantom{**}}_nT_m^{**}x^{**}\|=0
$$
and $T_n^{**}x^{**}=0$. Thus 
$$
\ker P=\bigcap_n\ker T_n^{**}.
$$
Let $M=(\ker P)_\perp$. Then $M^*=X$.\eck
\begin{corsub}
  Let \DE/.
  \begin{itemize}
  \item [$(i)$]If $\SLE{1}$
\index{unconditional approximation property!for spaces $\SLE{1}$} 
has \uap/, then $E$ is a Riesz\index{Riesz set}
  set.
  \item [$(ii)$]If $\SCE$
\index{unconditional approximation property!for spaces $\SCE$} 
has \uap/ and $\SCE\not\supseteq c_0$, then $E$
  is a Rosenthal set\index{Rosenthal set}.
  \end{itemize}
\end{corsub}
\dem
In both cases, Lemma \ref{uap:sep} shows that 
the two spaces are separable dual spaces and thus have the 
Radon--Nikodym property. We may now apply 
Lust-Piquard's\index{Lust-Piquard, Fran\c coise} characterization
\cite{lu76}.\eck\vskipb 

There are Riesz sets $E$ such that $\SLE{1}$ fails \uap/: indeed, the
family of Riesz sets is coanalytic \cite{ta88}
while the second condition
is in fact analytic. There are Rosenthal sets that cannot be
sup-norm-partitioned\index{sup-norm-partitioned sets} \cite{bl74}.

The converse of Proposition \ref{sbd:ssesp}$(iii)$ does not hold: 
$\SLE{1}$ may have \uap/ while $\SCE$ fails this property. We have 
\begin{prpsub}\label{uap:uncond}
  Let \DE/.
  \begin{itemize}
  \item [$(i)$]The Hardy space $H^1(\T)=\SL{1}{\sN}$ has \uap/.
  \item [$(ii)$]The disc algebra $A(\T)=\SC{\sN}$ fails \uap/. More
  generally, if $\Z\setminus E$ is a Riesz set, then $\SCE$ fails
  \uap/.
  \end{itemize}
\end{prpsub}
\dem
$(i)$
Indeed, $H^1(\T)$ has an unconditional basis \cite{ma80}. Note 
that the first unconditional \as/ for $H^1(\T)$ appears in 
\cite[\S II, introduction]{me68} with the help of 
Stein's\index{Stein, Elias} \cite{st66a,st66b} 
multiplier theorem (see also \cite{wo84}).

$(ii)$ 
Let $\Delta\subset\T$ be the Cantor set. By 
Bishop's\index{Bishop, Errett A.} 
improvement \cite{bi62} of 
Rudin\index{Rudin, Walter}--Carleson's\index{Carleson, Lennart} 
interpolation theorem, every function in 
$\script{C}(\Delta)$ extends to a function in $\SCE$ if 
$\Z\setminus E$ is a Riesz set. By 
\cite[main theorem]{pe64}, this implies that $\script{C}(\Delta)$ 
embeds in $\SCE$. Then $\SCE$ cannot have \uap/; otherwise
$\script{C}(\Delta)$ would 
embed in a space with an unconditional basis, which is false.\eck\vskipb 

\remsub
Recent studies of the Daugavet\index{Daugavet property}
Property by Kadets\index{Kadets, Vladimir M.}
and Werner\index{Werner, Dirk} 
generalize Proposition \ref{uap:uncond}$(ii)$. This property of a
Banach space $X$ states that for
every finite rank operator $T$ on $X$ $\|\Id+T\|=1+\|T\|$. 
By \cite[Th.\ 2.1]{ka96}, such an $X$ cannot have \uap/. Further, by
\cite[Th.\ 3.7]{we97}, $\SCE$ has the Daugavet Property if 
$\Z\setminus E$ is a
so-called semi-Riesz set\index{semi-Riesz set}, that is if 
all measures with
Fourier spectrum in 
$\Z\setminus E$ are diffuse. 

\vskipa\questsub
Is there some characterization of sets \DE/ such that $\SCE$ has
\uap/~? Only a few classes of such sets are known: Sidon sets and
sup-norm-partitioned sets, for which $\SCE$ even has \lap{1}; certain
Hilbert sets. Adapting the argument in \cite{ri00}, we get that $\SCE$
fails \uap/ if $E$ contains the sum of two infinite sets.

\subsection{Characterization of \umap/ and \lpmap/}

Let us introduce 
\begin{dfnsub}\label{block:block:def}
  Let \DE/ and $X$ be a homogeneous Banach space on \T.
  
  $E$ enjoys the Fourier block\index{Fourier block unconditionality}
  unconditionality property \UP/ in $X$ whenever, for any $\eps>0$ and
  finite $F\se E$, there is a finite $G\se E$ such that for $f\in
  B_{X_F}$ and $g\in B_{X_{E\setminus G}}$
\begin{equation}\label{block:bloc}
\osc_{\epsilon\in\sU}\|\epsilon f+g\|_X\le\eps.
\end{equation}
\end{dfnsub}
\begin{lemsub}\label{blockapp:lem}
  Let \DE/ and $X$ be a homogeneous Banach space on \T. The following
  are equivalent.
  \begin{itemize}
  \item [$(i)$]$X_E$ has $(u(\tau_f))$, where $\tau_f$ is the topology of
  pointwise convergence of the Fourier coefficients:
  $$
  x_n\hoch{\to}^{\tau_f}0\quad\Longleftrightarrow\quad \forall k\ 
  \widehat{x_n}(k)\to0.
  $$
  \item [$(ii)$]$E$ enjoys \UP/ in $X$.
  \item [$(iii)$]$X_E$ enjoys the property of block unconditionality for
  any, or equivalently for some, \as/ of multipliers $\{T_k\}$.
  \end{itemize}
\end{lemsub}
\dem
$(i)\imp(ii)$. 
Suppose that $(ii)$ fails: there are 
$\eps>0$ and a finite $F$ such that for each finite $G$, there 
are $x_G\in B_{X_F}$
and $y_G\in B_{X_{E\setminus G}}$ 
such that 
$$
\osc_{\epsilon\in\sU}\|\epsilon x_G+y_G\|>\eps.
$$
As $B_{X_F}$ is compact, we may suppose $x_G=x$. 
As $y_G\hoch{\to}^{\tau_f}0$, $(u(\tau_f))$ fails. 

$(ii)\imp(iii)$. 
Let $C$ be a uniform bound for $\|T_k\|$. 
Let $n\ge1$ and $\eps>0$. Let $F$ be the finite spectrum of 
$T_n$. Let $G$ be such that $(\ref{block:bloc})$ holds 
for all $f\in B_{X_F}$ 
and $g\in B_{X_{E\setminus G}}$. 
Now there is a term $V$ in de la Vall\'ee-Poussin's 
\as/ such that $V|_{X_G}=\Id|_{X_G}$ and 
$\|V\|_{\script{L}(X_E)}\le3$. 
As $V$ has finite rank, we may choose $m>n$ such that 
$\|(\Id -T_m)V\|_{\script{L}(X_E)}=\|V(\Id -T_m)\|_{\script{L}(X_E)}\le\eps$. 
Let then $x\in T_nB_{X_E}$ and $y\in(\Id -T_m)B_{X_E}$. 
We have 
\begin{eqnarray*}
\|\epsilon x+y\|&\le&\|\epsilon x+(\Id -V)y\|+\eps
\hoch{\le}^{(\ref{block:bloc})}\|x+(\Id -V)y\|+4(C+1)\eps+\eps\\
&\le&\|x+y\|+(4C+6)\eps.
\end{eqnarray*}

$(iii)\imp(i)$ is proved as Lemma \ref{block:lm} 
$(ii)\imp(i)$: note that if $y_j\hoch{\to}^{\tau_f}0$, 
then $\|Ty_j\|\to0$ for any finite rank multiplier $T$.
\eck
\vskipa 
We may now prove the main result of this section.
\begin{thmsub}\label{sbd:thm}
\index{metric unconditional approximation property!for homogeneous Banach spaces}
  Let \DE/ and $X$ be a homogeneous Banach space on \T.  If $X_E$ has
  \umap/, then $E$ enjoys \UP/ in $X$. Conversely, if $E$ enjoys \UP/
  in $X$ and furthermore $X_E$ has \uap/ and finite cotype, or simply
  \lap{1}, then $X_E$ has \umap/. In particular,
  \begin{itemize}
  \item [$(i)$]For $1<p<\infty$, $\SLPE$ has \umap/ if and only if $E$
  enjoys \UP/ in $\SLP{}$.
  \item [$(ii)$]$\SLE{1}$ has \umap/ if and only if $E$ enjoys \UP/ in
  $\SL{1}{}$ and $\SLE{1}$ has \uap/.
  \item [$(iii)$]If $E$ enjoys \UP/ in $\SC{}$ and $\SCE$ has \lap{1}, in
  particular if $E$ is a Sidon set, then $\SCE$ has \umap/.
  \end{itemize}
\end{thmsub}
\dem
Notice first that \umap/ implies \UP/ by Lemma \ref{blockapp:lem} 
$(iii)\imp(ii)$. 

$(i)$ 
Notice that $\SLPE$ $(1< p<\infty)$ has an unconditional \fdd/ of 
multipliers $\{\pi_{E\cap I_k}\}$ and cotype $2\vee p$. Thus \UP/
implies 
\umap/ by Theorem \ref{block:thm'}$(ii)\imp(iv)$. 

By Lemma \ref{blockapp:lem}, 
part $(ii)$ and $(iii)$ follow from 
Theorem \ref{block:thm}$(ii)\imp(iv)$ and 
Theorem \ref{block:thm'}$(ii)\imp(iv)$ respectively.
\eck\vskipb

\remsub\label{rem624}
Consider the special case $E=\{0\}\cup\{j^k\}_{k\ge0}$, $|j|\ge2$, and
suppose $X_E$ has complex \umap/. By Theorem \ref{sbd:thm},
$$
\osc_{\epsilon\in\sT}\|\epsilon a+b\e_{j^k}+c\e_{j^{k+1}}\|
\tol_{k\to\infty}0.
$$ 
Let us show that then $\{0,1,j\}$ is a
$1$-unconditional basic sequence in $X$. Indeed, for any 
$\epsilon,\mu,\nu\in\T$, and choosing $\kappa$ such that
$\mu\kappa=\nu\kappa^j$, 
\begin{eqnarray*}
%
%
\|\epsilon a+\mu b\e_1+\nu c\e_j\|&=&
\|\epsilon a+\mu\kappa b\e_1+\nu\kappa^jc\e_j\|\\&=&
\|\epsilon\overline{\mu\kappa}a+b\e_1+c\e_j\|=
\|\epsilon\overline{\mu\kappa}a+b\e_{j^k}+c\e_{j^{k+1}}\|
\end{eqnarray*}
whose oscillation tends to $0$ with $k$. By Proposition
\ref{mub:isom}$(i)$, $X_E$ fails complex \umap/ if $X$ 
is $\SLP{}$, $p$ not an even integer, or $\SC{}$. By Proposition
\ref{mub:isom}$(ii)$, $\SLE{2n}$, $n\ge1$ integer, fails complex 
\umap/ if $j$ is positive and $n\ge j$, 
or if $j$ is negative and $n\ge|j|+1$.\vskipb

The study of \lpmap/ in $X_E$ reduces to the trivial case $p=2$,
$X=\SL{2}{}$, and to the case $p=1$, $X=\SC{}$. To see this, note that we have
by a repetition of the arguments of Lemma \ref{blockapp:lem} 
\begin{lemsub}\label{sbd:lpmap:lem}
  Let \DE/ and $X$ be a homogeneous Banach space.  The following
  properties are equivalent.
  \begin{itemize}
  \item [$(i)$]  $X_E$ has $m_p(\tau_f)$.
  \item [$(ii)$]  $E$ enjoys the following property $\script{M}_p$ in $X$: for
  any $\eps>0$ and finite $F\se E$, there is a finite $G\se F$ such
  that for $f\in B_{X_F}$ and $g\in B_{X_{E\setminus G}}$
  $$
  \bigl|\|f+g\|_X-(\|f\|_X^p+\|g\|_X^p)^{1/p}\bigr|\le\eps
  $$
  \item [$(iii)$]  $X_E$ enjoys $m_p(T_k)$ for any, or equivalently for some,
  \as/ of multipliers.
  \end{itemize}
\end{lemsub}
\begin{prpsub}
\index{metric $p$-additive approximation property!for homogeneous Banach spaces}
  Let \DE/ and $X$ be a homogeneous Banach space.
  \begin{itemize}
  \item [$(i)$] If $X_E$ has \lpmap/, then $E$ enjoys $\script{M}_p$ in $X$.
  \item [$(ii)$] $\SLE{q}$ has \lpmap/ if and only if $p=q=2$.
  \item [$(iii)$] \index{metric $1$-additive approximation property!for spaces $\SCE$} 
$\SCE$ has \lmap{1} if and only if it has \lap{1}
  and $E$ enjoys $\script{M}_1$ in $\SC{}$: for all $\eps>0$ and
  finite $F\se E$, there is a finite $G\se E$ such that
  $$
  \forall f\in\SC{F}\ \forall g\in\SC{E\setminus G}\qquad
  \|f\|_\infty+\|g\|_\infty\le(1+\eps)\|f+g\|_\infty.
  $$
  \end{itemize}
\end{prpsub}
\dem 
$(i)$ 
Let $\eps>0$. Let $\{T_k\}$ be
an a.s.\ of multipliers that satisfies \Ref{block:def:str} with
$C<1+\eps$. By the
argument of \cite[Lemma 5]{li96}, we may assume that the $T_k$'s have 
their range in \PTE/.
Let $n\ge1$ be such that $\bigl(\sum_{k>n}\|\Delta
T_kf\|_X^p\bigr)^{1/p}<\eps$ for $f\in B_{X_F}$. Let $G$ be such that
$T_kg=0$ for $k\le n$ and $g\in X_{E\setminus G}$. Then successively 
$$
\Bigl|\Bigl(
\sum_{k\le n}\|\Delta T_k(f+g)\|_X^p\Bigr)^{1/p}-
\Bigl(\sum\|\Delta T_kf\|_X^p\Bigr)^{1/p}\Bigr|
\le\eps,
$$
$$
\Bigl|\Bigl(
\sum_{k> n}\|\Delta T_k(f+g)\|_X^p
\Bigr)^{1/p}-
\Bigl(\sum\|\Delta T_kg\|_X^p\Bigr)^{1/p}\Bigr|
\le\eps,
$$
$$
\Bigl|\Bigl(\sum
\|\Delta T_k(f+g)\|_X^p\Bigr)^{1/p}-
\Bigl(\sum\|\Delta T_kf\|_X^p+\sum\|\Delta T_kg\|_X^p\Bigr)^{1/p}
\Bigr|\le2^{1/p}\eps
$$
and 
$$
\bigl|
\|f+g\|_X-(\|f\|_X^p+\|g\|_X^p)^{1/p}
\bigr|\le2\eps(1+2^{1/p}).
$$
$(ii)$ 
By Corollary \ref{appen:uap}, we necessarily have $p=2$. Furthermore,
if $\SLE{q}$ has \lmap{2}, then by property $\script{M}_2$
$$
\bigl|\|\e_n+\e_m\|_q-\sqrt{2}\bigr|\tol_{m\to\infty}0.
$$
Now
$\|\e_n+\e_m\|_q=\|1+\e_1\|_q$ 
is constant and differs from $\|1+\e_1\|_2=\sqrt{2}$ unless $q=2$:
otherwise the only case of equality of the norms $\|\cdot\|_q$ and
$\|\cdot\|_2$ occurs for almost everywhere constant functions.

$(iii)$
Use Theorem \ref{block:strong:thm'}.
\eck

\section{Property \umap/ and arithmetical block independence}
\label{sect:umap}

We may now apply the technique used in the 
investigation of \umbs/ in order to obtain 
arithmetical conditions analogous to \I{n} (see Def.\ \ref{mub:def:ar})
for \umap/. According to Theorem \ref{sbd:thm}, it suffices 
to investigate property \UP/ of block unconditionality: we have 
to compute an expression of type $\|f+\epsilon g\|_p$, where the spectra 
of $f$ and $g$ are far apart and $\epsilon\in\U$. As before, $\U=\T$ 
(\vs $\U=\D$) is the complex (\vs real) choice of signs.

\subsection{Property of block independence}\label{ss:block}

To this end, 
we return to the notation of Computational lemmas 
\ref{mub:calcul} and \ref{mub:culcul}. Define 
\begin{eqnarray}
\lefteqn{\nonumber\Psi_r(\epsilon,z)\quad=\quad\Theta_r((\overbrace{1,\dots,1}^j,
\overbrace{\epsilon,\dots,\epsilon}^{m-j}), z)}&&\\
\nonumber&=&\!\int
\biggl|\e_{r_0}(t)+
\sum_{i=1}^jz_i\e_{r_i}(t)
+\epsilon\sum_{i=j+1}^mz_i\e_{r_i}(t)\biggr|^pdm(t)\\
&=&\sum_{\setbox0\hbox{$\scriptstyle\alpha\in\sN^{m}$}\wd0=0pt\box0}\,\,
{\bip{p/2}\alpha}^2\prod|z_i|^{2\alpha_i}
\label{umap:thet}+
\sum_{\scriptstyle\alpha\ne\beta\in\sN^m\atop\scriptstyle\alpha\sim\beta}
\bip{p/2}\alpha\bip{p/2}\beta
\epsilon^{\mathop{\lower1pt\hbox{$\Sigma$}}_{i>j}\!\alpha_i-\beta_i}
\prod z_i^{\alpha_i}
\overline{z_i}^{\beta_i}.
\end{eqnarray}
As in Computational lemma \ref{mub:culcul}, 
we make the following observation:
\begin{ldcsub}\label{umap:calcul}
  Let $\lstp\zeta0m\in\Z^*$ and $\gamma,\delta$ be as in
  $(\ref{mub:gammadelta})$. If the arithmetic relation
  $(\ref{mub:ArithRel})$ holds, then the coefficient of the term
  $\epsilon^{\lower1pt\hbox{$\Sigma$}_{i>j}\gamma_i-\delta_i} \prod
  z_i^{\gamma_i}\overline{z_i}^{\delta_i}$ in $(\ref{umap:thet})$ is
  $\bip{p/2}\gamma\bip{p/2}\delta$ and thus independent of $r$. If
  $\sum|\zeta_i|\le p$ or $p$ is not an even integer, this coefficient
  is nonzero. If $\smp\zeta0j$ is nonzero \(\vs odd\), then this term is
  nonconstant in $\epsilon\in\U$.
\end{ldcsub}
Thus the following arithmetical property shows up. It is 
similar to property \I{n}
of almost independence.
\begin{dfnsub}\label{arith:def}
  Let \DE/ and $n\ge1$.
  \begin{itemize}
  \item [$(i)$]$E$ enjoys the complex \(\vs real\) property \J{n} of block
  independence\index{block independent set of integers} if for any
  $\zeta\in\Zeta_n^m$ with $\sm\zeta j$ nonzero \(\vs odd\) and given
  $\lst pj\in E$, there is a finite $G\se E$ such that
  $\zeta_1p_1+\dots+\zeta_mp_m\ne0$ for all $\lstp p{j+1}m\in
  E\setminus G$.
  \item [$(ii)$]$E$ enjoys exactly complex \(\vs real\) \J{n} if furthermore
  it fails complex \(\vs real\) \J{n+1}.
  \item [$(iii)$]$E$ enjoys complex \(\vs real\) \J{\infty} if it enjoys
  complex \(\vs real\) \J{n} for all $n\ge1$.
  \end{itemize}
\end{dfnsub}
The complex (\vs real) property \J{n} means precisely the
following. ``For every finite $F\se E$ there is a finite $G\se E$ such
that for any two representations of any $k\in\Z$ as a sum of $n$
elements of $F\cup(E\setminus G)$ 
$$
\sm pn=k=\sm{p'}n
$$
one necessarily has
$$
\mes{\{j:p_j\in F\}}=\mes{\{j:p'_j\in F\}}\hbox{ in $\Z$ (\vs in $\Z/2\Z$).''}
$$
Thus property \J{n} has, unlike \I{n}, 
a complex and a real\index{real vs.\ complex}\index{complex vs.\ real} version. Real \J{n} is strictly 
weaker than complex \J{n}: see Section \ref{sect:arith}. 
Notice that \J{1} is void and $\J{n+1}\imp\J{n}$ in both 
complex and real cases. Also $\I{n}\not\imp\J{n}$: 
we shall see in the following section that $E=\{0\}\cup\{n^k\}_{k\ge0}$ 
provides a counterexample. 
The property \J{2} of 
real block independence appears implicitly in 
\cite[Lemma 12]{li96}.\vskipb

\remsub 
In spite of the intricate form of this 
arithmetical property, \J{n} is the ``simplest'' candidate, 
in some sense, that reflects the features of \UP/: 

\begin{itemize}
\item [\bloc ]it must hold for a set $E$ if and only if it holds 
for a translate $E+k$ of this set: this explains 
$\sum\zeta_i=0$ in Definition \ref{arith:def}$(i)$; 
\item [\bloc ]as for the property \UP/ of block independence, it must 
connect the break\index{break} of $E$ with its tail\index{tail}; 
\item [\bloc ]Li gives an example of a set $E$ whose pace does not 
tend to infinity while $\SCE$ has \lmap{1}. 
Thus no property \J{n} should forbid parallelogram 
relations of the type $p_2-p_1=p_4-p_3$, where 
$p_1,p_2$ are in the break of $E$ and 
$p_3,p_4$ in its tail. This 
explains the condition that $\sm\zeta j$ be nonzero 
(\vs odd) in 
Definition \ref{arith:def}$(i)$.
\end{itemize}

We now repeat the argument of Theorem \ref{mub:thm}
to obtain an analogous statement which relates 
property \UP/ of Definition \ref{block:block:def} 
with our new arithmetical conditions
\begin{lemsub}\label{umap:lem}
  \index{Fourier block unconditionality} 
  \index{block independent set of integers}
  Let \DEE/ and $1\le p<\infty$.
  \begin{itemize}
  \item [$(i)$]Suppose $p$ is an even integer. Then $E$ enjoys the complex
  \(\vs real\) Fourier block unconditionality property \UP/ in
  $\SLP{}$ if and only if $E$ enjoys complex \(\vs real\) \J{p/2}.
  \item [$(ii)$]If $p$ is not an even integer and $E$ enjoys complex \(\vs
  real\) \UP/ in $\SLP{}$, then $E$ enjoys complex \(\vs real\)
  \J{\infty}.
  \end{itemize}
\end{lemsub}
\dem
Let us first prove the necessity of the arithmetical 
property and assume $E$ fails \J{n}: then 
there are $\lstp\zeta0m\in\Z^*$ with $\sum\zeta_i=0$, 
$\sum|\zeta_i|\le2n$ and $\smp\zeta0j$ nonzero (\vs odd); 
there are $\lstp{r}0j\in E$ and sequences 
$\lstp{r^l}{j+1}m\in E\setminus\{\lstp n1l\}$ such that
$$
\zeta_0r_0+\dots+\zeta_jr_j+\zeta_{j+1}r_{j+1}^l+\dots+\zeta_mr_m^l=0.
$$ 

Assume $E$ enjoys \UP/ in $\SLP{}$. Then
the oscillation of $\Psi_r$ in $(\ref{umap:thet})$ satisfies
\begin{equation}\label{umap:ApplDef}
\osc_{\epsilon\in\sU}
\Psi_{r^l}(\epsilon,z)
\tol_{l\to\infty}0
\end{equation}
for each $z\in D^m$. The argument is now 
exactly the same as in Theorem \ref{mub:thm}: 
we may assume that 
the sequence of functions $\Psi_{r^l}$ converges 
in $\script{C}^\infty(\U\times D^m)$ to a function 
$\Psi$. Then by $(\ref{umap:ApplDef})$,
$\Psi(\epsilon,z)$ is constant in $\epsilon$ 
for each $z\in D^m$, and this is impossible by 
Computational lemma \ref{umap:calcul} 
if $p$ is either not an even integer or
$p\ge 2n$.

Let us now prove the sufficiency of \J{p/2} when 
$p$ is an even integer.
First, let $\Alpha_n^{k,l}=
\{\alpha\in \Alpha_n:\alpha_i=0\mbox{ for }k<i\le l\}$ 
($\Alpha_n$ is defined before Prop.\ \ref{mub:isom}), and
convince yourself that \J{p/2} is equivalent to
\begin{equation}\forall k\ \exists l\ge k\ 
\forall\alpha,\beta\in \Alpha_{p/2}^{k,l}\quad 
\sum\alpha_in_i=\sum\beta_in_i\ 
\label{umap:DefEq}\imp\ \sum_{i\le k}\alpha_i=
\sum_{i\le k}\beta_i\ (\mbox{\vs mod $2$}).
\end{equation}
Let $f=\sum a_i\e_{n_i}\in\PTE/$. Let $k\ge1$ and 
$\epsilon\in\U$. By the multinomial formula, 
\begin{eqnarray*}
\lefteqn{\|\epsilon \pi_kf+(f-\pi_{l}f)\|_p^p=
\int{\biggl|\sum_{\alpha\in \Alpha_{p/2}^{k,l}}
\bip {p/2}\alpha\epsilon^
{\lower1pt\hbox{$\Sigma$}_{p\le k}\alpha_i}
\biggl(\prod a_i^{\alpha_i}
\biggr)\e_{\lower1pt\hbox{$\Sigma$}\alpha_in_i}\biggr|}^2dm}\qquad\qquad\\
&&=\int{\biggl|\sum_{j=0}^n\epsilon^j
\sum_
{\scriptstyle\alpha\in \Alpha_{p/2}^{k,l}\atop
 \scriptstyle\sm\alpha k=j}
\bip{p/2}\alpha\biggl(\prod a_i^{\alpha_i}
\biggr)\e_{\lower1pt\hbox{$\Sigma$}\alpha_in_i}\biggr|}^2dm.
\end{eqnarray*}
$(\ref{umap:DefEq})$ now signifies that we may choose 
$l\ge k$ such that the terms of the above sum over $j$ 
(\vs the terms with $j$ odd 
and those with $j$ even) have 
disjoint spectrum. But then 
$\|\epsilon \pi_kf+(f-\pi_{l}f)\|_p$ is constant for 
$\epsilon\in\U$ and $E$ enjoys \UP/ in $\SLP{}$.
\eck\vskipb

Note that for even $p$, we have as in 
Proposition \ref{mub:isom} a constant $C_p>1$ such that either
$(\ref{block:bloc})$ holds for $\eps=0$ or fails 
for any $\eps\le C_p$. We thus get
\begin{corsub}
  Let \DE/ and $p$ be an even integer. If $E$ enjoys complex \(\vs
  real\) \UP/ in $\SLP{}$, then there is a partition $E=\bigcup E_k$ into
  finite sets such that for any coarser partition $E=\bigcup E'_k$
  $$
  \forall f\in\PTE/ \quad\osc_{\epsilon_k\in\sU}\Bigl\|\sum
  \epsilon_k\pi_{E'_{2k}}f\Bigr\|_p=0
  $$
\end{corsub}
Among other consequences, $E=E_1\cup E_2$ 
where the $\SLP{E_i}$ have a complex (\vs real) 
$1$-unconditional 
\fdd/\index{$1$-unconditional fdd@$1$-unconditional \fdd/}.\vskipb

\questsub
Is this rigidity proper to translation invariant subspaces of $\SLP{}$,
$p$ an even integer, or generic for all its subspaces 
(see 
\cite{djp})~?

\subsection{Main result}

Lemma \ref{umap:lem} and Theorem \ref{sbd:thm} yield the 
main result of this section.
\begin{thmsub}\label{umap:thm}
\index{metric unconditional approximation property!for spaces $\SLPE$, $p$ even}
  Let \DE/ and $1\le p<\infty$.
  \begin{itemize}
  \item [$(i)$]Suppose $p$ is an even integer. Then $\SLPE$ has complex
  \(\vs real\) \umap/ if and only if $E$ enjoys complex \(\vs real\)
  \J{p/2}.
  \item [$(ii)$]If $p$ is not an even integer and $\SLPE$ has complex \(\vs
  real\) \umap/, then $E$ enjoys complex \(\vs real\) \J{\infty}.
  \end{itemize}
\end{thmsub}
\begin{corsub}\label{umap:thm:cor}\index{metric unconditional approximation property!for $\SCE$ and $\SLPE$, $p\notin2\N$}

  Let \DE/.
  \begin{itemize}
  \item [$(i)$]If $\SCE$ has complex \(\vs real\) \umap/, then $E$ enjoys
  complex \(\vs real\) \J{\infty}.
  \item [$(ii)$] If any $\SLPE$, $p$ not an even integer, has complex \(\vs
  real\) \umap/, then all $\SLPE$ with $p$ an even integer have
  complex \(\vs real\) \umap/.
  \end{itemize}
\end{corsub}
Suppose $p$ is an even integer. Then Section \ref{sect:arith} 
gives various examples of sets such that $\SLPE$ 
has complex or real \umap/. Proposition \ref{comb:grow}
gives a general growth condition that ensures \umap/.

For $X=\SLP{}$, $p$ not an even integer, 
and $X=\SC{}$, however, we encounter the same obstacle 
as for \umbs/. Section \ref{sect:arith} only gives 
sets $E$ such that $X_E$ fails \umap/. Thus, we 
have to prove this property by direct means. This yields four 
types of examples of sets $E$ such that the space 
$\SCE$ --- and thus by \cite[Th.\ 7]{li96} all 
$\SLPE$ $(1\le p<\infty)$ as well --- have \umap/. 
\begin{itemize}
\item [\bloc]Sets found by Li\index{Li, Daniel} \cite{li96}: 
Kronecker's
theorem is used 
to construct a set containing arbitrarily long 
arithmetic sequences and a set whose pace does not 
tend to infinity. Meyer's\index{Meyer, Yves} \cite[VIII]{me72} techniques are 
used to construct a Hilbert set\index{Hilbert set}.
\item [\bloc]The sets that satisfy the growth condition of Theorem
\ref{positif:thm};
\item [\bloc]Sequences \DEE/ such that $n_{k+1}/n_k$ is an odd integer: 
see Proposition \ref{res:geo}.\vskipb
\end{itemize}
\questsub 
We know no example of a set $E$ such that 
some $\SLPE$, $p$ not an even integer, has \umap/ 
while $\SCE$ fails it.\vskipb

There is also a good arithmetical description of the case 
where $\{\pi_k\}$ or a subsequence 
thereof realizes \umap/. 
\begin{prpsub}\label{umap:prp:fdd}
\index{$1$-unconditional fdd@$1$-unconditional \fdd/!for spaces $\SLPE$, $p$ even} 
\index{metric unconditional fdd@metric unconditional \fdd/} 
Let \DEEE/. Consider a partition
  $E=\bigcup_{k\ge1}E_k$ into finite sets.
  \begin{itemize}
  \item [$(i)$]Suppose $p$ is an even integer. The series $\sum\pi_{E_k}$
  realizes complex \(\vs real\) \umap/ in $\SLPE$ if and only if there
  is an $l\ge1$ such that
\begin{equation}
\label{umap:fdd}
\left\{\begin{array}{l}
\lst pm\in E\\
\zeta_1p_1+\dots+\zeta_mp_m=0
\end{array}\right.
\quad\imp\quad
\forall k\ge l\suml_{p_j\in E_k}\zeta_j=0
\ (\mbox{\vs is even})
\end{equation}
for all $\zeta\in\Zeta_{p/2}^m$.  Then $\SLPE$ admits the series
$\pi_{\cup_{k<l}E_k}+\sum_{k\ge l}\pi_{E_k}$ as
$1$-un\-con\-di\-tio\-nal \fdd/. In particular, choose $E_k=\{n_k\}$.
The sequence $\{\pi_k\}$ realizes complex and real \umap/ in $\SLPE$
if and only if there is a finite $G$ such that for
$\zeta\in\Zeta_{p/2}^m$
\begin{equation}\label{umap:base}
\left\{\begin{array}{l}
\lst pm\in E\\
\zeta_1p_1+\dots+\zeta_mp_m=0
\end{array}\right.
\quad\imp\quad\lst pm\in G.
\end{equation}
Then $E\setminus G$ is a $1$-\ubs/ and $E$ enjoys \I{p/2}.
  \item [$(ii)$]Suppose $p$ is not an even integer. If $\sum\pi_{E_k}$
realizes complex \(\vs real\) \umap/ in $\SLPE$, then for each
$\zeta\in\Zeta^m$ there is an $l\ge1$ such that $(\ref{umap:fdd})$
holds.  In particular, if $\{\pi_k\}$ realizes either complex or real
\umap/ in $\SLPE$, then for all $\zeta\in\Zeta^m$ there is a finite
$G$ such that $(\ref{umap:base})$ holds. This is equivalent to \I{\infty}. 
   \end{itemize}
\end{prpsub}
\dem 
It is analogous to the proof of Lemma \ref{umap:lem}: 
suppose we have $\zeta\in\Zeta_n^m$ such that 
$(\ref{umap:fdd})$ fails for any $l\ge1$. 
Then there are $\lstp\zeta0m\in\Z^*$ with 
$\sum\zeta_i=0$, $\sum|\zeta_i|\le2n$ and 
$\smp\zeta0j$ nonzero (\vs odd) for some $j$; 
for each $l$, 
there are $\lstp{r^l}0j\in\cup_{k<l}E_k$ and 
$\lstp{r^l}{j+1}m\in\cup_{k\ge l}E_k$
such that $\zeta_0r_0^l+\dots+\zeta_mr_m^l=0$. 

But then $\sum\pi_{E_k}$ cannot 
realize complex (\vs real) \umap/: the function 
$\Psi_r$ in $(\ref{umap:thet})$ would satisfy 
$(\ref{umap:ApplDef})$ and we would 
obtain a contradiction 
as in Theorem \ref{mub:thm}.

Sufficiency in $(i)$ and $(i')$ is proved exactly as in 
Lemma \ref{umap:lem}$(i)$.
\eck\vskipb

In particular, suppose that the cardinal 
$\mes{E_k}$ is uniformly 
bounded by $M$ and $\{\pi_{E_k}\}$ realizes \umap/ in 
$\SLPE$. If $p\ne2$ is an even integer, then $E$ 
is a \EL{p}\index{Lambda(p) set@\EL{p} set} set as union of a 
finite set and $M$ $p/2$-independent sets 
(see Prop.\ \ref{mub:isom} and 
\cite[\protect{Th.\ 4.5(b)}]{ru60}). 
If $p$ is not an even 
integer, then $E$ is a \EL{q} set for all $q$ 
by the same argument. 



\section{Examples for \umap/: block independent sets of characters}
\label{sect:arith}

\subsection{General properties}

The pairing $\XE$ underlines the asymptotic nature of property $\J{n}$. 
It has been defined before Proposition \ref{mub:lim}, whose proof adapts to 
\begin{prpsub}\label{arith:lim}
  Let \DEE/.
  \begin{itemize}
  \item [$(i)$]If $\XE<\infty$ for $\lst\zeta m\in\Z^*$ with $\sum\zeta_i$
  nonzero \(\vs odd\), then $E$ fails complex \(\vs real\)
  \J{|\zeta_1|+\dots+|\zeta_m|}.  Conversely, if $E$ fails complex
  \(\vs real\) \J{n}, then there are $\lst\zeta m\in\Z^*$ with
  $\sum\zeta_i$ nonzero \(\vs odd\) and $\sum|\zeta_i|\le2n-1$ such
  that $\XE<\infty$.
  \item [$(ii)$]Thus $E$ enjoys complex \(\vs real\) \J{\infty} if and only
  if $\XE=\infty$ for all $\lst\zeta m\in\Z^*$ with $\sum\zeta_i$
  nonzero \(\vs odd\).
  \end{itemize}
\end{prpsub}
{\it Proof of the converse in $(i)$.\/} If $E$ fails complex (\vs real) 
\J{n}, then there are $\zeta\in\Zeta_n^m$ with 
$\sm\zeta j$ nonzero (\vs odd), $\lst pj\in E$ and 
sequences $\lstp{p^l}{j+1}m\in\{n_k\}_{k\ge l}$ 
such that 
$\sum_{i>j}\zeta_{i}p^l_i=
-\sum_{i\le j}\zeta_ip_i$. Let 
$\zeta'=(\zeta_{j+1},\dots,\zeta_m)$. Then 
$\sum|\zeta'_i|\le2n-1$ and 
$\langle\zeta',E\rangle<\infty$.\eck\vskipb

An immediate application is, as in Proposition \ref{mub:lim}, 
\begin{prpsub}\label{arith:csq}
  Let \DEE/.
  \begin{itemize}
  \item [$(i)$]Suppose $E$ enjoys \I{2n-1}.  Then $E$ enjoys complex \J{n}
  and actually there is a finite set $G$ such that $(\ref{umap:base})$
  holds for $\zeta\in\Zeta_n^m$.
  \item [$(ii)$]Suppose $E$ enjoys \I{\infty}. Then $E$ enjoys complex
  \J{\infty} and actually for all $\zeta\in\Zeta^m$ there is a finite
  $G$ such that $(\ref{umap:base})$ holds.
  \item [$(iii)$]Complex and real \J{\infty} are stable under bounded
  perturbations of $E$.
  \item [$(iv)$]Suppose there is $h\in\Z$ such that $E\cup\{h\}$ fails
  complex \(\vs real\) \J{n}. Then $E$ fails complex \(\vs real\)
  \J{2n-1}. Thus the complex and real properties \J{\infty} are stable
  under unions with an element: if $E$ enjoys it, then so does
  $E\cup\{h\}$.
  \item [$(v)$]Suppose $jF+s,kF+t\in E$ for an infinite $F$, $j\ne k\in\Z^*$
  and $s,t\in\Z$. Then $E$ fails complex \J{|j|+|k|}, and also real
  \J{|j|+|k|} if $j$ and $k$ have different parity.
  \end{itemize}
\end{prpsub}

We now turn to an arithmetical investigation of various sets $E$.

\subsection{Geometric sequences}\index{geometric sequences}

Let $G=\{j^k\}_{k\ge0}$ with 
$j\in\Z\setminus\{-1,0,1\}$. We resume Remark \ref{rem624}.

\vskipa{\bf (1) }
As $G,jG\se G$, $G$ 
fails complex \J{|j|+1}, and also real \J{|j|+1} if 
$j$ is even. The solutions $(\ref{arith:sol})$ 
to the Diophantine equation 
$(\ref{arith:dioph})$ show at once that $G$ enjoys complex
\J{|j|}, since there is no arithmetical relation
$\zeta\in\Zeta_{|j|}^m$ between the break and the tail
of $G$. If $j$ is odd, then $G$ enjoys in fact real \J{\infty}. 
Indeed, let $\lst\zeta m\in\Z^*$ and $k_1<\dots<k_m$: then 
$\sum\zeta_i j^{k_i}\in j^{k_1}\Z$ and either 
$\big|\sum\zeta_i j^{k_i}\big|\ge j^{k_1}$ or 
$\sum\zeta_i j^{k_i}=0$. 
Thus, if $\XE<\infty$ then $\XE=0$ and 
$\sum\zeta_i$ is even 
since $j$ is odd. Now apply Proposition \ref{arith:lim}$(iii)$. 
The same argument yields that even $G\cup-G\cup\{0\}$ 
enjoys real \J{\infty}. Actually much more is true: 
see Proposition \ref{res:geo}.

\vskipa{\bf (2) }
$G\cup\{0\}$ may behave differently 
than $G$ with respect to \J{n}: 
thus this property is not stable under unions with an 
element. Indeed, the first 
solution in $(\ref{arith:sol})$ may be written as 
$(-j+1)\cdot0+j\cdot j^k+(-1)\cdot j^{k+1}=0$. 
If $j$ is positive, $(-j+1)+j+(-1)\le2j$ and 
$G\cup\{0\}$ fails complex \J{j}. 
A look at $(\ref{arith:sol})$ shows that 
it nevertheless 
enjoys complex \J{j-1}. On the other hand, 
$G\cup\{0\}$ still enjoys complex 
\J{|j|} if $j$ is negative. 
In the real setting, our arguments yield 
the same if $j$ is even, but 
we already saw that $G\cup\{0\}$ 
still enjoys real \J{\infty} if $j$ is odd.

\subsection{Symmetric sets}\index{symmetric sets}

By Proposition \ref{mub:lim}$(iii)$ and \ref{arith:csq}$(vi)$, they do 
enjoy neither \I{2} nor complex \J{2}. They may 
nevertheless enjoy real \J{n}. Introduce property 
$(\script{J}^{\mbox{\tiny sym}}_n)$ for $E$: 
it holds if 
for all $\lst pj\in E$ and $\eta\in{\Z^*}^m$ with 
$\sum_1^m\eta_i$ even, $\sum_1^m|\eta_i|\le2n$ and 
$\sm\eta j$ odd, there is a finite set $G$ such that 
$\eta_1p_1+\dots+\eta_mp_m\ne0$ for 
any $\lstp p{j+1}m\in E\setminus G$. Then we obtain
\begin{prpsub}\label{arith:sym} 
  $E\cup-E$ has real \J{n} if and only if $E$ has
  $(\script{J}^{\mbox{\tiny sym}}_n)$.
\end{prpsub}
\dem
By definition, 
$E\cup-E$ has real \J{n} if and only if for all 
$\lst pj\in E$ 
and $\zeta^1,\zeta^2\in\Z^m$ with $\zeta^1+\zeta^2\in\Zeta_n^m$ 
and odd $\sum_{i\le k}\zeta^1_i-\zeta^2_i$, there is a 
finite set $G$ such that $\sum(\zeta^1_i-\zeta^2_i)p_i\ne0$ 
for any $\lstp p{j+1}m\in E\setminus G$ ---~and 
thus if and only if $E$ enjoys 
$(\script{J}^{\mbox{\tiny sym}}_n)$: 
just consider 
the mappings between arithmetical relations 
$(\zeta^1,\zeta^2)\mapsto\eta=\zeta^1-\zeta^2$
and 
$\eta\mapsto(\zeta^1,\zeta^2)$ such that
$\eta=\zeta^1-\zeta^2$, where
$\zeta^1_i=\eta_i/2$ if $\eta_i$ is even and,
noting that the number of odd $\eta_i$'s must be even, 
$\zeta^1_i=(\eta_i-1)/2$ and $\zeta^1_i=(\eta_i+1)/2$ 
respectively for each half of them.\eck

Consider again a geometric sequence $G=\{j^k\}$ with $j\ge2$. 
If $j$ is odd, we saw before that 
$G\cup-G$ and $G\cup-G\cup\{0\}$ enjoy real \J{\infty}.
If $j$ is even, then $G\cup-G$ fails 
real $\J{j+1}$ since $G$ does. 
$G\cup-G\cup\{0\}$ fails real \J{j/2+1}
by the arithmetical relation 
$1\cdot0+j\cdot j^k+(-1)\cdot j^{k+1}=0$
and Proposition \ref{arith:sym}. 
$G\cup-G$ enjoys real 
\J{j} and $G\cup-G\cup\{0\}$ enjoys real 
\J{j/2} as the 
solutions in $(\ref{arith:sol})$ show by a 
simple checking.

\subsection{Algebraic and transcendental numbers}
\index{transcendental numbers}

The proof of Proposition \ref{mub:trans} adapts to 
\begin{prpsub}\label{arith:trans}
  Let \DEE/.
  \begin{itemize}
  \item [$(i)$]If $n_{k+1}/n_k\to\sigma$ where $\sigma>1$ is transcendental,
  then $E$ enjoys complex \J{\infty}.
  \item [$(ii)$]Let $n_k=[\sigma^k]$ with $\sigma>1$ algebraic. Let
  $P(x)=\zeta_0+\dots+\zeta_dx^d$ be the corresponding polynomial of
  minimal degree. Then $E$ fails complex
  \J{|\zeta_0|+\dots+|\zeta_d|}, and also real
  \J{|\zeta_0|+\dots+|\zeta_d|} if $P(1)$ is odd.
  \end{itemize}
\end{prpsub}

\subsection{Polynomial sequences}\index{polynomial sequences}

Let $E=\{P(k)\}$ for a polynomial $P$ of degree $d$. 
The arithmetical relation $(\ref{mub:poly})$ does 
not adapt to property \J{n}. Notice, though, that 
$\{\Delta^jP\}_{j=1}^d$ is a basis for the space 
of polynomials of degree less than $d$ and that 
$2^dP(k)-P(2k)$ is a polynomial of degree at most $d-1$. 
Writing it in the basis $\{\Delta^jP\}_1^d$ yields 
an arithmetical relation 
$2^d\cdot P(k)-1\cdot P(2k)+
\sum_{j=0}^d\zeta_j\cdot P(k-j)=0$ such that 
$2^d-1+\sum\zeta_j$ is odd. By Proposition \ref{arith:lim}
$(ii)$, 
$E$ fails 
real \J{n} for a certain $n$. This $n$ may be bounded in certain cases:

\bloc The set of squares\index{polynomial sequences!squares} fails real \J{2}:
let $F_n$ be the Fibonacci\index{Fibonacci sequence} sequence defined by $F_0=F_1=1$ 
and $F_{n+2}=F_{n+1}+F_n$. As 
$\{F_{n+1}/F_n\}$ is the sequence of convergents of 
the 
continued fraction associated to an irrational 
(the golden 
ratio), $F_n\to\infty$ and 
$F_nF_{n+2}-F_{n+1}^2=(-1)^n$ (see \cite{eu62}). 
Inspired by
\cite[p.\ 15]{mo69}, we observe that
$$
(F_nF_{n+2}+F_{n+1}^2)^2+
(F_{n+1}^2)^2=
(F_nF_{n+1}+F_{n+1}F_{n+2})^2+1^2
$$

\bloc
The set of cubes\index{polynomial sequences!cubes} fails real $\J{2}$:
starting from Binet's\index{Binet, J. P. M.} 
\cite{bi41} simplified solution 
of Euler's\index{Euler, Leonhard} equation \cite{eu56}, we observe that
$p_n=9n^4$, $q_n=1+9n^3$, $r_n=3n(1+3n^3)$ satisfy
$p_n^3+q_n^3=r_n^3+1^3$ and tend to infinity.

\bloc
The set of biquadrates\index{polynomial sequences!biquadrates} fails real $\J{3}$:
by an equality of Ramanujan\index{Ramanujan, Srinivasa} (see \cite[p.\ 386]{ra57}),
$$
(4n^5-5n)^4+(6n^4-3)^4+(4n^4+1)^4=
(4n^5+n)^4+(2n^4-1)^4+3^4.
$$

As for \I{n}, a positive answer to Euler's 
conjecture\index{Euler's conjecture} would imply that 
the set of $k$th powers has complex \J{2} for $k\ge5$. \vskipb

{\bf Conclusion } 
By Theorem \ref{umap:thm}, 
property \J{n} yields directly 
\umap/ in the space $\SL{2p}{}$, $p\le n$ integer. 
But we do not know whether
\J{\infty} ensures \umap/ in spaces 
$\SLP{}$, $p$ not an even integer, or $\SC{}$. 

Nevertheless, 
the study of property \J{3} permits us 
to determine the density of sets 
such that $X_E$ enjoys \umap/ for some $X\ne\SL{2}{}, \SL{4}{}$: see 
Proposition \ref{comb:thm}. Other 
applications are given in Section \ref{sect:resume}.

\section{Positive results: parity and a sufficient 
growth condition}
\label{sect:positif}

\subsection[$\SC{\{3^k\}}$ has real \umap/ because $3$ is odd]
{${\hbox{\zwe C}}_{\{3^k\}}({\hbox{\zwm T}})$ has real \umap/ 
because $3$ is odd}

In the real case, parity plays an unexpected r\^ole.
\begin{prpsub}\label{res:geo}\index{geometric sequences}
  Let \DEE/ and suppose that $n_{k+1}/n_k$ is an odd integer for all
  sufficiently large $k$.  Then $\SCE$ has real \umap/.
\end{prpsub}
Then $X_E$ also has real \umap/ for every homogeneous Banach space $X$
on $\T$.\vskipb

\dem
Let us verify that 
real \UP/ holds. Let $\eps>0$ and 
$F\se E\cap[-n,n]$. Let $l$, to 
be chosen later, such that 
$n_{k+1}/n_k$ is an odd 
integer for $k\ge l$. Take 
$G\supseteq\{\lst nl\}$ finite. 
Let $f\in B_{\script{C}_F}$ and 
$g\in B_{\script{C}_{E\setminus G}}$. Then 
$g(u\exp{\ii\pi/n_l})=-g(u)$ 
and 
$$
|f(u\exp{\ii\pi/n_l})-f(u)|
\le\pi/|n_l|\cdot\|f'\|_\infty\le
\pi n/|n_l|\le\eps
$$ 
by Bernstein's inequality and for $l$ large enough. 
Thus, for some $u\in\T$,
\begin{eqnarray*}
\|f-g\|_\infty&=&|f(u)+g(u\exp{\ii\pi/n_l})|\\
&\le&|f(u\exp{\ii\pi/n_l})+g(u\exp{\ii\pi/n_l})|+\eps\\
&\le&\|f+g\|_\infty+\eps.
\end{eqnarray*}
As $E$ is a Sidon set, 
we may apply Theorem \ref{sbd:thm}$(iii)$.
\eck\vskipb

Furthermore, if $E$ satisfies the hypothesis
of Proposition \ref{res:geo}, so does
$E\cup-E=\{n_1,-n_1,n_2,-n_2,\dots\}$\index{symmetric sets}.
But $E\cup-E$ fails even complex \J{2}
and no $X_{E\cup-E}\ne\SL{2}{E\cup-E}$ has complex \umap/.
On the other hand, if there is an even integer $h$ such that 
$n_{k+1}/n_k=h$ infinitely often, 
then $E$ fails real \J{|h|+1} by Proposition
\ref{arith:csq}$(vi)$.\vskipb 

\remsub
Note that if $n_{k+1}/n_k$ is furthermore uniformly bounded, then the
a.s.\ that realizes \umap/ cannot be too simple. In particular, it
cannot be a \fdd/ in translation invariant spaces $\SC{E_i}$: let $k$ be
such that $n_k$ and $n_{k+1}$ are in distinct $E_i$; then
$n_{k+1}+(-n_{k+1}/n_k)\cdot n_k=0$ and we may apply Proposition
\ref{umap:prp:fdd}$(ii)$. This justifies the use of Theorem
\ref{sbd:thm}$(iii)$.

\subsection[Growth conditions: the case $\SLP{}$, $p$ an even integer]
{Growth conditions: the case ${\fam0 L}^p(\hbox{\zwm T})$, $p$
 an even integer}

For $X=\SLP{}$ with $p$ an even integer, 
a look at \I{n} and \J{n} gives by 
Theorems \ref{mub:thm} and \ref{umap:thm} 
the following general growth condition: 
\begin{prpsub}\label{comb:grow}
\index{$1$-unconditional approximation property!for spaces $\SLPE$, $p$ even} 
\index{$1$-unconditional basic sequence of characters!in spaces $\SLP{}$, $p$ even}
  Let \DEE/ and $p\ge1$ an integer. If
\begin{equation}\label{comb:grow:a}
\liminf|n_{k+1}/n_k|\ge p+1,
\end{equation} 
then $\{\pi_k\}$ realizes the complex \umap/ in $\SLE{2p}$ and there is
a finite $G\se E$ such that $E\setminus G$ is a $1$-unconditional
basic sequence in $\SL{2p}{}$.
\end{prpsub}
\dem
Suppose we have an arithmetical relation 
\begin{equation}\label{comb:eqn}
\zeta_1n_{k_1}+\dots+\zeta_mn_{k_m}=0\quad\mbox{with}\quad
\zeta\in\Zeta_p^m\mbox{ and }|n_{k_1}|<\dots<|n_{k_m}|.
\end{equation} 
Then 
$|\zeta_mn_{k_m}|\le|\zeta_1n_{k_1}|+\dots
+|\zeta_{m-1}n_{k_{m-1}}|$. 
The left hand 
side is smallest when $|\zeta_m|=1$. As 
$|\zeta_1|+\dots+|\zeta_m|\le 2p$ and necessarily 
$|\zeta_i|\le p$, the right hand side is largest 
when $|\zeta_{m-1}|=p$ and $|\zeta_{m-2}|=p-1$. 
Furthermore, it is largest when 
$k_m=k_{m-1}+1=k_{m-2}+2$. 
Thus, if $(\ref{comb:eqn})$ 
holds, then 
$$|n_{k_m}|\le p|n_{k_{m-1}}|+(p-1)|n_{k_{m-2}}|.$$
By $(\ref{comb:grow:a})$, this is impossible as soon 
as $m$ is chosen sufficiently large, because $p+1>p+(p-1)/(p+1)$.
\eck\vskipb

Note that Proposition \ref{comb:grow} is best 
possible: if $j$ is negative, then $\{j^k\}$ fails 
\I{|j|}. If $j$ is positive, 
then $\{j^k\}\cup\{0\}$ fails complex \J{j}.

\subsection{A general growth condition}

Although we could prove that $E$ enjoys \I{\infty} and \J{\infty} when 
$n_{k+1}/n_k\to\infty$, we need a direct 
argument in order 
to get the corresponding functional properties: we have 
\begin{thmsub}\label{positif:thm}
  Let \DEE/ such that $n_{k+1}/n_k\to\infty$.  Then $\SCE$ has
\lmap{1}\index{metric $1$-additive approximation property!for spaces $\SCE$} 
with $\{\pi_k\}$ and $E$ is a 
Sidon\index{Sidon set!with constant asymptotically $1$} 
    set with constant
  asymptotically $1$.  If the ratios $n_{k+1}/n_k$ are all integers,
  then the converse holds.
\end{thmsub}
Note that by Proposition \ref{mub:pinf}$(ii)$, $E$ is a
metric unconditional basic sequence
\index{metric unconditional basic sequence} 
in every homogeneous Banach space
$X$ on $\T$. Further $X_E$ has complex \umap/
\index{metric unconditional approximation property!for homogeneous Banach spaces} 
since $\SCE$ does.\vskipb

\dem
Suppose $|n_{j+1}/n_j|\ge q$ for $j\ge l$ and some $q>1$ to 
be fixed later. Let 
$f=\sum a_j\e_{n_j}\in\PTE/$ and $k\ge l$. We show by 
induction that for all $p\ge k$
\begin{equation}\label{positif:hyp}
\|\pi_pf\|_\infty\ge
\biggl(1-\frac{\pi^2}2\frac{1-q^{2(k-p)}}{{q^2}-1}\biggr)
\|\pi_kf\|_\infty+\sum_{j=k+1}^{p}
\biggl(1-\frac{\pi^2}2\frac{1-q^{2(j-p)}}{{q^2}-1}\biggr)
|a_j|.
\end{equation}

\bloc 
There is nothing to show for $p=k$.

\bloc 
By 
Bernstein's inequality 
applied to $\pi_kf''$ and separately to each 
$a_j\e_{n_j}''$, $j>k$,
\begin{equation}
  \label{positif:bern}
\|\pi_pf''\|_\infty\le n_k^2\|\pi_kf\|_\infty+\sum_{j=k+1}^pn_j^2|a_j|.
\end{equation}
Furthermore, by Lemmas $1$ and $2$ of \cite[\S VIII.4.2]{me72},
\begin{equation}\label{positif:mey}
\|\pi_{p+1}f\|_\infty\ge
\|\pi_pf\|_\infty+|a_{p+1}|-\pi^2/(2n_{p+1}^2)\|\pi_pf''\|_\infty.
\end{equation}
\Ref{positif:mey} together with \Ref{positif:hyp} and \Ref{positif:bern} yield
\Ref{positif:hyp} with $p$ 
replaced by $p+1$. %
%
%
%
Therefore 
\begin{equation}\label{positif:sidon}
\|f\|_\infty=\lim_{p\to\infty}\|\pi_pf\|_\infty\ge
\biggl(1-\frac{\pi^2}2\frac1{{q^2}-1}\biggr)
\biggl(\|\pi_kf\|_\infty+\sum_{j=k+1}^\infty|a_j|\biggr).
\end{equation}
Thus $\{\pi_j\}_{j\ge k}$ realizes 
\lap{1} with constant 
$1+\pi^2/(2q^2-2-\pi^2)$. As $q$ may be chosen 
arbitrarily large, $E$ has \lmap{1} 
with $\{\pi_j\}$. 
Additionally 
$(\ref{positif:sidon})$ shows by choosing $\pi_kf=0$
that $E$ is a \umbs/ in $\SC{}$.

Finally, the converse holds by Proposition \ref{arith:csq}$(vi)$: 
if $n_{k+1}/n_k$ does not tend to 
infinity while being integer, then there are 
$h\in\Z\setminus\{0,1\}$ and 
an infinite $F$ such that $F,hF\se E$.\eck\vskipb

\remsub 
The technique of Riesz products\index{Riesz product} 
as exposed in \cite[Appendix V,
\S1.II]{ks63} would have sufficed to prove Theorem
\ref{positif:thm}.\vskipb

\remsub 
Suppose still that \DEE/ with 
$n_{k+1}/n_k\to\infty$. A 
variation of the above argument yields 
that the space of {\it real\/} 
functions 
with spectrum in $E\cup-E$ has \lap{1}. \vskipb

\remsub 
Note however that there are sets $E$ that satisfy 
$n_{k+1}/n_k\to1$ and nevertheless enjoy 
\I{\infty} (see end of Section \ref{sect:comb}): 
they might be \umbs/ in $\SC{}$, but this is unknown.

\subsection{Sidon constant of Hadamard sets}

Recall that \DEE/ is a Hadamard\index{Hadamard set} set if there is a 
$q>1$ such that $n_{k+1}/n_k\ge q$ for all $k$. It is a 
classical fact that then $E$ is a Sidon set: 
Riesz products\index{Riesz product} (see \cite[Chapter 2]{lr75}) 
even yield effective 
bounds for its Sidon constant. In particular, if $q\ge3$, then $E$'s
Sidon constant is at most $2$. Our computations 
provide an alternative proof for 
$q>\sqrt{\pi^2/2+1}\approx 2.44$ and give a better bound for 
$q>\sqrt{\pi^2+1}\approx 3.30$.
Putting $k=1$ in $(\ref{positif:sidon})$, we obtain
\begin{corsub}\label{positif:cor}
\index{Sidon set!constant}%
Let \DEE/.
\begin{itemize}
\item [$(i)$]Let $q>\sqrt{\pi^2/2+1}$. If $|n_{k+1}|\ge q|n_k|$ for all $k$, then
the Sidon constant of $E$ is at most $1+\pi^2/(2q^2-2-\pi^2).$
\item [$(ii)$]\cite[Cor.~5.2]{ne01} Let $q\ge2$ be an integer. If
$E\supseteq\{n,n+k,n+qk\}$ for some $n$ and $k$, then the Sidon
constant of $E$ is at least
$\bigl(\cos(\pi/2q)\bigr)^{-1}\ge1+\pi^2/(8q^2)$.
\end{itemize}
\end{corsub}
In particular, we have the following bounds for the Sidon constant $C$ of $G=\{j^k\}$, $j\in\Z\setminus\{-1,0,1\}$:
$$
1+\pi^2/(8(j+1)^2)\le C\le 
1+\pi^2/(2j^2-2-\pi^2).$$

\section{Density conditions}
\label{sect:comb}

We apply combinatorial tools to find out how ``big'' 
a set $E$ may be while enjoying \I{n} or \J{n}, 
and how ``small'' it must be.

The coarsest notion of largeness is that of 
density. Recall that the maximal density\index{maximal density} of \DE/ 
is defined by 
$$
d^*(E)=\lim_{h\to\infty}\,\max_{\vphantom{h}a\in\sZ}
\frac{\mes{E\cap\{a+1,\dots,a+h\}}}h.
$$

Suppose $E$ enjoys \I{n}\index{maximal density!of independent sets} 
with $n\ge2$. Then $E$ is a 
\EL{2n} set by Theorem \ref{mub:thm}$(i)$. By 
\cite[Th.\ 3.5]{ru60} (see also \cite[\S1, Cor.\ 2]{mi75}), $d^*(E)=0$. 
Now suppose $E$ enjoys complex or real 
\J{n} with $n\ge2$. As Li\index{Li, Daniel} \cite[Th.\ 2]{li96} 
shows, there are sets $E$ 
such that $\SCE$ has \lmap{1} while $E$ contains 
arbitrarily 
long arithmetic sequences: we cannot apply Szemeredi's Theorem.

Kazhdan\index{Kazhdan, David A.} (see \cite[Th.\ 3.1]{hi82}) 
proved that if $d^*(E)>1/n$, 
then there is a $t\in\{1,\dots,n-1\}$ such that 
$d^*(E\cap E+t)>0$. One 
might hope that it should in fact suffice to choose 
$t$ in any interval of length $n$. 
However, Hindman\index{Hindman, Neil} \cite[Th.\ 3.2]{hi82} 
exhibits a counterexample: 
given $s\in\Z$ and positive 
$\eps$, there is a set $E$ with 
$d^*(E)>1/2-\eps$ and 
there are arbitrarily large $a$ such 
that $E\cap E-t=\emptyset$ for all 
$t\in\{a+1,\dots,a+s\}$. Thus, we have to 
be satisfied with
\begin{lem}\label{com:lem} 
  Let \DE/ with positive maximal density.  Then there is a $t\ge1$
  such that the following holds: for any $s\in\Z$ we have some $a$,
  $|a|\le t$, such that $d^*(E+a\cap E+s)>0$.
\end{lem}
\dem 
By a result of Erd\H{os}\index{Erd\H{o}s, Paul} 
(see \cite[Th.\ 3.8]{hi82}), there 
is a $t\ge1$ such that $F=E+1\cup\dots\cup E+t$ 
satisfies $d^*(F)>1/2$. But then, by \cite[Th.\ 3.4]{hi82}, 
$d^*(F\cap F+s)>0$ for any $s\in\Z$. 
This means that for 
any $s$ there are $1\le a,b\le t$ such that 
$d^*(E+a\cap E+s+b)>0$. 
\eck\vskipb

We are now able to prove 
\begin{prp}\label{comb:thm}
  \index{maximal density!of block independent sets}
  Let \DE/.
  \begin{itemize}
  \item [$(i)$]If $E$ has positive maximal density, then there is an
  $a\in\Z$ such that $E\cup\{a\}$ fails real \J{2}.  Therefore $E$
  fails real \J{3}.
  \item [$(ii)$]If $d^*(E)>1/2$, then $E$ fails real \J{2}.
  \end{itemize}
\end{prp}
\dem 
$(ii)$ is established in \cite[Prop.\ 14]{li96}. 
$(i)$ is a consequence 
of Lemma \ref{com:lem}: indeed, if 
$E$ has positive maximal density, then this lemma yields 
some 
$a\in\Z$ and an infinite $F\se E$ such that 
for all $s\in F$ there are arbitrarily large 
$k,l\in E$ such that $k+a=l+s$. Thus 
$E\cup\{a\}$ 
fails real \J{2}. Furthermore, $E$ fails 
real \J{3} by Proposition \ref{arith:csq}$(iv)$.
\eck\vskipb

\rem 
We may reformulate the remaining open case of \J{2}. 
Let us introduce the infinite difference\index{infinite difference set}
set of $E$: 
$\Delta E=\{t:\mes{E\cap(E-t)}=\infty\}$ (see 
\cite{st79} and \cite{ru78}). Then $E$ has real \J{2} if and only if, 
for any $a\in E$, $\Delta E$ meets $E-a$ 
finitely many times only. 
Thus our question is: are there sets with 
positive maximal density such that 
$E-a\cap\Delta E$ is finite for all $a\in E$~? \vskipb

Proposition \ref{comb:grow} 
and Theorem \ref{positif:thm} show that there is 
only one general condition of lacunarity on $E$ that 
ensures properties \I{n}, \J{n} or \I{\infty}, \J{\infty}: 
$E$ must grow exponentially\index{exponential growth} or 
superexponentially\index{superexponential growth}.
One may nevertheless 
construct inductively 
``large'' sets that enjoy these properties: 
they must only be 
sufficiently irregular to avoid all arithmetical relations. 
Thus there are sequences with growth slower than 
$k^{2n-1}$\index{polynomial growth} which nevertheless 
enjoy both \I{n} and complex and real 
\J{n}. See \cite[\S II, (3.52)]{hr83} 
for a proof in the case $n=2$: 
it can be easily adapted to 
$n\ge 2$ and shows also the way to construct, for any sequence
$n_k\to\infty$, 
sets that satisfy \I{\infty} and \J{\infty} and grow more 
slowly than $k^{n_k}$\index{superpolynomial growth}.

\section{Unconditionality vs.\ probabilistic
independen\-ce}
\label{sect:proba}

\subsection{Cantor group}

Let us first show how simple the 
problems of \umbs/ and \umap/ become when considered 
for independent uniformly distributed random variables 
and their span in some space.

Let $\D^\infty$ be the Cantor group\index{Cantor group} and $\Gamma$ 
its dual group of Walsh functions. 
Consider the set $R=\{r_i\}\se\Gamma$ 
of Rademacher\index{Rademacher functions} functions, \ie 
of the coordinate functions on $\D^\infty$: 
they form a family of independent random variables that 
take values $-1$ and $1$ with equal probability $\frac12$: 
Thus $\|\sum\epsilon_ia_ir_i\|_X$ does not depend on the choice of 
signs $\epsilon_i=\pm1$ for any 
homogeneous Banach space $X$ on $\D^\infty$ 
and $R$ is a real $1$-\ubs/
\index{$1$-unconditional basic sequence of characters!on the Cantor group}
in $X$. 

Clearly, $R$ is 
also a complex \ubs/ in all such $X$. But its complex unconditionality
constant is $\pi/2$\index{real vs.\ complex}\index{complex vs.\ real} 
\cite{se97} and ${\fam0 L}^p_W(\D^\infty)$ has complex \umap/
if and only if $p=2$
\index{metric unconditional approximation property!on the Cantor group}
or $W=\{w_i\}\se\Gamma$ is finite. 
Indeed, $W$ would have an analogue property \UP/ of 
block unconditionality in ${\fam0 L}^p(\D^\infty)$: for any $\eps>0$ 
there would be $n$ such that 
$$
\maxl_{w\in\sT}\|\epsilon a w_1+b w_n\|_p
\le(1+\eps)\|a w_1+b w_n\|_p.
$$
But this is false: 
for $1\le p<2$, take $a=b=1$, $\epsilon=\hbox{i}$:
$$
\maxl_{\epsilon\in\sT}\|\epsilon w_1+ w_n\|_p\ge
\bigl(\hbox{$\frac12$}
(|\hbox{i}+1|^p+|\hbox{i}-1|^p)\bigr)^{1/p}
=\sqrt{2}>\| w_1+ w_n\|_p=2^{1-1/p};
$$
for $2<p\le\infty$, take $a=1$, $b={\fam0 i}$, $\epsilon=\hbox{i}$:
$$
\maxl_{\epsilon\in\sT}\|\epsilon w_1+{\fam0 i} w_n\|_p\ge
\bigl(\hbox{$\frac12$}(|{\fam0 i}+{\fam0 i}|^p+
|{\fam0 i}-{\fam0 i}|^p)\bigr)^{1/p}
=2^{1-1/p}>\| w_1+{\fam0 i} w_n\|_p=\sqrt{2}.
$$

This is simply due to the fact that the image domain 
of the characters on $\D^\infty$ is 
too small. Take now 
the infinite torus $\T^\infty$ 
and consider the 
set $S=\{s_i\}$ of Steinhaus functions, 
\ie the coordinate 
functions on $\T^\infty$: 
they form again a family of independent 
random variables with values uniformly distributed in $\T$. 
Then $S$ is clearly a complex 
$1$-\ubs/ 
\index{$1$-unconditional basic sequence of characters!on the infinite torus}
in any homogeneous Banach space $X$ on $\T^\infty$.

\subsection{Two notions of approximate probabilistic independence}
\label{ss:proba:two}

As the random variables $\{\e_n\}$ also have their 
values uniformly distributed in $\T$, 
some sort of approximate independence 
should suffice to draw 
the same conclusions as in the case of $S$.

A first possibility is to look at the joint 
distribution of 
$(\e_{p_1},\dots,\e_{p_n})$, 
$\lst pn\in E$, and to 
ask it to be close to the product of the distributions 
of the $\e_{p_i}$. For example, 
Pisier\index{Pisier, Gilles}
\cite[Lemma 2.7]{pi83} 
gives the following characterization:
$E$ is a Sidon\index{Sidon set} set if and only if there are a 
neighbourhood $V$ of $1$ 
in $\T$ and 
$0<\varrho<1$ such that for any finite $F\se E$
\begin{equation}\label{proba:pisier}
m[\e_p\in V:p\in F]\le\varrho^{\smes{F}}.
\end{equation}
Murai \index{Murai, Takafumi}\cite[\S4.2]{mu82} 
calls 
\DE/ pseudo-independent\index{pseudo-independent set} if for all $\lst An\se\T$
\begin{equation}\label{proba:pi}
m[\e_{p_i}\in A_i:1\le i\le n]
\tol_{
\scriptstyle p_i\in E\atop
\scriptstyle p_i\to\infty}
\prod_{i=1}^n
m[\e_{p_i}\in A_i]=\prod_{i=1}^nm[A_i].
\end{equation}
We have
%
\begin{prpsub}\label{proba:murai}
  Let \DE/.  The following are equivalent. 
  \begin{itemize}
  \item [$(i)$]$E$ is pseudo-independent,
  \item [$(ii)$]$E$ enjoys
  \I{\infty},
  \item [$(iii)$]For every $\eps>0$ and $m\ge1$, there is a finite subset $G\se E$ such
that the Sidon constant of any subset of
$E\setminus G$ with $m$ elements is less than $1+\eps$. 
  \end{itemize}
\end{prpsub}
Note that by Corollary \ref{arith:cor}, 
$(\ref{proba:pi})$ does not imply $(\ref{proba:pisier})$.

\dem
$(i)\ssi(ii)$  follows by Proposition \ref{arith:lim}$(iii)$ and
\cite[Lemma 30]{mu82}. $(iii)\imp(ii)$ is true because $(iii)$ is just
what is needed to draw our conclusion in Corollary
\ref{arith:cor}. Let us prove $(i)\imp(iii)$. Let $\eps>0$, $m\ge1$
and $\script{A}$ be a covering of $\T$ with intervals of length
$\eps$. By \Ref{proba:pi}, there is a finite set $G\se E$ such that
for $\lst pm\in E\setminus G$ and $A_i\in\script{A}$ we have 
$m[\e_{p_i}\in A_i:A_i\in\script{A}]>0$. But then 
$$\Bigl\|\sum a_i\e_{p_i}\Bigr\|_\infty\ge\sum|a_i|\cdot(1-\eps).\eqno{\ecks}$$

\remsub
$(ii)\imp(iii)$ may be proved directly by the 
technique of Riesz products\index{Riesz product}: 
see \cite[Appendix V, \S1.II]{ks63}.\vskipb

Another possibility is to define some notion of 
almost independence. Berkes \cite{be87} introduces 
the following notion: let us call a sequence 
of random variables $\{X_n\}$ almost 
i.i.d.\index{almost i.i.d. sequence} 
(independent and 
identically distributed) if, after enlarging 
the probability space, there is an i.i.d.\ sequence $\{Y_n\}$ 
such that $\|X_n-Y_n\|_\infty\to0$. We have the straightforward 
\begin{prpsub}\label{proba:berkes}
  Let \DEE/. If $E$ is almost i.i.d., then $E$ is a Sidon
  set\index{Sidon set!with constant asymptotically $1$} with constant
  asymptotically $1$.
\end{prpsub}
\dem
Let $\{Y_j\}$ be an i.i.d sequence and suppose 
$\|\e_{n_j}-Y_j\|_\infty\le\eps$ for $j\ge k$. 
Then 
$$
\sum_{j\ge k}|a_j|= 
\Bigl\|\sum_{j\ge k}a_jY_j\Bigr\|_\infty\le
\Bigl\|\sum_{j\ge k}a_j\e_{n_j}\Bigr\|_\infty+
\eps\sum_{j\ge k}|a_j|
$$
and the unconditionality constant of 
$\{\lstf nk\}$ is less than $(1-\eps)^{-1}$.
\eck\vskipb

Suppose \DEE/ is such that $n_{k+1}/n_k$ 
is an integer for all $k$. 
In that case, Berkes \cite{be87} proves that $E$ is 
almost i.i.d.\ if and only if 
$n_{k+1}/n_k\to\infty$. 
We thus recover a part of Theorem \ref{positif:thm}.

\questsub 
What about the converse in Proposition \ref{proba:berkes}~?

\section{Summary of results. Remarks and questions}
\label{sect:resume}

For the convenience of the reader, we now reorder 
our results by putting together those which are 
relevant to a given class of Banach spaces. 

Let us first summarize our arithmetical results on the geometric sequence
$G=\{j^k\}_{k\ge0}$ ($j\in\Z\setminus\{-1,0,1\}$). The number given in the
first (\vs second, third) column is the value $n\ge1$ for which the set in the
corresponding row achieves exactly \I{n} (\vs complex \J{n}, real
\J{n}).\vskipa

\renewcommand{\arraystretch}{1.3}
\begin{center}
\begin{tabular}{|r||c|c|c|} \hline
$G=\{j^k\}_{k\ge0}$ with $|j|\ge2$&\I{n} &\C-\J{n}&\R-\J{n}\\\hline
$G$, $j>0$ odd\vphantom{\LARGE I} &$|j|$ &$|j|$ &$\infty$\\
$G$, $j>0$ even &$|j|$ &$|j|$ &$|j|$ \\
$G\cup\{0\}$, $j>0$ odd &$|j|$ &$|j|-1$ &$\infty$\\
$G\cup\{0\}$, $j>0$ even &$|j|$ &$|j|-1$ &$|j|-1$ \\
$G$, $G\cup\{0\}$, $j<0$ odd &$|j|-1$&$|j|$ &$\infty$\\
$G$, $G\cup\{0\}$, $j<0$ even &$|j|-1$&$|j|$ &$|j|$ \\
$G\cup-G$, $G\cup-G\cup\{0\}$, $j$ odd&$1$&$1$ &$\infty$\\
$G\cup-G$, $j$ even &$1$ &$1$ &$|j|$ \\
$G\cup-G\cup\{0\}$, $j$ even &$1$ &$1$ &$|j|/2$ \\\hline
\end{tabular}\nopagebreak\smallskip\\\nopagebreak
\refstepcounter{thm}Table \thethm
\end{center}\index{geometric sequences}

\subsection[The case $X=\SLP{}$ with $p$ an even integer]
{The case $X={\fam0 L}^p(\hbox{\zwm T})$ with $p$ an even integer}

Let $p\ge4$ be an even integer. 
We observed the following facts.
\begin{itemize}
\item [\bloc ]Real and complex \umap/ differ among subspaces $\SLPE$ for each $p$:
consider Proposition \ref{res:geo} or $\SLPE$ with $E=\{\pm(p/2)^k\}$.
\item [\bloc ]By Theorem \ref{umap:thm}, $\SLPE$ has complex (\vs real) 
\umap/ if so does $\SLE{p+2}$; 
\item [\bloc ]The converse is false for any $p$. In the complex case, 
$E=\{(p/2)^k\}$ is a counterexample. In the real 
case, take $E=\{0\}\cup\{\pm p^k\}$.
\item [\bloc ]Property \umap/ is 
not stable under unions with an element: for each $p$, 
there is a set $E$ 
such that $\SLPE$ has complex (\vs real) 
\umap/, but $\SLP{E\cup\{0\}}$ does not. 
In the complex case, consider $E=\{(p/2)^k\}$. 
In the real case, consider 
$E=\{\pm(2\lceil p/4\rceil)^k\}$. 
\item [\bloc ]If $E$ is a symmetric set and $p\ne2$, then $\SLPE$ fails complex \umap/. 
Proposition \ref{arith:sym} gives a criterion for real \umap/.
\end{itemize}

What is the relationship between \umbs/ and complex \umap/~? We have by 
Proposition \ref{arith:csq}$(i)$ and \ref{umap:prp:fdd}$(i)$
\begin{prpsub}\label{resume:mb_um}
  Let \DEE/ and $n\ge1$.
  \begin{itemize}
  \item [$(i)$]If $E$ is a \umbs/ in $\SL{4n-2}{}$, then $\SLE{2n}$ has complex
  \umap/.
  \item [$(ii)$]If $\{\pi_k\}$ realizes complex \umap/ in $\SLE{2n}$, then $E$
  is a \umbs/ in $\SL{2n}{}$.
  \end{itemize}
\end{prpsub}

We also have, by Proposition \ref{comb:thm}$(i)$
\begin{prpsub} 
  Let \DE/ and $p\ne2,4$ an even integer. If $\SLPE$ has real \umap/,
  then $d^*(E)=0$.
\end{prpsub}

Note also this consequence of Propositions 
\ref{mub:trans}, \ref{arith:trans}, \ref{proba:murai} and Theorems 
\ref{mub:thm}, \ref{umap:thm} 
\begin{prpsub} 
  Let $\sigma>1$ and $E=\{[\sigma^k]\}$. Then the following properties
  are equivalent:
  \begin{itemize}
  \item [$(i)$]$\sigma$ is transcendental\index{transcendental numbers};
  \item [$(ii)$]$\SLPE$ has complex \umap/ for any even integer $p$;
  \item [$(iii)$]$E$ is a \umbs/ in any $\SLP{}$, $p$ an even integer;
  \item [$(iv)$]$E$ is pseudo-independent.
  \item [$(v)$]For every $\eps>0$ and $m\ge1$, there is an $l$ such that for
  $\lst km\ge l$ the Sidon constant of
  $\{[\sigma^{k_1}],\dots,[\sigma^{k_m}]\}$ is less than $1+\eps$.
   \end{itemize}
\end{prpsub}

\subsection[Cases $X=\SLP{}$ with $p$ not an even integer and $X=\SC{}$]
{Cases $X={\fam0 L}^p(\hbox{\zwm T})$ 
with $p$ not an even integer and $X=\hbox{\zwe C}(\hbox{\zwm T})$}

In this section, $X$ denotes either 
$\SLP{}$, $p$ not an even integer, or $\SC{}$.

Theorems \ref{mub:thm} and \ref{umap:thm} 
only permit us to use the negative results of Section 
\ref{sect:arith}: thus, we can just gather 
negative results about the functional properties 
of $E$. 
For example, we know by 
Proposition \ref{arith:csq}$(iv)$ that \I{\infty} and \J{\infty} are 
stable under union with an element. Nevertheless, 
we cannot conclude that the same holds 
for \umap/. 
The negative results are (by Section \ref{sect:arith}):

\vskipa\bloc
for any infinite \DE/, 
$X_{E\cup2E}$ fails real \umap/. Thus \umap/ 
is not stable under unions; 

\vskipa\bloc 
if $E$ is a polynomial sequence (see Section \ref{sect:arith}), 
then $E$ is not a 
\umbs/ in $X$ and $X_E$ fails real 
\umap/;

\vskipa\bloc 
if $E$ is a symmetric set, then $E$ is not a 
\umbs/ in $X$ 
and $X_E$ fails complex \umap/. Proposition \ref{arith:sym} 
gives a criterion for real \umap/;

\vskipa\bloc 
if $E=\{[\sigma^k]\}$ with 
$\sigma>1$ an algebraic number ---~in particular 
if $E$ is a 
geometric sequence~---, then $E$ is not a 
\umbs/ in $X$ 
and $X_E$ fails complex \umap/.\vskipb 

Furthermore, by Proposition \ref{res:geo}, 
real and complex \umap/ differ in X. 

Theorem \ref{positif:thm} is the only but general 
positive result on \umbs/ and complex \umap/ in $X$.
Proposition \ref{res:geo} yields further examples 
for real \umap/.

What about the sets that satisfy 
\I{\infty} or \J{\infty}~? We only know that \I{\infty} 
does not even ensure Sidonicity by Corollary \ref{arith:cor}.

One might wonder whether for some reasonable class of sets 
$E$, $E$ is a finite 
union of sets that enjoy \I{\infty} or \J{\infty}. 
This is false 
even for Sidon sets: for example, let $E$ be the geometric
sequence $\{j^k\}_{k\ge0}$ with $j\in\Z\setminus\{-1,0,1\}$ 
and suppose $E=E_1\cup\dots\cup E_n$. 
Then $E_i=\{j^k\}_{k\in A_i}$, where the $A_i$'s are a 
partition of the set of positive integers. But then one 
of the $A_i$ contains arbitrarily large $a$ and $b$ 
such that $|a-b|\le n$. This means that there is an infinite 
subset $B\se A_i$ 
and an $h$, $1\le h\le n$, such that $h+B\se A_i$. 
We may apply Proposition \ref{arith:csq}$(vi)$: 
$E_i$ enjoys neither \I{j^h+1} nor complex 
\J{j^h+1} --- nor real \J{j^h+1} 
if furthermore $j$ is even. 

Does Proposition \ref{resume:mb_um}$(ii)$ remain true for general
$X$~? We do not know this. Suppose however that we know that 
$\{\pi_k\}$ realizes \umap/ in the following strong manner: for any
$\eps>0$, a tail $\{\pi_k\}_{k\ge l}$ is a $(1+\eps)$-unconditional
a.s.\ in $X_E$. 
Then $E$ is trivially a \umbs/ in $X$. In particular, this is the case
if 
$$1+\eps_n=
\sup_{\epsilon\in\sU}\|\Id -(1+\epsilon)\pi_n\|_{\script{L}(X)}$$ 
converges so rapidly to $1$ that $\sum\eps_n<\infty$. Indeed,
$$
\sup_{\epsilon_k\in\sU}
\|\pi_{n-1}+\sum_{k\ge n}\epsilon_k\Delta\pi_k\|
\le(1+\eps_n)
\sup_{\epsilon_k\in\sU}
\|\pi_n+\sum_{k>n}\epsilon_k\Delta\pi_k\|.
$$
and thus, for all $f\in\PTE/$,
$$
\sup_{\epsilon_k\in\sU}
\|\pi_lf+\sum_{k>l}\epsilon_k\Delta\pi_kf\|
\le\prod_{k>l}(1+\eps_k)\,\|f\|.
$$


Let us finally state 
\begin{prpsub}
  Let \DE/. If $X_E$ has real \umap/, then $d^*(E)=0$.
\end{prpsub}

\subsection{Questions}

The following questions remain open:

\vskipa{\bf Combinatorics } Regarding Proposition \ref{comb:thm}$(i)$,
is there a set $E$ enjoying \J{2} with positive maximal density, or
even with a uniformly bounded pace~? Furthermore, may a set $E$ with
positive maximal density admit a partition $E=\bigcup E_i$ in finite sets
such that all $E_i+E_j$, $i\le j$, are pairwise disjoint~? Then
$\SLE{4}$ would admit a $1$-unconditional 
\fdd/\index{$1$-unconditional fdd@$1$-unconditional \fdd/!for $\SLE{4}$} 
by Proposition
\ref{umap:prp:fdd}$(i)$.

\vskipa{\bf Functional analysis } Let $X\in\{\SL{1}{},\SC{}\}$ and
consider Theorem \ref{sbd:thm}. Is \UP/ sufficient for $X_E$ to share
\umap/~? Is there a set \DE/ such that some space $\SLPE$, $p$ not an
even integer, has \umap/, while $\SCE$ fails it~?

\vskipa{\bf Harmonic analysis } Is there a Sidon set \DEE/ of constant
asymptotically $1$ such that $n_{k+1}/n_k$ is uniformly bounded~? 
What about the case $E=[\sigma^k]$ for a
transcendental\index{transcendental numbers} 
$\sigma>1$~? If $E$ enjoys \I{\infty}, is $E$ a \umbs/
in $\SLP{}$ $(1\le p<\infty)$~? What about \J{\infty}~?

\vfill\pagebreak
{\addcontentsline{toc}{section}{Bibliography}

}

\vfill\pagebreak

\section*{Index of notation}
\markboth{index of notation}{index of notation}
\addcontentsline{toc}{section}{Index of notation}

\halign{\hfil#&\quad#\hfil\cr
$\mes{B}$&cardinal of $B$\cr
$X_E$&space of $X$-functions with spectrum in $E$\cr
$\widehat{f}$&Fourier transform of $f$: $\widehat{f}(n)=\int
 f(t)\e_{-n}(t)dm(t)$\cr 
$x\choose\alpha$&multinomial number, \S\ref{ss:almost}\cr
$\XE$&pairing of the arithmetical relation $\zeta$ against
 the spectrum $E$, \S\ref{ss:mubs:gen}\cr
$u_n\preccurlyeq v_n$&$|u_n|$ is bounded by $C|v_n|$ for some $C$\cr
\noalign{\vskip4pt minus 0.5pt}
$1$-\ubs/&$1$-unconditional basic sequence of characters, Def.\
\ref{mub:def}$(i)$\cr \noalign{\vskip4pt minus 0.5pt}
$A(\T)$&disc algebra $\SC{\sN}$\cr
$\Alpha_n,\Alpha_n^m$&sets of multi-indices viewed as arithmetic
 relations, \S\ref{ss:isom}\cr
a.s.&approximating sequence, Def.\ \ref{block:def}\cr\noalign{\vskip4pt minus 0.5pt}
$B_X$&unit ball of the Banach space $X$\cr\noalign{\vskip4pt minus 0.5pt}
$\SC{}$&space of continuous functions on $\T$\cr
\noalign{\vskip4pt minus 0.5pt}
$\D$&set of real signs $\{-1,1\}$\cr
$\Delta T_k$&difference sequence of the $T_k$: $\Delta
 T_k=T_k-T_{k-1}$ ($T_0=0$)\cr\noalign{\vskip4pt minus 0.5pt}
$\e_n$&character of $\T$: $\e_n(z)=z^n$ for $z\in\T$, $n\in\Z$\cr
\noalign{\vskip4pt minus 0.5pt}
\fdd/&finite dimensional decomposition, Def.\ \ref{block:def}\cr\noalign{\vskip4pt minus 0.5pt}
$H^1(\T)$&Hardy space $\SL{1}{\sN}$\cr\noalign{\vskip4pt minus 0.5pt}
\I{n}&arithmetical property of almost independence, Def.\
 \ref{mub:def:ar}\cr
$\Id$&identity\cr
i.i.d.&independent identically distributed, \S\ref{ss:proba:two}\cr\noalign{\vskip4pt minus 0.5pt}
\J{n}&arithmetical property of block independence, Def.\
 \ref{arith:def}\cr \noalign{\vskip4pt minus 0.5pt}
$\script{L}(X)$&space of bounded linear operators on the Banach space $X$\cr
$\SLP{}$&Lebesgue space of $p$-integrable functions on $\T$\cr
\lpap/&$p$-additive approximation property, Def.\ \ref{str:def}\cr
\lpmap/&metric $p$-additive approximation property, Def.\
 \ref{str:def}\cr 
$\Lambda(p)$&Rudin's class of lacunary sets, Def.\ \ref{mub:sido}\cr\noalign{\vskip4pt minus 0.5pt}
$\script{M}_p$&functional property of Fourier block $p$-additivity,
Lemma \ref{sbd:lpmap:lem}$(ii)$\cr
$\script{M}(\T)$&space of Radon measures on $\T$\cr
$m[A]$&measure of $A\se\T$\cr
$(m_p(\tau))$&functional property of $\tau$-$p$-additivity,
Def.\ \ref{block:strong:dfn}$(i)$\cr
$(m_p(T_k))$&functional property of commuting block $p$-additivity,
Def.\ \ref{block:strong:dfn}$(ii)$\cr\noalign{\vskip4pt minus 0.5pt}
$\osc f$&oscillation of $f$\cr\noalign{\vskip4pt minus 0.5pt}
$\PT{}$&space of trigonometric polynomials on $\T$\cr
$\pi_j$&projection of $X_E$, $E=\{n_k\}$, onto $X_{\{\lst n j\}}$\cr
$\pi_F$&projection of $X_E$ onto $X_F$\cr\noalign{\vskip4pt minus 0.5pt}
$\U$&real ($\U=\D$) or complex ($\U=\T$) choice of signs\cr
\noalign{\vskip4pt minus 0.5pt}
$(\T,dm)$&unit circle in $\C$ with its normalized Haar measure\cr
$\tau_f$&topology of pointwise convergence of the Fourier 
coefficients, Lemma\ \ref{blockapp:lem}$(i)$\cr\noalign{\vskip4pt minus 0.5pt}
\UP/&functional property of Fourier block unconditionality, Def.\
 \ref{block:block:def}\cr
$(u(\tau))$&functional property of $\tau$-unconditionality, Def.\
\ref{block:def:u}$(i)$\cr
$(u(T_k))$&functional property of commuting block 
unconditionality, Def.\ \ref{block:def:u}$(ii)$\cr
\uap/&unconditional approximation property, Def.\
 \ref{block:def}\cr 
\ubs/&unconditional basic sequence, Def.\ \ref{mub:def}\cr 
\umap/&metric unconditional approximation property, Def.\
 \ref{block:def}\cr
\umbs/&metric unconditional basic sequence, Def.\
 \ref{mub:def}\cr \noalign{\vskip4pt minus 0.5pt}
$\Zeta^m,\Zeta_n^m$&sets of multi-indices viewed as arithmetic relations,
 \S\ref{ss:isom}\cr
}


\addcontentsline{toc}{section}{Index}
\begin{theindex}

  \item $1$-unconditional approximation property
    \subitem for spaces $\SLPE$, $p$ even, 44
  \item $1$-unconditional basic sequence of characters, 16
    \subitem in $\SC{}$ and $\SLP{}$, $p\notin2\N$, 18, 19
    \subitem in spaces $\SLP{}$, $p$ even, 18, 44
    \subitem on the Cantor group, 47
    \subitem on the infinite torus, 48
  \item $1$-unconditional \fdd/, 40
    \subitem for $\SLE{4}$, 51
    \subitem for spaces $\SLPE$, $p$ even, 40

  \indexspace

  \item almost i.i.d. sequence, 48
  \item almost independence, 21
  \item approximating sequence, 25
  \item approximation property, 25
  \item arithmetical relation, 18, 20, 21

  \indexspace

  \item Binet, J. P. M., 43
  \item birelation, 18
  \item Bishop, Errett A., 35
  \item Blei, Ron C., 31
  \item block independent set of integers, 38, 39
  \item boundedly complete approximating sequence, 29
  \item Bourgain, Jean, 14
  \item break, 26, 39

  \indexspace

  \item Cantor group, 16, 19, 47
  \item Carleson, Lennart, 35
  \item Casazza, Peter G., 26
  \item commuting block unconditionality, 26
  \item complex vs.\ real, 16, 19, 38, 47
  \item cotype, 29

  \indexspace

  \item Daugavet property, 35

  \indexspace

  \item equimeasurability, 14
  \item Erd\H{o}s, Paul, 47
  \item Euler's conjecture, 24, 43
  \item Euler, Leonhard, 43
  \item exponential growth, 47

  \indexspace

  \item Fibonacci sequence, 43
  \item finite-dimensional decomposition, 25
  \item Forelli, Frank, 13
  \item Fourier block unconditionality, 35, 39
  \item Fr\'enicle de Bessy, Bernard, 24

  \indexspace

  \item geometric sequences, 23, 42, 44, 49
  \item Godefroy, Gilles, 27, 32

  \indexspace

  \item Hadamard set, 46
  \item Hilbert set, 31, 40
  \item Hindman, Neil, 46
  \item homogeneous Banach space, 15

  \indexspace

  \item independent set of integers, 18
  \item infinite difference set, 47
  \item isometries on ${\fam0 L}^p$, 20

  \indexspace

  \item Kadets, Vladimir M., 35
  \item Kalton, Nigel J., 13, 26, 27, 30, 32
  \item Kazhdan, David A., 46

  \indexspace

  \item \EL{p} set, 17, 41
    \subitem constant, 17
  \item Li, Daniel, 32, 40, 46
  \item Littlewood--Paley partition, 34
  \item Lust-Piquard, Fran\c coise, 35

  \indexspace

  \item maximal density, 46
    \subitem of block independent sets, 47
    \subitem of independent sets, 46
  \item metric $1$-additive approximation property
    \subitem for spaces $\SCE$, 37, 45
    \subitem for subspaces of ${\fam0 L}^1$, 33
  \item metric $p$-additive approximation property, 29, 32
    \subitem for homogeneous Banach spaces, 37
    \subitem for subspaces of ${\fam0 L}^p$, 33
  \item metric unconditional approximation property, 25, 27
    \subitem for $\SCE$ and $\SLPE$, $p\notin2\N$, 40
    \subitem for homogeneous Banach spaces, 36, 45
    \subitem for spaces $\SLPE$, $p$ even, 40
    \subitem on the Cantor group, 47
  \item metric unconditional basic sequence, 16, 22, 45
  \item metric unconditional \fdd/, 26, 40
  \item Meyer, Yves, 40
  \item Murai, Takafumi, 24, 48

  \indexspace

  \item oscillation, 15

  \indexspace

  \item $p$-additive approximation property, 29
    \subitem for spaces $\SCE$, 31
    \subitem for spaces $\SLPE$, 31
  \item Pe\l czy\'nski, Aleksander, 25
  \item Pisier, Gilles, 48
  \item Plotkin, A. I., 13
  \item polynomial growth, 47
  \item polynomial sequences, 24, 43
    \subitem biquadrates, 24, 43
    \subitem cubes, 24, 43
    \subitem squares, 24, 43
  \item pseudo-independent set, 48

  \indexspace

  \item Rademacher functions, 47
  \item Ramanujan, Srinivasa, 43
  \item real vs.\ complex, 16, 19, 38, 47
  \item relative multipliers, 17
    \subitem interpolation, 17
  \item renormings, 14
  \item Riesz product, 45, 46, 48
  \item Riesz set, 34
  \item Rosenblatt, Murray, 13
  \item Rosenthal set, 34
  \item Rosenthal, Haskell Paul, 12, 14
  \item Rudin, Walter, 18, 35

  \indexspace

  \item Schur property, $1$-strong, 32
  \item semi-Riesz set, 35
  \item shrinking approximating sequence, 29
  \item Sidon set, 17, 22, 31, 48
    \subitem constant, 17, 46
    \subitem with constant asymptotically $1$, 22, 45, 49
  \item smoothness, 14
  \item Stein, Elias, 35
  \item strong mixing, 13
  \item sup-norm-partitioned sets, 31, 35
  \item superexponential growth, 47
  \item superpolynomial growth, 47
  \item symmetric sets, 42, 44

  \indexspace

  \item tail, 26, 39
  \item $\tau$-unconditionality, 26
  \item transcendental numbers, 24, 43, 50, 51

  \indexspace

  \item unconditional approximation property, 25
    \subitem for spaces $\SCE$, 34
    \subitem for spaces $\SLE{1}$, 34
  \item unconditional basic sequence of characters, 16
  \item unconditional \fdd/, 25
  \item unconditional skipped blocking decompositions, 27
  \item unconditionality constant, 16
    \subitem in $\SL{4}{}$, 19

  \indexspace

  \item Werner, Dirk, 30, 32, 35
  \item Wojtaszczyk, Przemys\l aw, 25

\end{theindex}
\vfill


\begin{thebibliography}{10}

\bibitem{be87}
I.~Berkes, {\em On almost i.\,i.\,d. subsequences of the trigonometric system},
  in: Functional analysis (Austin, 1986--87), E.~W. Odell and H.~P. Rosenthal
  (eds.), Lect. Notes Math. 1332, Springer, 1988,  54--63.

\bibitem{be90}
\leavevmode\vrule height 2pt depth -1.6pt width 23pt, {\em Probability theory
  of the trigonometric system}, in: Limit theorems in probability and
  statistics (P\'ecs, 1989), I.~Berkes, E.~Cs\'aki and P.~R\'ev\'esz (eds.),
  Coll. Math. Soc. J\'anos Bolyai 57, North-Holland, 1990,  35--58.

\bibitem{bi41}
J.~P.~M. Binet, {\em Note sur une question relative \`a la th\'eorie des
  nombres}, C. R. Acad. Sci. Paris 12 (1841),  248--250.

\bibitem{bi62}
E.~Bishop, {\em A general {R}udin--{C}arleson theorem}, Proc. A. M. S. 13
  (1962),  140--143.

\bibitem{bl74}
R.~C. Blei, {\em A simple {D}iophantine condition in harmonic analysis}, Studia
  Math. 52 (1974/75),  195--202.

\bibitem{br80}
J.~Bourgain and H.~P. Rosenthal, {\em Geometrical implications of certain
  finite dimensional decompositions}, Bull. Soc. Math. Belg. (B) 32 (1980),
  57--82.

\bibitem{br96}
V.~Brouncker, {\em Letter to {J}ohn {W}allis}, in: {\OE}uvres de Fermat 3,
  Gauthier-Villars, 1896,  419--420.

\bibitem{ck91}
P.~G. Casazza and N.~J. Kalton, {\em Notes on approximation properties in
  separable {B}anach spaces}, in: Geometry of {B}anach spaces (Strobl, 1989),
  P.~F.~X. M{\"u}ller and W.~Schachermayer (eds.), London Math. Soc. Lect.
  Notes 158, Cambridge Univ. Press, 1991,  49--63.

\bibitem{djp}
F.~Delbaen, H.~Jarchow and A.~Pe{\l}czy\'nski, {\em Subspaces of ${L}_p$
  isometric to subspaces of $\ell_p$}.
\newblock To appear.

\bibitem{dgz93}
R.~Deville, G.~Godefroy and V.~Zizler, {\em Smoothness and renormings in
  {B}anach spaces}, Pitman Monographs and Surveys 64, Longman, 1993.

\bibitem{di59}
{Diophantus of Alexandria}, {\em Les six livres arithm\'etiques et le livre des
  nombres polygones}, Blanchard, 1959.

\bibitem{dr75}
S.~W. Drury, {\em Birelations and {S}idon sets}, Proc. Amer. Math. Soc. 53
  (1975),  123--128.

\bibitem{ek96}
R.~L. Ekl, {\em Equal sums of four seventh powers}, Math. Comp. 65 (1996),
  1755--1756.

\bibitem{ek98}
\leavevmode\vrule height 2pt depth -1.6pt width 23pt, {\em New results in equal
  sums of like powers}, Math. Comp. 67 (1998),  1309--1315.

\bibitem{eu56}
L.~Euler, {\em Solutio generalis quorundam problematum {D}iophanteorum, quae
  vulgo nonnisi solutiones speciales admittere videntur}, in: Op. Omnia (I) II,
  Teubner, 1915,  428--458.

\bibitem{eu72}
\leavevmode\vrule height 2pt depth -1.6pt width 23pt, {\em Observationes circa
  bina biquadrata, quorum summam in duo alia biquadrata resolvere liceat}, in:
  Op. Omnia (I) III, Teubner, 1917,  211--217.

\bibitem{eu62}
\leavevmode\vrule height 2pt depth -1.6pt width 23pt, {\em Specimen algorithmi
  singularis}, in: Op. Omnia (I) XV, Teubner, 1927,  31--49.

\bibitem{fe80}
M.~Feder, {\em On subspaces of spaces with an unconditional basis and spaces of
  operators}, Ill. J. Math. 24 (1980),  196--205.

\bibitem{fo64}
F.~Forelli, {\em The isometries of {$H^p$}}, Can. J. Math. 16 (1964),
  721--728.

\bibitem{fo82}
J.~J.~F. Fournier, {\em Two {UC}-sets whose union is not a {UC}-set}, Proc.
  Amer. Math. Soc. 84 (1982),  69--72.

\bibitem{gk95}
G.~Godefroy and N.~J. Kalton, {\em Approximating sequences and bidual
  projections}, Quart. J. Math. Oxford (2) 48 (1997),  179--202.

\bibitem{gkl96}
G.~Godefroy, N.~J. Kalton and D.~Li, {\em On subspaces of {$L^1$} which embed
  into $\ell_1$}, J. reine angew. Math. 471 (1996),  43--75.

\bibitem{gks93}
G.~Godefroy, N.~J. Kalton and P.~D. Saphar, {\em Unconditional ideals in
  {B}anach spaces}, Studia Math. 104 (1993),  13--59.

\bibitem{gs88}
G.~Godefroy and P.~D. Saphar, {\em Duality in spaces of operators and smooth
  norms on {B}anach spaces}, Ill. J. Math. 32 (1988),  672--695.

\bibitem{hr83}
H.~Halberstam and K.~F. Roth, {\em Sequences}, Springer, second~ed., 1983.

\bibitem{hww93}
P.~Harmand, D.~Werner and W.~Werner, {\em ${M}$-ideals in {B}anach spaces and
  {B}anach algebras}, Springer, 1993.

\bibitem{ha87}
S.~Hartman, {\em Some problems and remarks on relative multipliers}, Coll.
  Math. 54 (1987),  103--111.

\bibitem{hi82}
N.~Hindman, {\em On density, translates, and pairwise sums of integers}, J.
  Combin. Theory (A) 33 (1982),  147--157.

\bibitem{hmp86}
B.~Host, J.-F. M\'ela and F.~Parreau, {\em Analyse harmonique des mesures},
  Ast\'erisque 135--136, Soci\'et\'e math\'ematique de France, 1986.

\bibitem{jrz71}
W.~B. Johnson, H.~P. Rosenthal and M.~Zippin, {\em On bases, finite dimensional
  decompositions and weaker structures in {B}anach spaces}, Israel J. Math. 9
  (1971),  488--506.

\bibitem{ka96}
V.~M. Kadets, {\em Some remarks concerning the {D}augavet equation},
  Quaestiones Math. 19 (1996).

\bibitem{ka57}
J.-P. Kahane, {\em Sur les fonctions moyenne-p\'eriodiques born\'ees}, Ann.
  Inst. Fourier 7 (1957),  293--314.

\bibitem{ks63}
J.-P. Kahane and R.~Salem, {\em Ensembles parfaits et s\'eries
  trigonom\'etriques}, Hermann, 1963.
\newblock Actualit\'es Scientifiques et Industrielles 1301.

\bibitem{ka74}
N.~J. Kalton, {\em Spaces of compact operators}, Math. Ann. 208 (1974),
  267--278.

\bibitem{ka93}
\leavevmode\vrule height 2pt depth -1.6pt width 23pt, {\em ${M}$-ideals of
  compact operators}, Illinois J. Math. 37 (1993),  147--169.

\bibitem{kw95}
N.~J. Kalton and D.~Werner, {\em Property {$(M)$}, {$M$}-ideals, and almost
  isometric structure of {B}anach spaces}, J. reine angew. Math. 461 (1995),
  137--178.

\bibitem{ka48}
S.~Karlin, {\em Bases in {B}anach spaces}, Duke Math. J. 15 (1948),  971--985.

\bibitem{ka68}
Y.~Katznelson, {\em An introduction to harmonic analysis}, Wiley, 1968.

\bibitem{ke81}
T.~Ketonen, {\em On unconditionality in ${L}\sb{p}$ spaces}, Ann. Acad. Sci.
  Fenn. Ser. A I Math. Dissertationes 35 (1981).

\bibitem{ko91}
A.~L. Koldobsky, {\em Isometries of {$L_p(X;\,L_q)$} and equimeasurability},
  Indiana Math. J. 40 (1991),  677--705.

\bibitem{lps67}
L.~J. Lander, T.~R. Parkin and J.~L. Selfridge, {\em A survey of equal sums of
  like powers}, Math. Comp. 21 (1967),  446--459.

\bibitem{li95}
D.~Li, {\em On {H}ilbert sets and {$C_{\Lambda}(G)$}-spaces with no subspace
  isomorphic to $c_0$}, Coll. Math. 68 (1995),  67--77.
\newblock Addendum, ibid., p.\ 79.

\bibitem{li96}
\leavevmode\vrule height 2pt depth -1.6pt width 23pt, {\em Complex
  unconditional metric approximation property for {$C_{\Lambda}(\mathbf{T})$}
  spaces}, Studia Math. 121 (1996),  231--247.

\bibitem{lt77}
J.~Lindenstrauss and L.~Tzafriri, {\em Classical {B}anach spaces {I}. Sequence
  spaces}, Springer, 1977.

\bibitem{lp31}
J.~E. Littlewood and R.~E. A.~C. Paley, {\em Theorems on {F}ourier series and
  power series}, J. London Math. Soc. 6 (1931),  230--233.

\bibitem{lr75}
J.~M. Lopez and K.~A. Ross, {\em Sidon sets}, Dekker, 1975.

\bibitem{lu76}
F.~Lust-Piquard, {\em Ensembles de {R}osenthal et ensembles de {R}iesz}, C. R.
  Acad. Sci. Paris (A) 282 (1976),  833--835.

\bibitem{ma80}
B.~Maurey, {\em Isomorphismes entre espaces {$H_1$}}, Acta Math. 145 (1980),
  79--120.

\bibitem{me68}
Y.~Meyer, {\em Endomorphismes des id\'eaux ferm\'es de {$L^1(G)$}, classes de
  {H}ardy et s\'eries de {F}ourier lacunaires}, Ann. sci. \'Ecole Norm. Sup.
  (4) 1 (1968),  499--580.

\bibitem{me72}
\leavevmode\vrule height 2pt depth -1.6pt width 23pt, {\em Algebraic numbers
  and harmonic analysis}, North-Holland, 1972.

\bibitem{mi75}
I.~M. Mikheev, {\em On lacunary series}, Math. USSR--Sb. 27 (1975),  481--502.

\bibitem{mo39}
A.~Moessner, {\em Einige numerische {I}dentit\"aten}, Proc. Nat. Acad. Sci.
  India (A) 10 (1939),  296--306.

\bibitem{mo69}
L.~J. Mordell, {\em Diophantine equations}, Academic Press, 1969.

\bibitem{mu82}
T.~Murai, {\em On lacunary series}, Nagoya Math. J. 85 (1982),  87--154.

\bibitem{ne01}
S.~Neuwirth, {\em The {S}idon constant of sets with three elements}, electronic
  submission available at arXiv.org/abs/math.CA/0102145  (2001).

\bibitem{pe64}
A.~Pe{\l}czy\'nski, {\em On simultaneous extension of continuous functions. {A}
  ge\-ne\-ralization of theorems of {R}udin--{C}arleson and {B}ishop}, Studia
  Math. 24 (1964),  285--304.
\newblock Supplement, ibid.\ 25 (1965), 157--161.

\bibitem{pw71}
A.~Pe{\l}czy\'nski and P.~Wojtaszczyk, {\em Banach spaces with finite
  dimensional expansions of identity and universal bases of finite dimensional
  subspaces}, Studia Math. 40 (1971),  91--108.

\bibitem{pi81}
G.~Pisier, {\em De nouvelles caract\'erisations des ensembles de {S}idon}, in:
  Mathematical analysis and applications, part B, Advances in Math. Suppl.
  Series 7B, Academic Press, 1981,  685--726.

\bibitem{pi83}
\leavevmode\vrule height 2pt depth -1.6pt width 23pt, {\em Conditions
  d'entropie et caract\'erisations arithm\'etiques des ensembles de {S}idon},
  in: Topics in modern harmonic analysis II (Torino/Milano 1982), L.~{De
  Michele} and F.~Ricci (eds.), Ist. Naz. Alta Mat. Francesco Severi, 1983,
  911--944.

\bibitem{pl74}
A.~I. Plotkin, {\em Continuation of {$L^p$}-isometries}, J. Sov. Math. 2
  (1974),  143--165.

\bibitem{ra57}
S.~Ramanujan, {\em Notebooks}, Tata, 1957.

\bibitem{ra34}
S.~K. Rao, {\em On sums of sixth powers}, J. London Math. Soc. 9 (1934),
  172--173.

\bibitem{ri00}
E.~Ricard, {\em ${H}^1$ n'a pas de base compl\`etement inconditionnelle},
  electronic submission available at arXiv.org/abs/math.FA/0009073  (2000).

\bibitem{ro56}
M.~Rosenblatt, {\em A central limit theorem and a strong mixing condition},
  Proc. Nat. Acad. Sci. U. S. A. 42 (1956),  43--47.

\bibitem{ro79}
H.~P. Rosenthal, {\em Sous-espaces de {$L^1$}. Cours de troisi\`eme cycle},
  Universit\'e Paris 6, 1979.
\newblock Unpublished.

\bibitem{ru60}
W.~Rudin, {\em Trigonometric series with gaps}, J. Math. Mech. 9 (1960),
  203--228.

\bibitem{ru76}
\leavevmode\vrule height 2pt depth -1.6pt width 23pt, {\em {$L^p$-isometries
  and equimeasurability}}, Indiana Univ. Math. J. 25 (1976),  215--228.

\bibitem{ru78}
I.~Z. Ruzsa, {\em On difference sets}, Studia Sci. Math. Hung. 13 (1978),
  319--326.

\bibitem{se97}
J.~A. Seigner, {\em Rademacher variables in connection with complex scalars},
  Acta Math. Univ. Comenian. (N.S.) 66 (1997),  329--336.

\bibitem{si81}
I.~Singer, {\em Bases in {B}anach spaces {II}}, Springer, 1981.

\bibitem{sw95}
C.~M. Skinner and T.~D. Wooley, {\em On equal sums of two powers}, J. reine
  angew. Math. 462 (1995),  57--68.

\bibitem{st66a}
E.~Stein, {\em Classes ${H}\sp{p}$, multiplicateurs et fonctions de
  {L}ittlewood-{P}aley}, C. R. Acad. Sci. Paris S\'er. A-B 263 (1966),
  A716--A719.

\bibitem{st66b}
\leavevmode\vrule height 2pt depth -1.6pt width 23pt, {\em Classes ${H}\sp{p}$,
  multiplicateurs et fonctions de {L}ittlewood-{P}aley. {A}pplications de
  r\'esultats ant\'erieurs}, C. R. Acad. Sci. Paris S\'er. A-B 263 (1966),
  A780--A781.

\bibitem{st79}
C.~L. Stewart and R.~Tijdeman, {\em On infinite-difference sets}, Can. J. Math.
  31 (1979),  897--910.

\bibitem{ta88}
V.~Tardivel, {\em Ensembles de {R}iesz}, Trans. Amer. Math. Soc. 305 (1988),
  167--174.

\bibitem{we97}
D.~Werner, {\em The {D}augavet equation for operators on function spaces}, J.
  Funct. Anal. 143 (1997),  117--128.

\bibitem{wo84}
P.~Wojtaszczyk, {\em The {B}anach space ${H}\sb{1}$}, in: Functional analysis:
  surveys and recent results III (Paderborn, 1983), North-Holland, 1984,
  1--33.

\end{thebibliography}
\end{document}